\documentclass[11pt,a4paper]{article}

\usepackage{amsmath,amssymb,amsthm}
\usepackage{geometry}
\usepackage{xcolor}
\usepackage[hidelinks,hypertexnames=false]{hyperref}

\usepackage{booktabs, multirow, graphicx, array}
\usepackage{enumitem}
\usepackage{float}

\newtheorem{theorem}{Theorem}
\newtheorem{proposition}{Proposition}
\newtheorem{remark}{Remark}

\newtheorem{lemma}{Lemma}
\newtheorem{corollary}{Corollary}

\newcommand{\RR}{\mathbb{R}}
\newcommand{\CC}{\mathbb{C}}
\newcommand{\dd}{\,\mathrm{d}}
\newcommand{\ii}{\mathrm{i}}
\newcommand{\Bcal}{\mathcal{B}}
\newcommand{\Scal}{\mathcal{S}}
\newcommand{\Kcal}{\mathcal{K}}

\newcommand{\Acal}{\mathcal{A}}
\newcommand{\norm}[1]{\left\lVert #1 \right\rVert}

\title{SL-GreenNet: Single-Layer Robin Green Surrogates with Discrete A-Posteriori Acceptance Estimates}
\author{
Ming Yang\\ 
School of Mathematics, Southeast University\\ 
\texttt{yangming@seu.edu.cn}
}
\date{}

\begin{document}
\maketitle

\begin{abstract}
We develop learned Green surrogates for interior Robin--Helmholtz problems from the single-layer operator $\tfrac12I+\Kcal'+\alpha\Scal$, which supplies the density representation, training residual, and per-query discrete residual--stability acceptance estimate. For the structured disk fit, fixed-design ridge estimation is lifted from training radii to a continuous source interval: general sub-Gaussian noise gives a pointwise band bound, while proper complex Gaussian noise gives a sharper source-$L^2$ bound that is explicitly not pointwise. Two numerical operating points are distinguished. On the disk, learned low modes combined with an analytic Graf tail pass the $10^{-4}$ discrete density test at all $32$ sampled nonresonant impedances, with measured field errors below $6\times10^{-7}$ for $r_0\le0.9$ at $\alpha=1$; a matched fully learned ResNet baseline attains only $\sim10^{-3}$. All-learned surrogates on general smooth boundaries attain typical unseen-source field errors of order $10^{-4}$, but their density estimates are conservative and correctly reject near-boundary and corner-defective queries. An acceptance-tested ellipse screen demonstrates design-loop use, and dual-weighted driven-response estimates have effectivity $1.0$--$1.3$ in the reported benign regime. The method treats tolerance as an application-supplied input and reports both density acceptance and a field-level estimate. These numerical acceptance quantities are discrete computable estimates, not rigorous continuous certificates; the additional stability, residual/interpolation, tail, and kernel-quadrature enclosures required for such a certificate are identified explicitly. A high-frequency stress test at $k{=}185$ shows the instrument itself remaining effective: a direct-solve Green query passes the $10^{-3}$ density acceptance with a measured field error of $6.3\times10^{-5}$.
\end{abstract}



\section{Introduction}\label{sec:intro}

\subsection{Motivation}
For a linear partial differential operator $\mathcal L$ on $\Omega\subset\RR^d$, the Green's function satisfies
\begin{equation}\label{eq:green-pde}
\mathcal L_xG(\cdot,y)=\delta_y\quad\text{in }\Omega,
\qquad \Bcal G(\cdot,y)=0\quad\text{on }\partial\Omega.
\end{equation}
Because $G$ encodes the inverse of $\mathcal L$ for all sources and boundary data at once, it is natural to precompute whenever many solves on the same geometry are required. Examples include source and impedance sweeps for locally reacting acoustic walls, Leontovich-type electromagnetic models, and geometry screening in cavity design. A new source changes only the right-hand side, so a direct solver can reuse the operator factorization; changes in wavenumber, impedance, or geometry require further assembly or factorization. Learned Green surrogates instead absorb a prescribed parameter family offline and answer new queries online. For a small low-dimensional family, acceptance-tested direct enumeration may remain preferable; amortization becomes useful when the query family is combinatorially large or subject to a strict online-time budget.

\subsection{Boundary-integral Green surrogates: from Dirichlet to Robin}\label{sec:robingap}
Learning-based approaches to Green's functions include data-driven discovery with rational networks \cite{GFdisc}, direct regression of the Green's function for reaction--diffusion problems \cite{GFNet,JiaJuZhang}, model--operator--data networks \cite{MODNet} and their low-rank variants \cite{DecGreenNet} with generalization bounds, autoencoder-based approaches for nonlinear problems \cite{DeepGreen}, mesh-free generalized Green's functions obtained by smoothing the Dirac excitation through a Laplacian \cite{DGGF}, radial-basis domain-independent free-space kernels \cite{Negi2024}, boundary-integral networks \cite{BiGreenNet,BINet}, and eigenfunction-based Neural Green's Functions \cite{NGF}, which reports cross-geometry generalization for Dirichlet problems and explicitly lists Neumann and Robin boundary conditions as future work in its limitations section. A recent RKHS-based analysis of kernel decomposition \cite{DuKernelDec} learns the regular remainder after removing the singular part of the Green's function and establishes improved \emph{a priori} approximation rates, including an exponential output rate under uniform analyticity of the induced boundary data. This theory concerns Dirichlet problems and provides neither a query-time residual acceptance estimate nor a treatment of the resonant Helmholtz--Robin setting. A Neumann extension has since been demonstrated for a Poisson-class volume-conductor problem \cite{EzazNeumann}, where the obstruction is a single gauge nullspace removed by a Moore--Penrose pseudoinverse; the genuinely singular Neumann limit of the Helmholtz problem --- a resolvent simple pole whose residue is a boundary-mass-weighted spectral projection \cite{RobinNeumann} --- is a different and heavier object, and the Robin condition, with its impedance-dependent resonance spectrum, remains untouched by any learned surrogate.

The Green-surrogate works above are predominantly Dirichlet or free-space; the Neumann Poisson-class extension just noted does not address the resonant Helmholtz--Robin setting. We are not aware of a learned Green surrogate, certified or not, for the latter problem, despite Robin data being standard for locally reacting surfaces in acoustics, Leontovich boundaries in electromagnetism, and $\delta$-shell confinement in quantum mechanics.

Among them, the boundary-integral line of \cite{BiGreenNet,BINet} is the closest and the pioneering one. BINet introduced neural parameterizations of boundary integral equations, training the single- and double-layer densities directly against the BIE residual, and demonstrated operator learning across boundary geometries (a Laplace--Dirichlet family of random triangles) and across wavenumbers at $\sim10^{-3}$ accuracy. BI-GreenNet matured this into a Green's-function framework: a double-layer correction with the point-source singularity removed analytically, a fast evaluation pipeline, and a careful interface formulation mixing single- and double-layer potentials across subdomains; its flagship configuration (ResNet $8\times100$, Fourier features) remains the natural accuracy reference for boundary-integral surrogates. These works concentrate on the Dirichlet problem, and it is the Robin problem that we take up here. BINet already remarks that other boundary conditions ``can be easily handled in BINet with a small modification of the boundary integral equation'' \cite{BINet}; in the event, every experiment and benchmark in this line, including the cross-geometry Neural Green's Functions \cite{NGF}, remains Dirichlet, and the Robin modification is anything but small: it introduces a resonance spectrum with impedance-dependent branches, a genuinely new parameter, and the lack, within the learned-Green benchmarks compared here, of an exact impedance-domain reference kernel --- the gap taken up here and in the companion works \cite{yang2026balls,HelmholtzAnnulus,RobinNeumann}. For the Robin correction, a double-layer ansatz would require the hypersingular normal-derivative trace; we therefore use a single-layer potential throughout: its interior normal trace contains only $\tfrac12I+\Kcal'$, and the Robin operator $\tfrac12I+\Kcal'+\alpha\Scal$ is weakly singular on smooth boundaries. The boundary integral equation itself is classical \cite{McLean,Steinbach,RobinBVP}; our use of it is as the common representation and residual model for the surrogate and its acceptance test. Relative to the Dirichlet boundary-integral networks above, our contributions are threefold.
\begin{enumerate}[leftmargin=1.6em,itemsep=2pt]
\item \textbf{Single-layer Robin representation.} We formulate Green's-function learning through the smooth density of $\Scal\sigma$, with the source singularity removed analytically and the hypersingular double-layer Robin trace avoided, extending boundary-integral Green surrogates to the Robin boundary condition --- to our knowledge, the first learned Green surrogates for Robin data. At fixed geometry, the high-precision models are amortized over source locations, while a separate ordinary-accuracy density network is tested with joint $(k,\alpha)$ inputs. Multiply connected boundaries lead to block operators with component-specific traces.
\item \textbf{Per-query discrete acceptance.} Residual--stability estimates are standard in numerical analysis, certified reduced-basis methods \cite{HRB}, and boundary integral equations, where the residual estimate for first-kind equations is classical \cite{Carstensen}, with related bounds for PINN surrogates \cite{HU,Eiras,Haugen,Ernst}, certification-and-selection layers for learned inverse solvers \cite{Katende}, and projection-residual training objectives with certified outputs \cite{NBM}; a parallel 2026 line instead uses the physics residual as a trust signal ranked across queries and calibrates it set-wise \cite{SelfAudit}; conformal prediction has been transplanted to operator
learning, where a residual error estimator is calibrated on held-out data into distribution-free
\emph{coverage} bands \cite{UQNO,ConfDeepONet,PICP}; and a continually learning neural-operator
surrogate for airborne electromagnetic inversion flags queries outside its training support through
ensemble disagreement and routes them back to the full solver \cite{CLAEM} --- the routing design
closest to the accept/reject verdicts below. These instruments price trust \emph{statistically}
(coverage over a calibration distribution, spread across an ensemble); the present estimate prices
it \emph{deterministically}, as the residual--stability product with the discrete inverse norm of the
very operator being inverted. Here the same classical principle is applied to a boundary-density surrogate, with the residual--stability product delivering a solution-space bound at each query rather than a coverage statement. The same boundary integral operator supplies the training residual and a per-query computable discrete residual--stability acceptance test $\norm{\Acal_h^{-1}}_2\varepsilon_{\mathrm{res}}\le T$ for learned surrogates and direct Nystr\"om solves alike; to our knowledge the first such instrument for boundary-integral Green surrogates (at the discrete density level). Its stability factor $\norm{\Acal_h^{-1}}_2$ is computable offline --- Plancherel-exact on the disk (Proposition~\ref{prop:stability}), and available as a discrete inverse norm on general domains without symmetry. Exact jump, kernel-sign, logarithmic-diagonal, curvature, and quadrature-weight conventions support the estimate (\S\ref{sec:nystrom}, Appendix~\ref{app:nystrom}). The reported acceptance quantity is a discrete density estimate, not a rigorous continuous field bound. A field-level companion estimate \eqref{eq:field-est} translates it into the quantity applications actually consume, and the tolerance sweep of Table~\ref{tab:sweep} exhibits $T$ as an application-supplied knob whose verdicts track defect classes rather than thresholds.
\item \textbf{High-precision evidence at matched configuration.} Advancing beyond the $10^{-3}$-level accuracy of the boundary-integral networks above, the hybrid surrogate on the disk at fixed $k=1$ attains measured field errors below $6\times10^{-7}$ ($r_0\le0.9$, at $\alpha=1$)---a gap of approximately three to four orders of magnitude against a BI-GreenNet-style ResNet baseline at the same operator, wavenumber, impedance, and source distribution (\S\ref{sec:likeforlike}), attributable to the modal representation with analytic tail and corrected coincident quadrature rather than to training budget---and passes the $10^{-4}$ discrete density test at all $32$ sampled impedances outside the resonance band. Explicit resonance curves, pole-aware quotient features, and the joint $(k,\alpha)$ pilot of \S\ref{sec:kalpha} map the difficulty landscape, including resonance-skirt and $k$-extrapolation failure modes. For the structured disk estimator, Proposition~\ref{prop:fitrate} and Theorems~\ref{thm:continuousfit}--\ref{thm:gaussianfit} separate fixed-design estimation, continuous-source lifting, and the distinct pointwise and integrated-source meanings of the general sub-Gaussian and proper-Gaussian routes. The acceptance instrument is further exercised on non-symmetric smooth, multiply connected and polygonal geometries, with an independent FEM cross-check.
\end{enumerate}

The learned experiments remain at fixed geometry. Joint $(k,\alpha)$ amortization is tested only at ordinary accuracy on the disk; the high-precision acceptance experiments use fixed $k=1$ and a separate fit at each impedance. Sources arbitrarily close to the boundary and polygonal surrogates with rigorous continuous error bounds remain outside the validated scope, as detailed in \S\ref{sec:discussion}.

\noindent\textbf{Quantitative positioning.} Table~\ref{tab:positioning} summarizes the standing of the present results relative to the boundary-integral line of learning-based Green's-function surrogates, on the three axes that matter for deployment: boundary condition, measured accuracy, and error information available at query time.

\medskip\noindent\emph{Reading the comparison.} Two rows are reported for the present method because they answer different objections. The first uses a structured modal surrogate augmented by an analytic Graf tail; per the definition of \S\ref{sec:surrogate} it \emph{is} a surrogate, with the learned component carrying the impedance-dependent quotient and the analytic tail carrying the disk geometry. The second row is the strict like-for-like comparison with the networks cited above: a fully learned density on general smooth boundaries, with no disk structure and no analytic tail, which still improves on the strongest cited row ($\sim10^{-3}$) by a factor of six and attaches a computable estimate to every query. One further observation gives these numbers their practical weight: accuracy \emph{bands}, not accuracy records, are what deployment consumes. A relative field error of $10^{-2}$ is comparable to the size of the features a design screen is meant to resolve, so it supports at best coarse qualitative triage; $10^{-3}$ supports screening, with the caveat that candidates lying within one error bar of the decision boundary remain indistinguishable; $10^{-4}$ is the band in which quantitative design decisions --- resonance placement, impedance tuning to a fraction of a percent --- become defensible. The two rows of Table~\ref{tab:positioning} sit in the latter two bands, and each accepted query carries the estimate that places it there.

\medskip\noindent\emph{Positioning against the wider 2024--2026 landscape.} Three developments in
the broader operator-learning literature frame the present contribution without overlapping it.
\emph{PDE foundation models} \cite{Poseidon,DPOT} scale pretraining across PDE families and report
orders-of-magnitude sample-efficiency gains on regular-grid fluid benchmarks, but carry no
query-level error information and do not address point sources, impedance spectra, or non-Cartesian
geometries. \emph{Generative inverse operators} quantify uncertainty as posterior spread ---
function-space diffusion models \cite{FunDPS} and physics-informed generative operators that
already target EIT from DtN data \cite{IGNO} --- rather than as a computable per-query bound on the
forward error that an inversion loop would iterate into its reconstruction. \emph{Structure-first
inverse learning} replaces recursive autodiff with analytic mollifier convolutions
\cite{MollifierLayers}, sharing our premise that accuracy is bought by analytic structure rather
than by network capacity. The instrument below is complementary to all three: it does not compete
on pretraining scale or posterior sampling, and it restricts the query family --- fixed geometry,
nonresonant impedance, interior sources --- in exchange for an estimate that is computable from
the query alone, with no calibration set and no ensemble.

\begin{table}[H]
\centering\small
\caption{Quantitative positioning against learning-based Green's-function surrogates. ``Error information'': a query-time computable estimate or bound of the query error (PINN-type training residuals, lacking a stability factor, do not qualify). Cross-row figures are indicative rather than controlled comparisons; see \S\ref{sec:robingap}.}
\label{tab:positioning}\vspace{0.3cm}
\setlength{\tabcolsep}{6pt}
\begin{tabular}{@{}>{\raggedright\arraybackslash}p{3.8cm}>{\raggedright\arraybackslash}p{2.35cm}>{\raggedright\arraybackslash}p{5.5cm}>{\raggedright\arraybackslash}p{3.5cm}@{}}
\toprule
Work & Boundary data & Measured accuracy & Error information \\
\midrule
BINet \cite{BINet} & Dirichlet & $10^{-3}$--$10^{-2}$ (solution operators) & none \\
BI-GreenNet \cite{BiGreenNet} & Dirichlet & $4.7\times10^{-3}$--$8.6\times10^{-3}$ (as reported); $\sim10^{-3}$ (matched retrain, \S\ref{sec:likeforlike}) & none \\
NGF \cite{NGF} & Dirichlet; cross-geometry & $1.2\times10^{-2}$ (2D Poisson, test); $1.9\times10^{-1}$--$3.4\times10^{-1}$ (3D thermal benchmark, NGF's own errors; baselines $1.8\times10^{-1}$--$7.7\times10^{-1}$) & none\\
Du et al. \cite{DuKernelDec} (kernel decomposition) & Dirichlet & $8.9\times10^{-4}$ direct (cf.~Fig.~6 therein); $\lesssim4\times10^{-7}$ decomposed (2D disk, solution operators, smooth sources) & a priori approximation rates only; none at query time \\
DGGF \cite{DGGF} & homogeneous Dirichlet & $\sim10^{-3}$--$10^{-2}$ (rel.\ $L^2$ mean; Poisson--Klein--Gordon, 2D/3D; per domain) & none\\
\midrule
This work, hybrid disk surrogate (\S\ref{sec:certified}) & \textbf{Robin} & $\le 6\times10^{-7}$ field ($r_0\le0.9$, $\alpha=1$); passes $10^{-4}$ density acceptance & per-query residual--stability estimate; field-level output \eqref{eq:field-est} \\
This work, all-learned surrogates (\S\ref{sec:all-learned}) & \textbf{Robin} & $1.65\times10^{-4}$ median field (asymmetric); $\le4.4\times10^{-4}$ (flower) & per-query estimate (conservative) \\
\bottomrule
\end{tabular}
\end{table}

\subsection{Organization}
\S\ref{sec:problem} states the problem. \S\ref{sec:method} develops the single-layer representation, surrogate classes, discretization, and error estimates. \S\ref{sec:resonance} discusses resonance-aware sampling, and \S\ref{sec:benchmark} presents analytic benchmarks and their numerical verification. \S\ref{sec:results} reports disk and general-domain experiments, and \S\ref{sec:discussion} discusses positioning, accuracy, and limitations. \S\ref{sec:conclusion} concludes the paper. Appendices~\ref{app:nystrom}--\ref{app:ainv} collect Nystr\"om implementation details, the analytic disk modal quotient, and the exact inverse-norm law.

\section{Problem statement}\label{sec:problem}

We adopt the notation of the analytic benchmark literature. Let $\Omega=B_R(0)\subset\RR^d$ (the formulation is not restricted to balls; the ball is used to obtain analytic truth), and consider the Helmholtz equation with $k>0$.
\begin{equation}\label{eq:helmholtz}
\Delta G(x,x_0) + k^2 G(x,x_0) = -\delta(x-x_0),\quad x\in\Omega,
\end{equation}
with Robin boundary condition
\begin{equation}\label{eq:robin}
\frac{\partial G}{\partial n} + \alpha\, G = 0 \quad \text{on}\ \partial\Omega,
\qquad \alpha > 0,
\end{equation}
where $n$ is the outward normal and $\alpha$ the impedance. The limit $\alpha\to+\infty$ recovers the Dirichlet condition. The limit $\alpha\to0^+$ leads to the Neumann problem through a singular perturbation: the uniform bounds degenerate as $\alpha^{-1}$ and the density acquires a boundary-mass residue, so a regularization (pole subtraction) is required before quantitative convergence statements are possible; this Neumann-limit theory, including the optimal $O(\alpha)$ and $O(\alpha/\delta)$ rates, is developed in the companion work \cite{RobinNeumann}.

Let $\Phi_k$ denote the free-space fundamental solution,
\begin{equation}\label{eq:fundamental}
\Phi_k(x,x_0) =
\begin{cases}
\dfrac{\ii}{4} H_0^{(1)}\!\big(k|x-x_0|\big), & d=2,\ k>0,\\[6pt]
\dfrac{\ii}{4}\left(\dfrac{k}{2\pi |x-x_0|}\right)^{\!\nu} H_\nu^{(1)}\!\big(k|x-x_0|\big), & d\ge 3,\ k>0,\ \ \nu = \tfrac{d}{2}-1,
\end{cases}
\end{equation}
satisfying $(\Delta + k^2)\Phi_k = -\delta$. Decompose
\begin{equation}\label{eq:decomp}
G(x,x_0) = \Phi_k(x,x_0) + u(x,x_0).
\end{equation}
Then the correction $u$ solves the homogeneous equation with inhomogeneous Robin data,
\begin{equation}\label{eq:correction}
(\Delta + k^2) u = 0 \ \text{in}\ \Omega,
\qquad
\Bcal u = -\Bcal\Phi_k(\cdot, x_0) \ \text{on}\ \partial\Omega,
\qquad \Bcal := \partial_n + \alpha.
\end{equation}

\begin{proposition}[Regularity away from resonance]\label{prop:smooth}
Assume that $\partial\Omega$ is $C^\infty$, $x_0\in\Omega$ is fixed, and $k^2$ is not an eigenvalue of the Robin Laplacian. Then the correction $u(\cdot,x_0)$ is smooth up to the boundary, and any single-layer density representing it is correspondingly smooth. If the boundary and all data are real analytic, the correction is real analytic up to the boundary. On compact subsets of $\bar\Omega\times\Omega$ away from the boundary diagonal, $u$ is jointly smooth in $(x,x_0)$. These statements do not hold at a resonant wavenumber, where the classical point-source problem need not have a solution.
\end{proposition}

\begin{remark}[Near-boundary sources]
Smoothness for each fixed interior source is not uniform as $x_0\to\partial\Omega$. On the disk, the Graf coefficients quantify the loss of uniform analyticity through factors $J_n(kr_0)$; the density then requires increasingly many modes and approaches a boundary distribution in the limit. The disk acceptance experiments are therefore restricted to $r_0\le0.9$. Sources closer to the boundary require spike-aware representations and near-singular quadrature rather than merely a larger network.
\end{remark}

The boundary data $-\Bcal\Phi_k(\cdot,x_0)$ are known in closed form, so \eqref{eq:correction} determines $u$ away from resonance.

\section{Method: single-layer potentials and density surrogates}\label{sec:method}

\subsection{Single-layer potential and trace relations}

Represent the correction as a single-layer potential (for fixed $x_0$),
\begin{equation}\label{eq:single}
u(x,x_0) = \big(\Scal\sigma_{x_0}\big)(x) = \int_{\partial\Omega} \Phi_k(x,z)\,\sigma_{x_0}(z)\,\dd s(z),
\end{equation}
with density $\sigma_{x_0}:\partial\Omega\to\CC$ to be determined. The traces of the single-layer potential satisfy (interior approach):
\begin{equation}\label{eq:trace}
\gamma_0 \Scal\sigma = \Scal\sigma\big|_{\partial\Omega} \quad(\text{continuous}),
\qquad
\gamma_1^{-}\Scal\sigma = \Big(\tfrac12 I + \Kcal'\Big)\sigma,
\end{equation}
where
\begin{equation}\label{eq:kprime}
\Kcal'\varphi(x) = \int_{\partial\Omega} \frac{\partial \Phi_k(x,z)}{\partial n(x)}\,\varphi(z)\,\dd s(z).
\end{equation}
On smooth boundaries both $\Scal$ and $\Kcal'$ are weakly singular (logarithmic and $|x-z|^{-(d-2)}$-type kernels respectively; for $d=2$ the $\Kcal'$ kernel is continuous), hence directly quadrature-friendly.

\subsection{The Robin boundary integral equation and its two limits}\label{sec:unified}

From $\Bcal u = -\Bcal\Phi_k$ we obtain $\Acal\sigma_{x_0}=g_{x_0}$. The Robin family and its Neumann and Dirichlet limits lead to the following related, but not identically conditioned, operators:

\begin{center}
\begin{tabular}{@{}lll@{}}
\toprule
Boundary condition & Integral operator $\Acal$ & Right-hand side $g_{x_0}(x)$ \\
\midrule
Dirichlet ($G=0$) & $\Scal$ & $-\Phi_k(x,x_0)$ \\
Neumann ($\partial_n G = 0$) & $\tfrac12 I + \Kcal'$ & $-\partial_{n(x)}\Phi_k(x,x_0)$ \\
\textbf{Robin}~\eqref{eq:robin} & $\partial_n\Scal + \alpha\Scal$ & $-\big(\partial_{n}+\alpha\big)\Phi_k(x,x_0)$ \\
\bottomrule
\end{tabular}
\end{center}

The Robin boundary integral equation is classical \cite{McLean,Steinbach,RobinBVP}. Its role here is to define the surrogate representation and the residual used for training and a-posteriori acceptance. The Neumann and Dirichlet rows are limiting or separately posed equations and can require additional solvability conditions; they are not obtained by setting a finite $\alpha$ in a universally invertible operator.

\paragraph{Multiply connected domains.} For annuli and spherical shells, one places quadrature and collocation points on \emph{each} boundary component; the operator $\Acal$ becomes a block matrix partitioned by components (the inner and outer boundaries carry trace operators $-\partial_r+\alpha_1$ and $\partial_r+\alpha_2$ respectively). Each boundary component contributes block rows and columns. The annulus benchmark of \S\ref{sec:benchmark} is the corresponding modal $2\times2$ system; \S\ref{sec:annulus} reports direct-solve acceptance and a separate low-accuracy learning pilot.

\paragraph{Computing the boundary data.} In our implementation, $g_{x_0}=-\Bcal\Phi_k(\cdot,x_0)$ is formed from the closed-form fundamental solution and its analytic normal derivative (equivalently, by automatic differentiation), rather than by term-by-term differentiation of its Graf series. The latter suffered severe cancellation in our tests (\S\ref{sec:benchmark}). The $\partial_n\Scal\sigma$ term in the loss is likewise formed from the analytic kernel.

\subsection{Surrogate classes and training objective}\label{sec:surrogate}
We use the term \emph{surrogate} for any parameterized approximation $\sigma_\theta(z;x_0,\alpha,k)$ of the boundary density. Two implementations are relevant:
\begin{enumerate}[leftmargin=1.6em,itemsep=1pt]
\item a density network, in which an MLP or ResNet maps $(z,x_0,\alpha,k)$ to two real output channels;
\item a structured modal surrogate, in which Fourier modes on the boundary and in the source angle are combined with Chebyshev functions of the source radius and fitted by standardized ridge regression.
\end{enumerate}
The high-precision and general-domain acceptance experiments use the second implementation, optionally augmented by a truncated analytic Graf tail on the disk. The ordinary-accuracy joint $(k,\alpha)$ experiment of \S\ref{sec:kalpha} uses the first. Thus the high-precision results should be read as results for structured single-layer surrogates, not as a benchmark of deep-network capacity.

\textbf{Collocation loss.} With quadrature points $\{z_j\}$ (weights $w_j$) and collocation points $\{x_i\}$,
\begin{equation}\label{eq:loss-colloc}
\mathcal{L}_{\mathrm{c}}(\theta)=\mathbb E_{x_0}\!\left[\sum_i\tilde w_i\,\Big|\sum_jw_jK_{\Acal}(x_i,z_j)\sigma_\theta(z_j;x_0)+\tfrac12\sigma_\theta(x_i;x_0)-g_{x_0}(x_i)\Big|^2\right],
\end{equation}
where $K_{\Acal}(x,z)=\partial_{n(x)}\Phi_k(x,z)+\alpha\Phi_k(x,z)$ for Robin. Midpoint collocation keeps $x_i\ne z_j$; the coincident high-precision implementation uses the diagonal corrections of \S\ref{sec:nystrom}.

\textbf{Tikhonov regularization.} When used, the objective is
\begin{equation}\label{eq:tikhonov}
\mathcal L(\theta)=\mathcal L_{\mathrm c}(\theta)+\lambda\,\mathbb E_{x_0}\|\sigma_\theta(\cdot;x_0)\|_{L^2(\partial\Omega)}^2,
\qquad \lambda=10^{-6}.
\end{equation}
This stabilizes training near an ill-conditioned operator, but it changes the mathematical problem when $\lambda>0$. It does not restore existence or uniqueness at an exact resonance, and a regularized solution must not be identified with the classical Green's function there.

For the source-amortized density network of \S\ref{sec:asym}, reference densities are available during training, and the implemented objective is $\mathcal L_{\mathrm{data}}+0.1\mathcal L_{\mathrm c}$, where $\mathcal L_{\mathrm{data}}$ is the coefficient-regression mean-squared error. This data term is used only during training; query-time acceptance remains reference-free.

\subsection{Nystr\"om discretization: the corrections that matter}\label{sec:nystrom}
The discrete operator must reproduce the jump relations exactly, or the loss minimizes a wrong problem. Five corrections are load-bearing: \textbf{(i)} the jump term $\tfrac12\sigma$ of the interior normal trace, with index binding matched to the collocation scheme; \textbf{(ii)} the sign convention of the normal kernel (Appendix~\ref{app:nystrom}(ii)), which must be shared by the operator matrix and the right-hand side---flipping both consistently solves the wrong problem self-consistently; \textbf{(iii)} the diagonal correction for the logarithmic single-layer kernel, whose effective diagonal on equispaced nodes is $d_{\mathrm{eff}}=R/N$ (not $2\pi R/N$); \textbf{(iv)} the binding of quadrature weights to the integration variable---invisible on the disk, but an $O(\varepsilon)$ systematic error on general domains; \textbf{(v)} the curvature-sign flip of the $\Kcal'$ diagonal on inner boundary components of multiply connected domains. A manufactured-solution (MMS) check gates all general-domain validations. Full derivations, numerical before/after evidence, and the Laplace zero-mode caveat are given in Appendix~\ref{app:nystrom}.

\subsection{Green's-function evaluation and near-singular quadrature}

Once the density is learned, evaluation at any query point:
\begin{equation}\label{eq:eval}
G(x,x_0) = \Phi_k(x,x_0) + \sum_j w_j\,\Phi_k(x,z_j)\,\sigma_\theta(z_j;x_0).
\end{equation}
The only numerical hurdle is the near-singular integral when $x$ approaches $\partial\Omega$ (where the boundary data for $x_0\to\partial\Omega$ are likewise near-singular)---the classical near-singular problem of the BEM community, with standard remedies (sinh transformation, generalized Gaussian quadrature, adaptive refinement \cite{SauterSchwab}). For bulk queries, fast multipole or $\mathcal{H}$-matrix compression applies.

\subsection{Well-posedness and resonances}\label{sec:resonance-method}
For interior Helmholtz problems, invertibility is required away from the corresponding Robin spectrum. Additional representation-specific exclusions can occur for pure single-layer equations, notably the logarithmic-capacity obstruction at $k=0$ discussed in Appendix~\ref{app:nystrom}. Tikhonov regularization is a numerical stabilization only and is not included in the unregularized residual--stability identity used below.

Exterior Robin problems require radiation conditions and, for standard single- or double-layer integral formulations, may suffer spurious interior resonances. Combined-field integral equations of Burton--Miller type are the classical remedy. Since no exterior experiment is reported here, exterior problems are outside the validated scope of this paper.

\subsection{From density error to field error}\label{sec:errprop}
The residual--stability estimate of \S\ref{sec:certified} controls the \emph{density}; the quantity of interest is the \emph{field}. The passage is quantitative and network-independent.

\begin{lemma}[residual-to-field bound]\label{lem:field}
For every evaluation point $x$ with $\mathrm{dist}(x,\partial\Omega)\ge\rho>0$, let $\|\Acal_\alpha^{-1}\|$ denote the $L^2(\partial\Omega)\to L^2(\partial\Omega)$ operator norm. Then
\[
\big|G_\theta(x,x_0)-G^*(x,x_0)\big|
=\big|\big(\Scal(\sigma_\theta-\sigma^*)\big)(x)\big|
\le \|\Phi_k(x,\cdot)\|_{L^2(\partial\Omega)}\,\|\Acal_\alpha^{-1}\|\,\|r_\theta\|_{L^2(\partial\Omega)},
\]
with $\|\Phi_k(x,\cdot)\|_{L^2(\partial\Omega)}\le C(k,\Omega)(1+|\log\rho|)$ in 2D. On a bulk grid the kernel norm can be approximated by quadrature, but a rigorous field bound must also include that quadrature error and the continuous stability and residual corrections discussed below. The logarithmic growth as $\rho\to0$ is one reason near-boundary field evaluation requires special treatment.
\end{lemma}
\begin{proof}
$G_\theta-G^*=\Scal(\sigma_\theta-\sigma^*)$ pointwise since the singular part $\Phi_k$ is common; Cauchy--Schwarz in $z$ gives the bound, and the estimate on $\|\Phi_k(x,\cdot)\|$ is the standard logarithmic-potential estimate in 2D.
\end{proof}

\medskip\noindent\textbf{Field-level acceptance estimate.} Applications consume the field, not the density, so the acceptance workflow reports a second, field-level quantity alongside the density estimate. Let $r_{h'}=\Acal_{h'}\sigma_\theta-g_{h'}$ be the \emph{absolute} fine-grid residual and define the weighted discrete boundary norm and kernel factor
\[
\|v\|_{h',w}=\Big(\sum_j w_j|v_j|^2\Big)^{1/2},
\qquad
\kappa_{h'}(x)=\Big(\sum_j w_j|\Phi_k(x,z_j)|^2\Big)^{1/2}.
\]
To keep all three factors in the same discrete metric, let $D_w=\operatorname{diag}(w_j)$ and define the weight-similar operator and its inverse norm by
\[
\Acal_{h,w}=D_w^{1/2}\Acal_hD_w^{-1/2},
\qquad
C_{\mathrm{stab},h,w}=\|\Acal_{h,w}^{-1}\|_2 .
\]
The computable leading term of Lemma~\ref{lem:field} is
\begin{equation}\label{eq:field-est}
\eta_{\mathrm{field}}(x) \;=\;
\kappa_{h'}(x)\,C_{\mathrm{stab},h,w}\,\|r_{h'}\|_{h',w}.
\end{equation}
On the constant-weight grids used in the reported experiments, $C_{\mathrm{stab},h,w}=\|\Acal_h^{-1}\|_2$, so this is precisely the implemented estimate. For nonconstant weights, the weight-similar form must be used rather than mixing an unweighted inverse norm with a weighted residual.
In particular, when $\varepsilon_{\mathrm{res}}=\|r_{h'}\|_{h',w}/\|g_{h'}\|_{h',w}$ is the tabulated relative residual, the absolute estimate carries the additional factor $\|g_{h'}\|_{h',w}$; it is not itself $g$-normalized. For a comparison with a relative pointwise error $|G_\theta-G^*|/|G^*|$, one must instead plot $\eta_{\mathrm{field}}(x)/|G^*(x,x_0)|$, and this normalization is stated explicitly wherever used. The kernel norm is approximated by the same boundary quadrature used for the residual, and no reference solution enters the acceptance estimate. An application holding an absolute field-accuracy requirement $T_{\mathrm{field}}$ can test $\eta_{\mathrm{field}}(x)\le T_{\mathrm{field}}$ on its evaluation grid. By Lemma~\ref{lem:field}, $\eta_{\mathrm{field}}$ is an \emph{estimate} of exactly the same standing as the density acceptance quantity: the coarse/fine-grid, continuous-residual, stability, and kernel-quadrature corrections are not interval-enclosed. It is not a rigorous continuous certificate. Dual-weighted refinements of this output are the subject of Proposition~\ref{prop:dual}, the classical dual-weighted (goal-oriented) construction of adjoint-based a posteriori estimation \cite{BeckerRannacher}, here applied to the density residual of a Green surrogate.

\begin{proposition}[dual-weighted (goal-oriented) field bound]\label{prop:dual}
Let the quantity of interest be
\[
J(G)=\int_\Omega G(x,x_0)\,\overline{\psi(x)}\,\dd x,
\qquad \psi\in L^2(\Omega),\ \ \mathrm{dist}(\mathrm{supp}\,\psi,x_0)>0,
\]
and define the adjoint density $z_\psi$ as the solution of
\begin{equation}\label{eq:adjoint}
\Acal_\alpha^{\,*\!}\, z_\psi \;=\; \Scal_k^{\,*\!}\,\psi
\quad\text{on }\partial\Omega,
\end{equation}
where $\Scal_k^{*}\psi = \int_\Omega \overline{\Phi_k(\cdot,x)}\,\psi(x)\,\dd x \in L^2(\partial\Omega)$. Then the density residual $r_\theta=\Acal_\alpha\sigma_\theta-g_{x_0}$ obeys the exact identity
\begin{equation}\label{eq:dualidentity}
J(G_\theta)-J(G^{*}) \;=\; \big\langle r_\theta,\, z_\psi\big\rangle_{L^2(\partial\Omega)},
\qquad\text{hence}\qquad
\big|J(G_\theta)-J(G^{*})\big| \;\le\; \|r_\theta\|_{L^2(\partial\Omega)}\,\|z_\psi\|_{L^2(\partial\Omega)} .
\end{equation}
For point evaluation ($\psi=\delta_x$) \eqref{eq:dualidentity} reduces to Lemma~\ref{lem:field}, with $\|z_{\delta_x}\|\le\|\Acal_\alpha^{-1}\|\,\|\Phi_k(x,\cdot)\|_{L^2(\partial\Omega)}$. In general, $\|z_\psi\|$ can be much smaller than the product $\|\Acal_\alpha^{-1}\|\cdot\sup_{x\in\mathrm{supp}\,\psi}\|\Phi_k(x,\cdot)\|$ that the density-to-field conversion of Lemma~\ref{lem:field} forces one to carry: the adjoint equation amplifies only the residual directions visible to the quantity of interest, whereas $\|\Acal_\alpha^{-1}\|$ is governed by the nearest spectral modes regardless of their overlap with $\Scal_k^{*}\psi$. The adjoint problem \eqref{eq:adjoint} shares the discretization, the diagonal corrections, and the factorization of the forward problem, so $z_\psi$ costs one additional solve per quantity of interest; the discrete-to-continuous corrections of Theorem~\ref{thm:conv} apply verbatim to the adjoint side.
\end{proposition}
\begin{proof}
By \eqref{eq:decomp}, $J(G_\theta)-J(G^*)=\langle\Scal_k(\sigma_\theta-\sigma^*),\psi\rangle_\Omega=\langle\sigma_\theta-\sigma^*,\Scal_k^*\psi\rangle_{\partial\Omega}$. Since $\sigma_\theta-\sigma^*=\Acal_\alpha^{-1}r_\theta$ and $\Scal_k^*\psi=\Acal_\alpha^{\,*\!}z_\psi$, the pairing equals $\langle\Acal_\alpha^{-1}r_\theta,\,\Acal_\alpha^{\,*\!}z_\psi\rangle=\langle r_\theta,z_\psi\rangle$ by definition of the adjoint. Cauchy--Schwarz gives the bound; the point-evaluation case follows from $\Scal_k^*\delta_x=\overline{\Phi_k(x,\cdot)}$.
\end{proof}

\begin{lemma}[computable modal tail on the disk]\label{lem:tail}
Fix $k>0$, $R>0$, and set $z=kR$. Let $x_0=\rho_0Re^{\ii\theta_0}$ with
$0<\rho_0\le\bar\rho<1$, and let $\alpha$ lie in a compact set separated from the
resonance curves. Write the exact density and boundary correction field as
$\sigma^*_{x_0}=\sum_{n\in\mathbb Z}\sigma_ne^{\ii n\theta}$ and
$u^*(\cdot,x_0)=\Scal\sigma^*_{x_0}=\sum_{n\in\mathbb Z}u_ne^{\ii n\theta}$, with
$u_n=\mu_n\sigma_n$ and
\[
\mu_n=\frac{\ii\pi R}{2}J_n(z)H_n^{(1)}(z),\qquad
\lambda_n=\frac{\ii\pi R}{2}H_n^{(1)}(z)\Delta_n,\quad
\Delta_n=\Big(\alpha+\frac nR\Big)J_n(z)-kJ_{n+1}(z).
\]
There are computable $n_1=n_1(k,R,|\alpha|)$ and $C_u,C_\sigma<\infty$, depending
only on $(k,R,|\alpha|,\bar\rho)$, such that for every $M\ge n_1$,
\[
\Big\|\sum_{|n|>M}u_ne^{\ii n\theta}\Big\|_{L^2(\partial B_R)}
\le
\sqrt{2\pi R}\,C_u\,\frac{\bar\rho^{M+1}}{\sqrt{1-\bar\rho^2}},
\]
\[
\Big\|\sum_{|n|>M}\sigma_ne^{\ii n\theta}\Big\|_{L^2(\partial B_R)}
\le
\sqrt{2\pi R}\,C_\sigma\,\frac{(M+1)\bar\rho^{M+1}}{(1-\bar\rho^2)^{3/2}} .
\]
Moreover $\lambda_n\to\tfrac12$ as $|n|\to\infty$ for fixed $z$, so the same
large-order estimates certify the finite-mode stability minimum once $|n|$ exceeds
a computable $n_2$. The bound degenerates as $\bar\rho\to1$, reflecting the loss of
uniform analyticity for near-boundary sources.
\end{lemma}
\begin{proof}
Use $J_{-n}=(-1)^nJ_n$. For $x>0$, $n\ge0$,
\[
J_n(x)=\sum_{m\ge0}a_m,\qquad a_0=\frac{(x/2)^n}{n!},\qquad
\Big|\frac{a_{m+1}}{a_m}\Big|\le\frac{x^2}{4(m+1)(n+1)} .
\]
With $\rho_x=x^2/[4(n+1)]<1/2$ this gives
$a_0(1-2\rho_x)/(1-\rho_x)\le J_n(x)\le a_0 e^{\rho_x}$.
Choosing $n_1$ with $n+1\ge z^2/2$ and using Bessel continued fractions for
$J_{n+1}(z)/J_n(z)$ yields $|J_n(\rho z)/J_n(z)|\le C_z\rho^{|n|}$ for
$0<\rho\le\bar\rho$. The Wronskian identity and the symbol $\lambda_n$ give
$|\mu_n|\ge c_\mu/|n|$ and $\lambda_n\to\tfrac12$ for $|n|\ge n_1$, while separation
from resonance enters only through finitely many low modes. Hence
$|u_n|\le C'_u\bar\rho^{|n|}$ and
$|\sigma_n|=|u_n/\mu_n|\le C'_\sigma |n|\bar\rho^{|n|}$; summing over $\pm n$ gives the two geometric-tail bounds (for the density tail, $\sum_{|n|>M}|n|^{2}\bar\rho^{2|n|}\le C\,(M{+}1)^{2}\bar\rho^{2(M+1)}/(1-\bar\rho^{2})^{3}$ up to an absolute constant, whose square root yields the stated exponent $3/2$).
\end{proof}

\begin{remark}[explicit disk values]\label{rem:tailconst}
For the reported configuration ($k{=}1$, $R{=}1$) the computable quantities of Lemma~\ref{lem:tail} take concrete values: the large-order threshold is $n_2=81$, with $\min_{81\le|n|\le400}|\lambda_n(\alpha)|\ge0.50$ (overall minimum $0.50006$) verified at $40$-digit precision at every sampled impedance (\S\ref{sec:accept}), and the omitted-mode sums are directly evaluable from the closed-form coefficients \eqref{eq:normalized-coeff} without any enclosure constants (Appendix~\ref{app:tail}). The large-order asymptotics $|\lambda_n|=\tfrac12+\tfrac{\alpha R}{2|n|}+O(|n|^{-2})$ place the minimum at the smallest sampled impedance and the largest checked order --- exactly where the tabulated $0.50006$ is attained --- and make the threshold $n_2$ transparent. The scanned tail minima increase monotonically with the sampled impedance ($0.500063$ at $\alpha=0.05$, the overall minimum), consistent with the asymptotics.
\end{remark}

\begin{remark}[quadrature-format-dependent consistency]\label{rem:spectral}
For periodic analytic integrands, the periodized trapezoid rule is exponentially convergent \cite{Kress}, and a product-integration treatment of the logarithmic kernel preserves that geometric consistency mechanism. The implemented coincident Nystr\"om format is different: it replaces the logarithmic diagonal by $d_{\mathrm{eff}}=R/N$ and the normal-kernel diagonal by $-\kappa/(4\pi)$. These exact diagonal values repair the coincident entries but do not, by themselves, justify applying the analytic product-integration rate to the resulting nonsmooth matrix entry pattern. We therefore make no universal spectral-accuracy claim for that implementation. Under the non-resonance and stability assumptions below, consistency and $\|\Acal_h^{-1}\|_2\to\|\Acal^{-1}\|$ still follow from the collectively compact Nystr\"om framework; the MMS and refinement studies of \S\ref{sec:specquad} measure the realized consistency rate. Lemma~\ref{lem:rate} separates the geometric product-integration statement from the observed $O(h^3)$ law of the diagonal-corrected disk implementation.
\end{remark}
\begin{theorem}[Residual decomposition and effectivity for the disk surrogate]\label{thm:conv}
Fix $k>0$, $R$, and a margin $\eta>0$. Let $K\subset(0,\infty)$ be compact and restrict to the benign parameter set
\[
\mathcal{P}_{\eta,K}=\Big\{(\alpha,x_0):\ \alpha\in K,\quad
\mathrm{dist}\big(\alpha,\{\alpha_\ell(k)\}_\ell\big)\ge\eta,\quad |x_0|\le\bar r\Big\},
\qquad \bar r<1 .
\]
Let $\sigma^*$ be the exact density, $\sigma_h$ the coincident Nystr\"om solution with the exact diagonal corrections of \S\ref{sec:specquad}, and $\sigma_\theta$ the hybrid surrogate of \S\ref{sec:hybrid}: modes $|n|\le M_{\mathrm{fit}}$ are learned, modes $M_{\mathrm{fit}}<|n|\le M_{\mathrm{tail}}$ use the analytic Graf quotient, and modes $|n|>M_{\mathrm{tail}}$ are omitted. The experiments use $(M_{\mathrm{fit}},M_{\mathrm{tail}})=(48,220)$. Then, with constants depending only on $(k,R,\eta,K,\bar r)$:
\begin{enumerate}[leftmargin=1.6em, itemsep=1pt]
\item \textbf{discrete stability.} There is $N_0=N_0(k,R,\eta,K)$ such that for all $N\ge N_0$, $\Acal_h$ is invertible and
\[
\|\Acal_h^{-1}\|_2\le 2\|\Acal^{-1}\|_2\le \frac{2}{\min_n|\lambda_n(\alpha,k)|}\le C(\eta,K),
\]
where $\lambda_n$ is the explicit symbol of Proposition~\ref{prop:stability};
\item \textbf{format-dependent quadrature consistency.} Let $\mathcal E_{\mathrm{quad}}(h)$ denote a consistency bound for the exact solution,
$\big\|\Acal_h\,\sigma^*|_h-g_h\big\|_2\le \mathcal E_{\mathrm{quad}}(h)$.
If $\sigma^*$ continues analytically to a strip and the logarithmic kernel is treated by analytic product integration, then $\mathcal E_{\mathrm{quad}}(h)\le C h^{-3/2}\mathrm{e}^{-\gamma/h}$ for some $\gamma>0$. For the diagonal-corrected implementation used in the experiments, consistency is retained but this geometric rate is not asserted; its realized disk rate is measured in Lemma~\ref{lem:rate}(iii);
\item \textbf{surrogate residual decomposition.} The residual of \S\ref{sec:hybrid}, evaluated on a fine grid with meshwidth $h'$, satisfies
\[
\varepsilon_{\mathrm{res}}(\sigma_\theta)\ \le\ C_1\,\mathcal E_{\mathrm{quad}}(h')\ +\ C_2\frac{(k\bar r/2)^{M_{\mathrm{tail}}+1}}{M_{\mathrm{tail}}!}\ +\ C_3(M_{\mathrm{fit}})\,\delta_{\mathrm{fit}},
\]
where the constants absorb the uniformly positive and bounded $g$-normalization factors on $\mathcal P_{\eta,K}$, and $\delta_{\mathrm{fit}}$ is the aggregate learned-mode error
\[
\delta_{\mathrm{fit}}^2=\sum_{|n|\le M_{\mathrm{fit}}}\|\widehat c_n-c_n\|_{\Sigma_n}^2,
\]
with $\Sigma_n$ the fixed-design source-grid metric made explicit in Proposition~\ref{prop:fitrate};
\item \textbf{uniform effectivity in the benign zone.} for $(\alpha,x_0)\in\mathcal P_{\eta,K}$, with $r_\theta=\Acal_\alpha\sigma_\theta-g_{x_0}$ the exact residual,
\[
1\ \le\ \frac{\|\Acal^{-1}\|\,\|r_\theta\|}{\|\sigma_\theta-\sigma^*\|}\ \le\ \mathrm{cond}_2\big(\Acal_\alpha(k)\big)=\frac{\max_n|\lambda_n|}{\min_n|\lambda_n|}\ \le\ C(k,R,\eta,K)\ <\ \infty .
\]
With the $g$-normalized residual convention $\varepsilon_{\mathrm{res}}=\|r_\theta\|/\|g_{x_0}\|$ of \S\ref{sec:protocol}, the same statement reads
\[
1\ \le\ \frac{\|\Acal^{-1}\|\,\varepsilon_{\mathrm{res}}}{\|\sigma_\theta-\sigma^*\|/\|g_{x_0}\|}\ \le\ \mathrm{cond}_2\big(\Acal_\alpha(k)\big):
\]
the acceptance product bounds the $g$-normalized density error, and the computable factor $\|g_{x_0}\|/\|\sigma^*\|\in[\|\Acal^{-1}\|^{-1},\|\Acal\|]$ of \S\ref{sec:protocol} converts it to the $\sigma$-normalization.
\end{enumerate}
\textbf{Discrete-to-continuous bridge.} For a common $g$-normalized density metric, let $\eta_{\mathrm{res}}(h')$ bound the contribution of residual evaluation and interpolation, $\eta_{\mathrm{stab}}(h)$ the additional acceptance contribution caused by replacing the continuous stability constant by its discrete estimator, and $\eta_{\mathrm{tail}}(M_{\mathrm{tail}})$ the omitted-mode contribution. Thus all four terms below have the same normalization and units. A continuous acceptance test must use a stability upper bound and include these corrections, for example
\[
\widehat C_{\mathrm{stab}}\,\widehat\varepsilon_{\mathrm{res}}
+\eta_{\mathrm{res}}(h')+\eta_{\mathrm{stab}}(h)+\eta_{\mathrm{tail}}(M_{\mathrm{tail}})\le T.
\]
On the disk, Lemma~\ref{lem:tail} supplies a computable $\eta_{\mathrm{tail}}(M_{\mathrm{tail}})$ for fixed nonresonant $(\alpha,k)$ and $\bar\rho<1$; on general smooth boundaries the analogous term requires an analyticity-strip or coefficient-decay enclosure, and at corners it is not available in this form.
The experiments below report the leading discrete term $\|\Acal_h^{-1}\|_2\varepsilon_{\mathrm{res}}$, with the correction terms controlled asymptotically by the consistency analysis but not interval-enclosed. We therefore call it a \emph{computable discrete a-posteriori acceptance estimate}; it becomes a rigorous continuous bound only after computable correction bounds are added.

The three principal channels---quadrature, omitted-mode truncation, and surrogate fit---are nevertheless separately monitorable, and the continuous error--residual equivalence gives the stated effectivity once the exact residual and stability constant are used.
\end{theorem}
\begin{proof}[proof sketch]
\textbf{(i)} The corrected trapezoid Nystr\"om operators $\Acal_h$ converge to $\Acal$ in the collectively compact sense \cite{Kress}; the logarithmic part is a compact perturbation with its coincident diagonal treated separately. We do \emph{not} claim norm convergence of $\Acal_h\to\Acal$ on the weighted space. Away from the resonance curves, Proposition~\ref{prop:stability} gives a positive lower bound on $|\lambda_n|$; its large-$|n|$ asymptotics and continuity in $\alpha$ yield the uniform bound on the separated set with $\alpha\in K$. Discrete invertibility and the resolvent estimate then follow from collectively compact perturbation theory.

\textbf{(ii)} In the product-integration variant, the logarithmic singularity is integrated exactly against the trigonometric basis. For each collocation node the remaining consistency functional is the periodic trapezoid error for an analytic periodic integrand, giving geometric convergence uniformly in the nodes. This argument does not apply unchanged to the implemented diagonal-entry replacement; for that format we use collective compactness for consistency and determine the algebraic rate separately by refinement, rather than claiming that the diagonal corrections remove every non-analytic feature of the discrete scheme.

\textbf{(iii)} Subtract and add the actual hybrid approximation
\[
\sigma_\theta=\sum_{|n|\le M_{\mathrm{fit}}}\widehat c_ne^{\ii n\theta}
+\sum_{M_{\mathrm{fit}}<|n|\le M_{\mathrm{tail}}}c_ne^{\ii n\theta}
\]
through $\Acal_{h'}$. The $\Acal_{h'}\sigma^*$-versus-$g_{h'}$ term is bounded by $\mathcal E_{\mathrm{quad}}(h')$ as defined in (ii); the omitted modes $|n|>M_{\mathrm{tail}}$ are controlled by Lemma~\ref{lem:tail}; and the learned modes contribute $\|\Acal_{h'}\|\,\delta_{\mathrm{fit}}$, using $\sup_h\|\Acal_h\|<\infty$ from (i). The analytic modes $M_{\mathrm{fit}}<|n|\le M_{\mathrm{tail}}$ cancel exactly between $\sigma^*$ and $\sigma_\theta$.

\textbf{(iv)} The lower bound is the error--residual equivalence $\|\sigma_\theta-\sigma^*\|=\|\Acal^{-1}r_\theta\|\le\|\Acal^{-1}\|\,\|r_\theta\|$; the upper bound is $\|r_\theta\|=\|\Acal(\sigma_\theta-\sigma^*)\|\le\|\Acal\|\,\|\sigma_\theta-\sigma^*\|$. On the disk both norms are Plancherel-exact in the symbol $\lambda_n$ (Proposition~\ref{prop:stability}), and $\min_n|\lambda_n|\gtrsim\eta$ on $\mathcal P_{\eta,K}$ while $\max_n|\lambda_n|\le C(k,R,K)$ on the compact impedance window, giving the uniform bound.
\end{proof}
\begin{remark}[scope]
Theorem~\ref{thm:conv} is not a training or approximation theorem: $\delta_{\mathrm{fit}}$ is supplied by the surrogate. Its purpose is to locate that fit error within the residual budget and to state precisely which additional terms are required for continuous error control. The constants $C_1,C_2,C_3(M_{\mathrm{fit}})$, $\gamma$, and $N_0$ are not assigned numerical values in this paper; the effectivity statement concerns the exact continuous residual and stability constant, not the finite-mode acceptance estimate used in the experiments. Lemma~\ref{lem:rate} and Remark~\ref{rem:channels} make the quadrature and tail channels explicit on the disk. Proposition~\ref{prop:fitrate} and Theorems~\ref{thm:continuousfit}--\ref{thm:gaussianfit} separately give fixed-design and continuous-source estimation rates in their stated stochastic noise models; they are not unconditional training theorems for deterministic direct-solve data.
\end{remark}

\begin{lemma}[analyticity strip and quadrature-dependent consistency on the disk]\label{lem:rate}
Under the hypotheses of Theorem~\ref{thm:conv}, with $|x_0|\le\bar r$ and $\bar\rho:=\bar r/R<1$:
\begin{enumerate}[leftmargin=1.6em,itemsep=1pt]
\item the exact density continues analytically to the strip $|\mathrm{Im}\,\theta|<\gamma_0$ with the \emph{explicit} width $\gamma_0=\log(1/\bar\rho)$, and for every $0<\gamma<\gamma_0$,
\[
\sup_{|y|\le\gamma}\|\sigma^*(\cdot+\mathrm{i}y)\|_{L^\infty}\le C(k,R,|\alpha|,\gamma)<\infty;
\]
no uniform bound on the whole open strip is asserted;
\item \textbf{product-integration format.} If the logarithmic kernel is treated by analytic product integration, the exact solution's discrete residual obeys the explicit geometric estimate
\[
\big\|\Acal_h^{\mathrm{PI}}\,\sigma^*|_h-g_h\big\|_2\ \le\ C_1\,N^{3/2}\,e^{-\gamma_0 N},
\qquad
C_1=C_1(k,R,\eta,\bar\rho)<\infty\ \text{computable}.
\]
Equivalently, the standard strip formulation gives $C_\gamma N^{3/2}e^{-\gamma N}$ for every $\gamma<\gamma_0$; the endpoint rate displayed here follows directly from the Fourier aliasing sum below;
\item \textbf{implemented diagonal-corrected format: measured rate.} For the coincident Nystr\"om discretization actually used here, with logarithmic diagonal $d_{\mathrm{eff}}=R/N$ and normal-kernel diagonal $-\kappa/(4\pi)$, we do not assert the geometric estimate in (ii). A refinement experiment performed for the present discretization at $k=\alpha=R=1$ instead exhibits $\mathcal E_{\mathrm{quad}}(h)\asymp C_h h^3$ at $r_0=0.45$: over five successive dyadic levels the error ratio is approximately $8$, with $C_h\approx1.26465\times10^{-2}$ for $h=2\pi/N$ and less than $0.08\%$ drift. At $r_0=0.9$ the curve saturates near $3.4\times10^{-7}$ for $N\ge1024$, so the same asymptotic rate is not claimed at that operating point. Table~\ref{tab:quadrefine} summarizes these numerical consistency observations; they are not an interval-enclosed residual bound.
\end{enumerate}
\end{lemma}
\begin{proof}[proof sketch and numerical protocol]
\textbf{(i)} The coefficient bound of the proof of Lemma~\ref{lem:tail}, $|\sigma_n|\le C'_\sigma |n|\bar\rho^{|n|}$, shows that the Fourier series extends to $|\mathrm{Im}\,\theta|<\log(1/\bar\rho)$. On each closed substrip its sum is bounded by a finite $\gamma$-dependent geometric series; the constant can deteriorate as $\gamma\uparrow\gamma_0$.

\textbf{(ii)} After analytic product integration of the logarithmic part, each collocation consistency functional is the trapezoid-rule error for an integrand $f_i$ that is analytic and periodic in the strip $|\mathrm{Im}\,\theta|<\gamma_0$. Since $\Acal_\alpha$ acts diagonally with $|\lambda_n|$ bounded (Proposition~\ref{prop:stability}), the Fourier coefficients obey $|\widehat f_i(n)|\le C|n|\bar\rho^{|n|}$. The aliasing identity then gives, per node, $|T_N[f_i]-I[f_i]|\le 2\pi\sum_{j\ne0}|\widehat f_i(jN)|\le 4\pi C\,N\bar\rho^{N}/(1-\bar\rho)^2$; summing over the nodes in $\ell^2$ and absorbing the discrete quadrature weights gives the stated bound.

\textbf{(iii)} The present refinement protocol evaluates the same fixed continuous residual at successive dyadic grids without refitting the density or changing the acceptance tolerance. It reports the ratios and fitted prefactor stated above. This protocol calibrates a refinement bridge for the implemented format; it is not a proof of geometric convergence and is not used to replace the unenclosed continuous correction by zero.
\end{proof}
\begin{remark}[two explicit channels]\label{rem:channels}
For the product-integration variant, Lemma~\ref{lem:rate} upgrades the quadrature channel of Theorem~\ref{thm:conv} from an existence statement to an explicit geometric rate: the ratio is $\bar\rho=r_0/R$ --- the Graf ratio itself --- and the polynomial factor is $N^{3/2}$. For the diagonal-corrected format used in the experiments, the corresponding channel is instead calibrated by the measured $h^3$ law, with the documented near-boundary saturation; no geometric rate is imported from the product-integration proof. The tail channel is directly evaluable from the closed-form coefficients without enclosure constants (Appendix~\ref{app:tail}); for $M_{\mathrm{tail}}\gtrsim kR$ it sharpens to $(M_{\mathrm{tail}}+1)(k\bar r/2)^{M_{\mathrm{tail}}+1}/(M_{\mathrm{tail}}+1)!$ through the factorial decay of $J_n(k\bar r)$; the prefactor reflects the density bound $|\sigma_n|\le C|n|\bar\rho^{|n|}$ and cannot generally be absorbed into an $M_{\mathrm{tail}}$-independent constant. The remaining channel is the surrogate fit. Which channel leads the residual budget depends on the operating point: for the hybrid disk surrogate of \S\ref{sec:hybrid}, the tabulated fine-grid residuals sit at the level of the quadrature floor calibrated by the refinement study of \S\ref{sec:specquad} ($4.6\times10^{-7}$ at $N_{\mathrm{eval}}=512$), while the tail channel is closed explicitly ($9.5\times10^{-11}$ at $r_0=0.9$) and the learned modes are fitted well below this level, so the quadrature channel leads there; for the all-learned surrogates of \S\ref{sec:all-learned}, the fit channel sets the residual floor, as diagnosed by the mode-cutoff sweep of \S\ref{sec:all-learned} and the truncated-direct comparison of \S\ref{sec:asym}. The leading discrete term $\norm{\Acal_h^{-1}}_2\varepsilon_{\mathrm{res}}$ reported below is therefore quadrature-led on the disk hybrid operating point and fit-led on the general smooth boundaries of \S\ref{sec:all-learned}.
\end{remark}

The following stochastic estimation results specialize the continuous-source
error budget of the companion theory \cite{ErrorBudget} to the present
single-layer disk model.  They are deliberately separate from the
deterministic query-time acceptance estimate.

\begin{proposition}[aggregate ridge rate on the identifiable band]\label{prop:fitrate}
Fix a source interval $I$, a carrier threshold $\eta_c>0$, and the learned
band $|n|\le M_{\mathrm{fit}}$.  Consider the following \emph{stochastic data
model}, separate from the deterministic direct-solve pipeline.  At
$N_{\mathrm{tr}}$ fixed source radii $\rho_i\in I$, suppose
\[
\widetilde s_n(\rho_i)=s_n(\rho_i)+\zeta_{n,i},
\]
where the $\zeta_{n,i}$ are independent centered complex sub-Gaussian errors.
On the identifiable band
\[
\mathcal B_I(\eta_c)=\{n:\ |n|\le M_{\mathrm{fit}},\
|J_n(k\rho)|\ge\eta_c\ \text{for all }\rho\in I\},
\]
assume the common variance proxy
$\sigma_{n,i}^2\le\sigma_{\mathcal B}^2:=C\varepsilon_d^2/\eta_c^2$.
For each $n$, let $X$ be the common degree-$m$ Chebyshev design matrix and let
\[
\widehat\beta_n=(X^\top X+\mu I)^{-1}X^\top\widetilde s_n,
\qquad \mu\ge0,
\]
with $X$ of full column rank when $\mu=0$.  Write
$s_n=X\beta_n+e_{n,m}$ for the fixed-design representation error, let
$\widehat s_n$ denote the fitted polynomial, and define
\[
\Delta_{\mathrm{fit},\mathcal B}^2:=
\sum_{n\in\mathcal B_I(\eta_c)}\frac1{N_{\mathrm{tr}}}
\sum_{i=1}^{N_{\mathrm{tr}}}|J_n(k\rho_i)|^2\,
|\widehat s_n(\rho_i)-s_n(\rho_i)|^2,
\]
\[
\mathrm{dof}_{\mu,n}=\operatorname{tr}\big(X(X^\top X+\mu I)^{-1}X^\top\big),
\qquad d_{\mathrm{eff}}=\sum_{n\in\mathcal B}\mathrm{dof}_{\mu,n},
\]
\[
B_\beta^2=\sum_{n\in\mathcal B}\|\beta_n\|_2^2,
\qquad
\mathcal E_m^2=\sum_{n\in\mathcal B}\frac1{N_{\mathrm{tr}}}
\sum_{i=1}^{N_{\mathrm{tr}}}|J_n(k\rho_i)e_{n,m}(\rho_i)|^2 .
\]
Then, for nonempty $\mathcal B$, with probability at least $1-\delta$,
\[
\Delta_{\mathrm{fit},\mathcal B}^2
\le C_{\mathrm{SG}}\left[
\frac{\varepsilon_d^2}{\eta_c^2}
\frac{d_{\mathrm{eff}}\log(e|\mathcal B|/\delta)}{N_{\mathrm{tr}}}
+\frac{\mu}{N_{\mathrm{tr}}}B_\beta^2
+\mathcal E_m^2\right],
\]
where $C_{\mathrm{SG}}$ is the constant of the fixed-design sub-Gaussian
quadratic-form tail and is not silently set to one.  For analytic quotient
envelopes, $\mathcal E_m^2\le C|\mathcal B|\bar\rho_C^{-2m}$ for a Bernstein
parameter $\bar\rho_C>1$; on the constant-coefficient disk the exact quotient
is constant in $r_0$, so this approximation term can be zero.

The complementary low-mode blind set
$\{|n|\le M_{\mathrm{fit}}\}\setminus\mathcal B_I(\eta_c)$ is not controlled
by this estimate.  Its fit contribution is denoted by $\delta_{\mathrm{blind}}$
and must be independently bounded, excluded from the approximation class, or
retained as an explicit floor.  Only unidentifiable modes in the analytic
high-order range can be absorbed into Lemma~\ref{lem:tail}; a low-order zero
of $J_n(kr_0)$ cannot.
\end{proposition}
\begin{proof}
For each mode, the ridge prediction error decomposes into noise, bias, and
representation terms.  Since $|J_n|\le1$ on the real carrier grid, restoring
the carrier cannot increase the aggregate prediction norm.  The noise term is
bounded by $\|H_{\mu,n}\zeta_n\|^2/N_{\mathrm{tr}}$, where
$H_{\mu,n}=X(X^\top X+\mu I)^{-1}X^\top$; the fixed-design sub-Gaussian bound
and a union bound give the effective-dimension term.  For the bias,
\[
\frac1{N_{\mathrm{tr}}}\|(H_{\mu,n}-I)X\beta_n\|^2
=\frac{\mu^2}{N_{\mathrm{tr}}^2}\beta_n^\top
(\Sigma+\mu N_{\mathrm{tr}}^{-1}I)^{-1}\Sigma
(\Sigma+\mu N_{\mathrm{tr}}^{-1}I)^{-1}\beta_n
\le \frac{\mu}{N_{\mathrm{tr}}}\|\beta_n\|^2,
\]
with $\Sigma=X^\top X/N_{\mathrm{tr}}$.  The representation term is
$\mathcal E_m^2$ by definition.  Summation over $\mathcal B$ completes the
proof.
\end{proof}

Fixed-design prediction alone does not control the error between training
radii.  The following explicit design factor supplies the missing lifting.

\begin{lemma}[Chebyshev design lifting]\label{lem:designlift}
Set $q=m+1$, assume that $X$ has full column rank, and define
\[
\lambda_X:=\lambda_{\min}(X^\top X)>0,
\qquad \mathfrak L_X:=\sqrt{qN_{\mathrm{tr}}/\lambda_X}.
\]
For every degree-$m$ polynomial $p$ on $I$,
\[
\|p\|_{L^\infty(I)}\le
\mathfrak L_XN_{\mathrm{tr}}^{-1/2}\|Xp\|_2 .
\]
At the Gauss--Chebyshev nodes
$t_i=\cos((i-\tfrac12)\pi/N_{\mathrm{tr}})$ with $N_{\mathrm{tr}}\ge q$,
\[
X^\top X=N_{\mathrm{tr}}\operatorname{diag}(1,\tfrac12,\ldots,\tfrac12),
\qquad \mathfrak L_X=\sqrt{2q}=:\mathfrak L_* .
\]
\end{lemma}
\begin{proof}
If $a$ is the Chebyshev coefficient vector of $p$, then
$\|a\|_2\le\lambda_X^{-1/2}\|Xa\|_2$ and
$\|p\|_\infty\le\sqrt q\,\|a\|_2$.  Exact Gauss--Chebyshev quadrature of all
basis inner products gives the displayed Gram matrix and hence
$\lambda_X=N_{\mathrm{tr}}/2$.
\end{proof}

\begin{theorem}[continuous-source band bounds]\label{thm:continuousfit}
For each $n\in\mathcal B_I(\eta_c)$, let $\bar s_{n,m}$ be a best degree-$m$
sup-norm approximation of $s_n$ on $I$, with coefficient vector
$\bar\beta_n$.  Define
\[
\mathcal E_{m,\infty}^2:=
\sum_{n\in\mathcal B}\sup_{\rho\in I}
|J_n(k\rho)(s_n-\bar s_{n,m})(\rho)|^2,
\qquad B_s^2:=\sum_{n\in\mathcal B}\|\bar\beta_n\|_2^2,
\]
and distinguish the pointwise, integrated, and uniform band errors
\[
\begin{aligned}
\Delta_{\mathcal B}(\rho)^2&:=
\sum_{n\in\mathcal B}|J_n(k\rho)|^2
|\widehat s_n(\rho)-s_n(\rho)|^2,\\
\Delta_{\mathcal B,I}^2&:=\frac1{|I|}\int_I
\Delta_{\mathcal B}(\rho)^2\,\dd\rho,
\qquad
\Delta_{\mathcal B,I}^{\mathrm{sup}}:=\sup_{\rho\in I}
\Delta_{\mathcal B}(\rho).
\end{aligned}
\]
Under the general sub-Gaussian hypotheses of
Proposition~\ref{prop:fitrate}, with probability at least $1-\delta$,
\begin{equation}\label{eq:continuousfit}
\Delta_{\mathcal B,I}^{\mathrm{sup}}\le \mathcal E_{m,\infty}
+\mathfrak L_X\left[
\frac{C\varepsilon_d}{\eta_c}
\sqrt{\frac{C_{\mathrm{SG}}d_{\mathrm{eff}}
\log(e|\mathcal B|/\delta)}{N_{\mathrm{tr}}}}
+\sqrt{\frac{\mu}{N_{\mathrm{tr}}}}\,B_s\right].
\end{equation}
This is simultaneous in $\rho\in I$; the probability concerns only the
training noise.  For Gauss--Chebyshev nodes, $\mathfrak L_X=\mathfrak L_*$.
\end{theorem}
\begin{proof}
Split $\widehat s_n-s_n$ into the polynomial estimation/ridge error
$\widehat s_n-\bar s_{n,m}$ and the best-approximation error.  Apply the
per-mode fixed-design bound underlying Proposition~\ref{prop:fitrate} at
confidence $\delta/|\mathcal B|$, lift each polynomial error by
Lemma~\ref{lem:designlift}, and sum the squared modal errors.  The
approximation term is controlled directly in sup norm.
\end{proof}

A sharper bound is available under a more specific noise model, but its
query meaning is different.

\begin{theorem}[joint proper-Gaussian integrated-source bound]\label{thm:gaussianfit}
Use ordinary least squares ($\mu=0$) at $N_{\mathrm{tr}}\ge q$
Gauss--Chebyshev nodes, written in the normalized basis $\check T_0=1$ and
$\check T_j=\sqrt2\,T_j$ for $j\ge1$.  Suppose the independent noise variables
are proper complex Gaussian with common variance proxy
$\sigma_{\mathcal B}^2$.  Then, with probability at least $1-\delta$,
\begin{equation}\label{eq:gaussianfit}
\Delta_{\mathcal B,I}\le \mathcal E_{m,\infty}
+\sigma_{\mathcal B}
\sqrt{\frac{\pi\chi^2_{2|\mathcal B|q,\,1-\delta}}{4N_{\mathrm{tr}}}},
\end{equation}
where $\chi^2_{\nu,1-\delta}$ is the $(1-\delta)$-quantile of a chi-square
variable with $\nu$ real degrees of freedom.  This is a source-$L^2(I)$
statement: it does \emph{not} bound $\Delta_{\mathcal B,I}^{\mathrm{sup}}$ or
$\Delta_{\mathcal B}(\rho)$ at every source radius.
\end{theorem}
\begin{proof}
In the normalized basis, $X^\top X=N_{\mathrm{tr}}I_q$ and exact
Gauss--Chebyshev quadrature gives
$\frac1{N_{\mathrm{tr}}}\sum_i|p(t_i)|^2=\|\beta\|_2^2$.  Since
$(1-t^2)^{-1/2}\ge1$,
\[
\frac12\int_{-1}^{1}|p(t)|^2\,\dd t\le\frac\pi2\|\beta\|_2^2 .
\]
Thus the integrated stochastic band error is controlled without the sup-norm
factor.  The matrices $XX^\top/N_{\mathrm{tr}}$ are complex rank-$q$
orthogonal projections, so the joint noise quadratic form has law
$(\sigma_{\mathcal B}^2/2N_{\mathrm{tr}})\chi^2_{2|\mathcal B|q}$.  Adding the
deterministic approximation term proves \eqref{eq:gaussianfit}.
\end{proof}

\begin{remark}[noise model and norm conventions]\label{rem:fitnorms}
Direct Nystr\"om solves in the numerical pipeline produce deterministic
discretization errors, not automatically independent centered sub-Gaussian or
Gaussian observations.  The results above therefore describe an estimation
mechanism under their stated stochastic model; they do not convert
Lemma~\ref{lem:rate} into a probabilistic noise assumption or claim an
unconditional training guarantee for the experiments.  In particular,
\eqref{eq:continuousfit} is pointwise in the source radius, whereas
\eqref{eq:gaussianfit} is integrated over the source interval.  Neither is a
substitute for the deterministic per-query residual--stability estimate.
\end{remark}

\begin{corollary}[conditional residual rates with an explicit blind-mode floor]\label{cor:endtoend}
Under Theorem~\ref{thm:conv}, the product-integration channel of
Lemma~\ref{lem:rate}, and the general sub-Gaussian hypotheses above,
with probability at least $1-\delta$, simultaneously for every $\rho\in I$,
\[
\begin{aligned}
\|\Acal_h^{-1}\|_2\,\varepsilon_{\mathrm{res}}(\rho)
\le C\Bigg[&
N^{3/2}\bar\rho^{N}
+(M_{\mathrm{tail}}+1)\frac{(k\bar r/2)^{M_{\mathrm{tail}}+1}}{(M_{\mathrm{tail}}+1)!}
+\mathcal E_{m,\infty}+\delta_{\mathrm{blind}}\\
&+\mathfrak L_X\frac{\varepsilon_d}{\eta_c}
\Big(\frac{d_{\mathrm{eff}}\log(e|\mathcal B|/\delta)}{N_{\mathrm{tr}}}\Big)^{1/2}
+\mathfrak L_X\Big(\frac{\mu}{N_{\mathrm{tr}}}\Big)^{1/2}B_s\Bigg].
\end{aligned}
\]
For the proper-Gaussian least-squares protocol, define the integrated
acceptance error
$\overline B_{\mathrm{res}}^2:=|I|^{-1}\int_I
(\|\Acal_h^{-1}\|_2\varepsilon_{\mathrm{res}}(\rho))^2\,\dd\rho$.
Minkowski's inequality gives the analogous bound for
$\overline B_{\mathrm{res}}$, with the last two statistical terms replaced by
\[
\frac{\varepsilon_d}{\eta_c}
\sqrt{\frac{\pi\chi^2_{2|\mathcal B|q,\,1-\delta}}{4N_{\mathrm{tr}}}} .
\]
This Gaussian consequence is integrated in the source variable, not
pointwise.  In either route, a vanishing bound requires control of the
representation and blind-mode floors in addition to the deterministic
quadrature and tail hypotheses.  No unconditional end-to-end guarantee for
the deterministic training pipeline is claimed.
\end{corollary}

\begin{remark}[scope of the rates]\label{rem:ratescope}
The rates above are disk-specific in exactly two places: the explicit strip
width $\gamma_0=\log(1/\bar\rho)$ (the Graf ratio) and the diagonal symbol
$\lambda_n$.  On general smooth boundaries the same channel structure holds
only after an analyticity-strip or coefficient-decay enclosure for
$\sigma(\cdot;x_0)$ is supplied.  The certified numerical sample budgets of
\cite{ErrorBudget} concern a different two-layer parametrix-stripped model
and are not transferred to the present experiments.  General-domain rates and
joint training in $(\rho,\alpha,k)$ remain open.
\end{remark}

\section{Resonance-aware sampling}\label{sec:resonance}

\subsection{Operator interpretation}
On a general domain, a resonant wavenumber satisfies the Dirichlet-to-Neumann condition
\begin{equation}\label{eq:dtn}
(\Lambda(k)+\alpha)u=0\quad\text{on }\partial\Omega,
\end{equation}
and the associated single-layer operator can become singular. Away from these points, the unregularized residual--stability estimate is available; at an exact resonance, the classical Green's function need not exist and no finite stability constant can certify it.

For self-adjoint Robin Laplacians, ordered eigenvalues satisfy the Dirichlet--Neumann comparison
\begin{equation}\label{eq:squeeze}
\mu_j^N<\mu_j(\alpha)<\mu_j^D,
\qquad
\frac{\mathrm d\mu_j(\alpha)}{\mathrm d\alpha}
=\frac{\int_{\partial\Omega}|u_j|^2\,\dd s}{\int_\Omega|u_j|^2\,\dd x}>0
\end{equation}
under standard regularity and simple-branch assumptions \cite{RobinBVP,Kennedy}. This motivates avoiding spectral neighborhoods in training, but does not by itself provide a uniform forbidden-band width for surrogate accuracy. Weyl asymptotics concern spectral counting and average density on general domains; the more explicit spacing formulas used below are specific to balls and shells \cite{yang2026balls,HelmholtzAnnulus}.

\subsection{Explicit disk curves and screening}
On balls, the characteristic equation can be solved for the impedance:
\begin{equation}\label{eq:alpha-res}
\alpha_\ell(k)=\frac{\nu}{R}-\frac{kJ_{\mu}'(kR)}{J_{\mu}(kR)},
\qquad \mu=\ell+\nu.
\end{equation}
In two dimensions at $k=R=1$, the radial curve gives $\alpha_0=J_1(1)/J_0(1)\approx0.5751$. These curves are used to exclude neighborhoods during training and to construct pole-aware features. Their derivation and the associated spectral asymptotics are given in the companion analytic works \cite{yang2026balls}; they are not new results of the present surrogate paper.

For parameters away from the exact curves, a small direct Nystr\"om solve can be used as an empirical screening probe: exceptionally large density norms indicate an ill-conditioned sample. The thresholds used in our disk experiments are calibration heuristics, not geometry-independent error bounds. The discrete acceptance estimate, rather than the screening rule, remains the final acceptance instrument.

\noindent\emph{Theoretical calibration.} In the companion theory \cite{Ainv} these heuristics are upgraded to quantitative laws: the family $\Acal_{k,\alpha}$ is uniformly bounded in a safe zone $0<\alpha\le\alpha_{0}(\Omega,k)$, and near every resonance satisfies the two-sided law $c/\mathrm{dist}\le\norm{\Acal^{-1}}\le C/\mathrm{dist}$ with $\alpha$-independent constants, which on the disk sharpens, whenever the distance is attained at a unique nearest spectral mode, to the exact equality (proved self-contained in Appendix~\ref{app:ainv})
\[
\big\|\Acal_{k,\alpha}^{-1}\big\|=\frac{2}{\pi R\,\big|H_{n^{*}}^{(1)}(kR)\,J_{n^{*}}(kR)\big|}\cdot\frac{1}{\mathrm{dist}\big(\alpha,\Sigma(k)\big)},
\]
with $n^{*}$ the nearest spectral mode. The screening probe used above is the measured channel of $\norm{\Acal^{-1}g}$: its skirt exponent is universally $1$, including at exceptional points of multiply connected modal blocks \cite{Ainv,EPcomplex}, so the safe/degraded/failed grading is a theoretical calibration rather than a fitted threshold.

\subsection{Disk calibration and scope}
For the disk experiment at $k=1$, excluding $(0.555,0.60)$ during training produced large errors near $\alpha_0$ for a smooth quotient model. Density-norm amplification and increasing residual provided useful empirical signatures of this degradation. The corresponding safe/degraded/failed labels are specific to the reported model, source distribution, and discretization. On other domains we use the same idea of screening by conditioning, but do not transfer numerical thresholds without recalibration. General-domain sampling by coarse spectral computations and homotopy tracking is left for future work.

\medskip\noindent\emph{Validity scope of the screening probe.} The probe is the measured channel of $\norm{\Acal^{-1}g}$ \cite{Ainv}, and it is informative exactly when the failure mechanism is spectral: near a resonance skirt, $\norm{\Acal^{-1}}$ --- and hence the density norm --- amplifies by $\mathrm{dist}^{-1}$, while the approximation error stays bounded. The $k$-extrapolation failure of \S\ref{sec:kalpha} ($k{=}3.2$: field error $1.46$ with normalized density norm $1.03$) shows the converse: when the failure is approximation-theoretic --- the modal quotients depend on $k$ through oscillatory Bessel/Hankel factors rather than through the low-order rational structure they have in $\alpha$ --- no conditioning signature need accompany the error. We therefore use the density-norm probe only as an optional warning screen (Algorithm~1, step 2), never as a pass/fail instrument, and keep the residual--stability estimate as the sole acceptance test.

\section{Analytic benchmarks: formulas and validation}\label{sec:benchmark}

\noindent\emph{Why exact Robin kernels are a distinguishing feature.} The explicit series for Robin Green's functions on balls in all dimensions \cite{yang2026balls} and on annuli/spherical shells with two independent impedances \cite{HelmholtzAnnulus} are, to our knowledge, the only computable non-Dirichlet Green kernels available in the learned-Green literature: the benchmarks used by boundary-integral networks (e.g.\ the Poisson disk formula) are Dirichlet, and every surrogate cited in \S\ref{sec:robingap} is trained and validated without any exact impedance-domain truth. Here these kernels serve three roles that no refined solve can provide: \textbf{(a)} machine-precision reference fields, making the measured errors of \S\ref{sec:results} genuine accuracy statements rather than cross-discretization comparisons; \textbf{(b)} exact resonance spectra \eqref{eq:char}--\eqref{eq:delta-n}, which drive the sampling exclusions, pole-aware features, and stability constants of \S\ref{sec:resonance}; \textbf{(c)} the Plancherel-exact stability constant of Proposition~\ref{prop:stability}; and \textbf{(d)} a mode-wise \emph{error-source decomposition}: since every modal coefficient $c_n$ is available in closed form (Appendix~\ref{app:tail}), the error $\hat c_n-c_n$ of any computed density can be resolved mode by mode into low-mode (fit or solve error), mid-mode (quadrature/consistency), and tail (omitted-mode) contributions --- directly answering whether an observed error lives in the low frequencies, the high frequencies, or the tail, a diagnostic none of the Dirichlet-only references of \S\ref{sec:robingap} can perform without an exact kernel. The pipeline is therefore self-contained on the model domains: representation, training residual, acceptance test, error diagnosis, and reference truth all derive from the same explicit symbol.

Recent works provide these explicit series for Robin Green's functions on balls \cite{yang2026balls} and annuli \cite{HelmholtzAnnulus}; with stable modal evaluation and controlled truncation, they furnish independent analytic benchmarks in addition to refined boundary-integral calculations.

\subsection[Helmholtz on the ball]{Helmholtz on the ball \cite{yang2026balls}}

For $B_R(0)$, define the Robin trace
\begin{equation}\label{eq:btrace}
\Bcal\big[Z_n(kr)\big] = \Big(\alpha + \frac{n}{R}\Big)Z_n(kR) - k\,Z_{n+1}(kR),
\end{equation}
then
\begin{equation}\label{eq:helmholtz-series}
G(x,x_0) = \frac{\ii}{4}H_0^{(1)}\!\big(k|x-x_0|\big) - \frac{\ii}{4}\sum_{n=-\infty}^{\infty}\frac{\Bcal\big[H_n^{(1)}\big]}{\Bcal\big[J_n\big]}\,J_n(kr_0)\,J_n(kr)\,e^{\ii n(\theta-\theta_0)}.
\end{equation}
The series converges absolutely on compact interior subsets and is symmetric. The truncation order depends on $kR$, the source and evaluation radii, and the requested tolerance; it is not determined by $kR$ alone. The poles of the Green's function are given by the roots of
\begin{equation}\label{eq:char}
\Delta_\ell(k;\alpha) = \Big(\alpha + \frac{\ell}{R}\Big)J_{\mu}(kR) - k\,J_{\mu+1}(kR) = 0,
\qquad \mu = \ell + \nu,
\end{equation}
Here $\nu=d/2-1$ in dimension $d$; in 2D, $\nu=0$ and $\mu=\ell$, consistent with the Fourier--Bessel series above. The roots define the resonance condition used for training-time sampling (cf.\ \S\ref{sec:resonance}).

\begin{proposition}[Exact stability constant on the two-dimensional disk]\label{prop:stability}
On the disk $B_R\subset\RR^2$, the Robin operator $\Acal_\alpha=\tfrac12 I+\Kcal'+\alpha\Scal$ is diagonalized by the Fourier modes $e^{\ii n\theta}$ with eigenvalues
$\lambda_n(\alpha,k)=\tfrac{\ii\pi R}{2}H_n^{(1)}(kR)\,\Delta_n(k;\alpha)$, where $\Delta_n$ is the characteristic function \eqref{eq:char} (with $n$ in place of $\ell$). Consequently, for $k$ not a resonance,
\[
\big\|\Acal_\alpha^{-1}\big\|_{L^2(\partial\Omega)\to L^2(\partial\Omega)} = \frac{1}{\displaystyle\min_{n}\big|\lambda_n(\alpha,k)\big|}
\quad\text{exactly},
\]
A finite-mode numerical evaluation is only an estimate unless a separate tail bound shows that the discarded modes cannot contain the minimum; the reported experiments use $|n|\le80$ together with the finite-range high-precision scan and analytic large-order argument of \S\ref{sec:accept}. Near a resonance the minimum tends to zero and the stability constant diverges. The safe-zone uniform bound, the two-sided law in resonance neighborhoods, and the exact constant of the degenerate zone are established in the companion theory \cite{Ainv}.
\end{proposition}
\begin{proof}
The single and double layer potentials of a Fourier mode on a circle are Fourier modes of the same order (Graf addition theorem), so $\Acal_\alpha$ acts diagonally; the eigenvalue identity follows by evaluating $(\partial_n+\alpha)\Scal[e^{\ii n\theta}]$ at $r=R$ using the boundary traces \eqref{eq:trace} and the Wronskian $H_n^{(1)}J_{n+1}-H_{n+1}^{(1)}J_n = 2\ii/(\pi kR)$. Norm-exactness follows from Plancherel's identity on $L^2(\partial\Omega)$.
\end{proof}

\subsection[Annuli: Helmholtz and the two-parameter resonance spectrum]{Annuli: Helmholtz and the two-parameter resonance spectrum \cite{HelmholtzAnnulus}}

On the annulus, each angular mode is governed by a $2\times2$ system, and the basis must contain both regular and singular Bessel functions (the origin is excluded by the hole):
\begin{equation}\label{eq:helmholtz-annulus}
G = \frac{\ii}{4}H_0^{(1)}\!\big(k|x-x_0|\big) + \sum_{n=-\infty}^{\infty}\big(A_n J_n(kr) + B_n Y_n(kr)\big)e^{\ii n(\theta-\theta_0)},
\end{equation}
with Robin traces ($\mathcal{B}_1 = -\partial_r + \alpha_1$ at $r=R_1$; $\mathcal{B}_2 = \partial_r + \alpha_2$ at $r=R_2$):
\begin{equation}\label{eq:robintrace}
\mathcal{B}_1[Z_n] = \Big(\alpha_1 - \frac{n}{R_1}\Big)Z_n(kR_1) + k Z_{n+1}(kR_1),\qquad
\mathcal{B}_2[Z_n] = \Big(\alpha_2 + \frac{n}{R_2}\Big)Z_n(kR_2) - k Z_{n+1}(kR_2).
\end{equation}
The coefficients are given by a Cramer formula, and the zeros of the characteristic determinant
\begin{equation}\label{eq:delta-n}
\Delta_n = \mathcal{B}_1[J_n]\,\mathcal{B}_2[Y_n] - \mathcal{B}_2[J_n]\,\mathcal{B}_1[Y_n]
\end{equation}
are the resonance wavenumbers $k_{n,p}(\alpha_1,\alpha_2)$. The zeros of this determinant give the shell resonances. The companion work establishes their monotonicity, low-frequency resonance-free interval, and high-frequency spacing on annuli and spherical shells \cite{HelmholtzAnnulus}. Here we use only the forward modal system as an independent benchmark; the impedance-recovery theory is outside the scope of this paper.

\textbf{Our implementation and verification} ($k=5$, $\alpha_1=1$, $\alpha_2=2$, $R_1=1$, $R_2=2$, $r_0=1.5$): modal back-substitution residual $1.6\times10^{-17}$; interior residual of $(\Delta+k^2)v$ for the correction $v$ at $4.5\times10^{-7}$ (difference-limit); double-boundary Robin residual of $G=\varphi+v$ at $2.1\times10^{-4}$ (first-order difference limit, improvable to $10^{-8}$ by refinement); two-sided series versus folded implementation differ by $2.4\times10^{-16}$. The learned density network against this benchmark, and its acceptance-tested extension, are reported in \S\ref{sec:annulus}.

\subsection{Implementation, triple self-checks, and a numerical lesson}\label{sec:triplecheck}

We implemented \eqref{eq:helmholtz-series} and \eqref{eq:helmholtz-annulus} fully (two-sided series folded into cosine form) and applied triple self-checks ($k=3$, $\alpha=1$, $R=1$; Helmholtz truncation $N=60$):
\begin{center}\small
\begin{tabular}{@{}ll@{}}
\toprule
Check (Helmholtz disk, $k=3$, $\alpha=1$; truncation $N=60$) & value \\
\midrule
Equation residual $\Delta G + k^2G$ (interior, centered differences) & $7.5\times10^{-8}$ \\
Robin residual $\partial_nG + \alpha G$ (boundary, 4000 points) & $8.7\times10^{-6}$ \\
Series convergence ($N{=}20$ vs $N{=}60$) & \textbf{machine precision} \\
\midrule
Annulus Helmholtz (dual impedance, $k{=}5$): modal back-substitution & $1.6\times10^{-17}$ \\
harmonic residual / $G$ Robin residual (diff.\ limit) & $4.5\times10^{-7}$ / $2.1\times10^{-4}$ \\
\bottomrule
\end{tabular}
\end{center}

The series-convergence and modal back-substitution checks reach roundoff in the reported cases. The PDE and Robin residuals are limited by finite-difference verification and are therefore not machine-precision checks. Together with the independent Nystr\"om comparison, they provide the reference accuracy used below.

\textbf{Numerical evaluation of the boundary data.} In our implementation, the Graf-expanded function-value series was stable in the reported range, but direct evaluation of its term-by-term derivative series suffered severe cancellation from the factorial growth of $Y_n$. We therefore form $g$ from the closed-form fundamental solution and its analytic normal derivative (automatic differentiation gives the same expression), rather than differentiating the Graf series numerically. For $r<r_0$, the Graf convention places $J$ at the boundary radius and $H$ at the source radius; reversing this ordering gives the wrong expansion.

\textbf{Floating-point range.} Although the analytic series converges on interior compact subsets, its direct double-precision evaluation can fail when $J_n(kR)$ underflows and the ratio $\mathcal{B}[H_n^{(1)}]/\mathcal{B}[J_n]$ overflows. At $k=1$ and $r_0=0.9$, our unnormalized implementation distorted $|G|$ by more than two orders of magnitude beyond 50 modes. This is a numerical evaluation failure, not divergence of the exact series; stable evaluation requires normalized ratios or logarithmic scaling with a truncation/remainder check, and refined BEM cross-checks near-boundary sources.

\textbf{Stabilized implementation.} To enlarge the usable numerical range, we normalize modal coefficients by interface values: $J_n$ ratios $J_n(kr_0)/J_0(kr_0)$ and $J_n(kR)/J_0(kR)$ are computed by Miller's backward recurrence with periodic rescaling, $Y_n$ is represented logarithmically ($\ln\Gamma$ asymptotics), and the modal coefficients
\begin{equation}\label{eq:normalized-coeff}
c_n(r_0) = -\,\frac{J_n(kr_0)}{2\pi R\,H_n^{(1)}(kR)}\cdot\Big(1+\mathrm{i}\,\frac{Y_n(kR)}{J_n(kR)}\, b_n\Big),\qquad b_n = \frac{\big(\alpha + \tfrac{n}{R}\big) - k\,\tfrac{Y_{n+1}(kR)}{Y_n(kR)}}{\big(\alpha + \tfrac{n}{R}\big) - k\,\tfrac{J_{n+1}(kR)}{J_n(kR)}}
\end{equation}
are reconstructed in fully logarithmic combinations so no intermediate quantity leaves double precision. Implementation points (all found in actual debugging): $Y_n$ is uniformly negative when $kR<y_{n,1}$ (the sign does not alternate), and the signs of $J_n/Y_n$ and $Y_{n+1}/Y_n$ must be retained inside the combination factors. The displayed combination is algebraically identical to the quotient form $c_n(r_0)=q_n(\alpha)J_n(kr_0)$ of Appendix~\ref{app:tail}, since $\big(\alpha+\tfrac nR\big)-kH_{n+1}^{(1)}/H_n^{(1)}=\big(1+\mathrm{i}\,\tfrac{Y_n}{J_n}b_n\big)\,\mathcal B_\alpha[J_n](kR)/H_n^{(1)}(kR)$ (equivalently, $\frac{1+\mathrm{i}(Y_n/J_n)b_n}{1+\mathrm{i}Y_n/J_n}\,\mathcal B_\alpha[J_n](kR)/J_n(kR)$, using $H_n^{(1)}=(1+\mathrm{i}Y_n/J_n)J_n$); the two forms agree to roundoff for $n\le80$ (approximately $10^{-16}$ relative discrepancy in the tested range). In the practical configuration ($N{=}100$), field points $r<0.95$ verify at $8.7\times10^{-8}$ (source $r_0{=}0.45$) and $1.0\times10^{-7}$ (source $r_0{=}0.90$; Robin boundary residual $1.9\times10^{-4}$, difference limit). Only when both field point and source are extremely close to the boundary ($r\ge0.95$ and $r_0\ge0.9$) does the series slow down; there we still cross-validate against refined BEM.

\section{Numerical results}\label{sec:results}

\noindent\emph{Reading guide.} \S\ref{sec:kalpha} reports ordinary-accuracy joint $(k,\alpha)$ amortization with a density network and a matched-configuration baseline, while \S\ref{sec:certified} reports the per-impedance learned disk surrogate and its sampled acceptance test. \S\ref{sec:cavityscreen} then turns the same protocol to the application-driven screening raised in \S\ref{sec:intro}: an acceptance-tested scan of an elliptical cavity family with a dual-weighted bound for driven-response quantities and an eigen-frequency anchor read against the TESLA design point. The remaining direct-solve examples --- the ellipse, non-convex flower, polygon, annulus, and crescent in \S\ref{sec:ellipse}--\S\ref{sec:annulus} --- are validation or stress tests of the acceptance instrument, not learned-surrogate results. \S\ref{sec:all-learned} separately reports the general-geometry learning experiments: a concise secondary modal-flower check with an independent FEM cross-check, followed by the main source-point neural surrogate study on an asymmetric boundary.

\begin{center}\small
\begin{tabular}{@{}>{\raggedright\arraybackslash}p{3.4cm}>{\raggedright\arraybackslash}p{3.4cm}>{\raggedright\arraybackslash}p{2.6cm}>{\raggedright\arraybackslash}p{4.6cm}@{}}
\toprule
Operating point & Representation & Amortization & Accuracy status \\
\midrule
\S\ref{sec:kalpha} joint $(k,\alpha)$ network & MLP $4{\times}256$ + Fourier features & joint $(k,\alpha)$ and sources & measured field $4.5{\times}10^{-2}$--$1.4{\times}10^{-1}$ (in-domain); no acceptance claim \\
\S\ref{sec:certified} hybrid disk surrogate & structured modal ($|n|{\le}48$) + analytic Graf tail & sources; one fit per impedance & passes $10^{-4}$ discrete density test at $32$ impedances; measured field ${\le}6{\times}10^{-7}$ for $r_0{\le}0.9$ ($\alpha=1$) \\
\S\ref{sec:all-learned} all-learned surrogates & modal-on-flower / source-point network & sources (fixed $k,\alpha$) & measured field ${\le}4.4{\times}10^{-4}$ (flower), median $1.65{\times}10^{-4}$ (asymmetric); density bounds $5{\times}10^{-3}$--$3{\times}10^{-2}$ and $10^{-2}$--$10^{-1}$: conservative, no $10^{-4}$ acceptance \\
\bottomrule
\end{tabular}
\end{center}

\subsection{Experimental protocol}\label{sec:protocol}

We distinguish three quantities: (i) relative $L^2$ density error against an analytic series or a refined direct solve; (ii) relative $L^2$ field error on an interior grid, excluding $|x-x_0|<0.05$; and (iii) the discrete residual and inverse norm entering the acceptance estimate. Unless stated otherwise, the field grid on the disk is $0.05\le r\le0.97$. Residuals are reported as relative discrete $\ell^2$ quantities, $\varepsilon_{\mathrm{res}}=\norm{\Acal_h\sigma-g}_2/\norm{g}_2$. The field-level output additionally forms the absolute weighted residual norm $\|r_{h'}\|_{h',w}$ of \eqref{eq:field-est}; on constant-weight grids the two norms differ only by the common quadrature-weight factor, which cancels in the relative residual but not in the absolute field estimate. Hence $\norm{\Acal_h^{-1}}_2\varepsilon_{\mathrm{res}}$ bounds the $g$-normalized density error, which differs from the more common $\sigma$-normalization by the computable factor $\norm{g}_2/\norm{\sigma_{\mathrm{ref}}}_2\in[\norm{\Acal^{-1}}^{-1},\norm{\Acal}]$, $O(1)$--$O(10)$ in the reported regimes; applications holding a $\sigma$-normalized density tolerance should read the tabulated bounds with this conversion factor in mind. Measured errors are used only for validation; the acceptance estimate itself uses no reference solution. Throughout the experiments, the stability factor is computed on the coarse grid while the residual is evaluated on the fine grid; both grids lie away from resonance, so the two discrete inverse norms converge under refinement at the format-dependent consistency rate (Remark~\ref{rem:spectral}), and the MMS gate monitors exactly this channel.

The accuracy and timing baseline is a direct Nystr\"om solve at the same fixed geometry. Direct solves can be assessed by the same residual--stability estimate, so acceptance testing is not unique to learning. The surrogate's potential benefit is online amortization over source locations. In our single-core CPU implementation, a direct query at $N=512$, including assembly, factorization, solve, and evaluation at $10^4$ points, took approximately $210$--$260$ ms, compared with $10.9$ ms for a disk-surrogate query; training took approximately $100$--$115$ ms per impedance. These implementation-specific measurements correspond to approximately $20\times$ online speedup after offline training, not a complexity-theoretic comparison with an optimized BEM library. The amortized point is decisive for design loops: a $10^{3}$-candidate screen remains feasible by direct solves ($\approx$3--4 minutes at $\approx$210--260\,ms per candidate), but a $2.5\times10^{5}$-combination family (source$\times$impedance$\times$frequency) would require $\approx$15 hours, whereas the trained surrogate with a cached fine-grid operator answers the same family in $\approx$46 minutes at $\approx$11\,ms per candidate, with the acceptance estimate attached to every accepted candidate. On the flower, the learned surrogate's measured field errors are small but its discrete density bounds are much more conservative (\S\ref{sec:all-learned}).

\medskip\noindent\textbf{Tolerance as an application knob.} The acceptance tolerance $T$ is an application-supplied requirement, not a method constant chosen to make the surrogate pass; the verdicts at three reference tolerances are summarized in Table~\ref{tab:sweep} (density level; ranges as reported in \S\ref{sec:results}--\S\ref{sec:annulus}). Two facts matter. \emph{First}, at the application-level tolerance $T=10^{-3}$ the entire direct-solve suite (ellipse, flower, annulus, benign cascade) passes, while corner-defective polygonal solves and the ultra-near query remain rejected: the verdict tracks the \emph{defect class}, not the threshold. \emph{Second}, loosening $T$ does not convert the learned-surrogate verdicts --- the asymmetric surrogate's estimates ($5.2\times10^{-3}$--$3.0\times10^{-2}$) and the flower surrogate's ($4.6\times10^{-2}$--$1.8\times10^{-1}$) exceed even $T=10^{-2}$ in their tails, so the conclusions of \S\ref{sec:asym} and \S\ref{sec:all-learned} are unchanged by any plausible tolerance. For field-accuracy requirements we recommend the field-level output \eqref{eq:field-est}.

\medskip\noindent\textbf{Three-zone verdicts.} Throughout \S\ref{sec:results} we use the grading of the companion theory \cite{Ainv}: \emph{pass} (the bound meets the stated tolerance), \emph{degraded} (the bound sits within the next decade), and \emph{declared} (a defect class --- near-boundary density spike, corner-discretization defect, or resonance skirt --- is reported rather than smoothed over). The symbol $\times$ in the tables below covers both degraded and declared rows.

All experiments use CPU NumPy scripts, with PyTorch for the density-network runs. We report only configurations actually run; proposed near-boundary, three-dimensional, and cross-geometry extensions are not included as validated results.

\begin{figure}[t]
\centering
\includegraphics[width=0.74\textwidth]{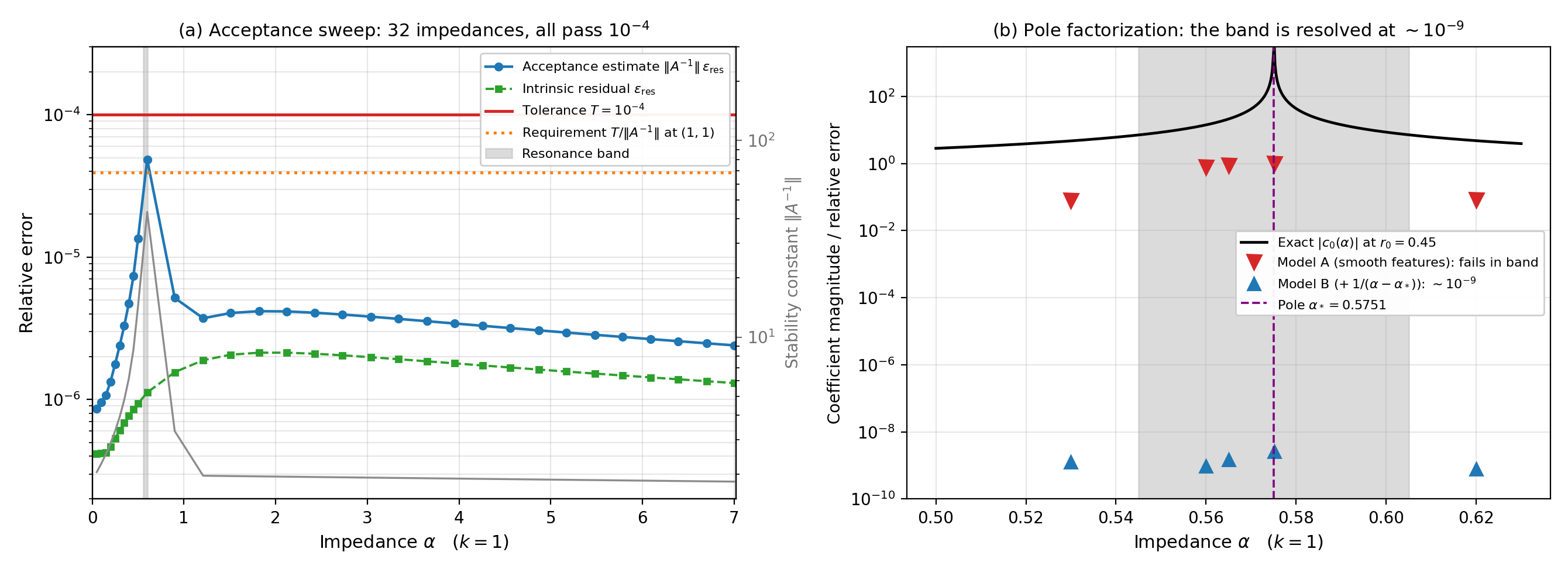}
\caption{\textbf{(a)} Discrete acceptance estimates at 32 sampled impedances outside the excluded resonance band: residual, acceptance estimate, target $T=10^{-4}$, and inverse-norm factor. \textbf{(b)} Coefficient magnitude and relative approximation error for the smooth and pole-aware quotient models. The tests are near-pole evaluations, not an evaluation of a classical Green's function at the exact resonance.}
\label{fig:accept}
\end{figure}

\medskip\noindent\emph{Reproduction coverage.} All tables and figures of this section except two are reproduced by the single-file release accompanying the paper (\texttt{slgreennet\_all\_in\_one.py}: fixed seeds, per-mode timing and validation anchors; the mode inventory is listed in Appendix~\ref{app:repro} and the script header): \textbf{(i)} the model A/B training behind the pole-model comparison of \S\ref{sec:polefact} lies outside the release (the pole panel of Figure~\ref{fig:accept}(b) is redrawn by the companion figure script from the recorded test errors), and \textbf{(ii)} the hand-placed cascade queries of Table~\ref{tab:asymcascade} (the regenerated counterparts used in Figure~\ref{fig:asymcascade} follow the protocol documented in the reproduction guide). The ellipse experiments of \S\ref{sec:ellipse}--\S\ref{sec:cavityscreen} are reproduced by the \texttt{ellipse} mode, and Table~\ref{tab:quadrefine} is reproduced by the \texttt{quad-refine} mode.

\begin{table}[H]
\centering\small
\caption{Tolerance sweep of the acceptance verdicts (density level). $\checkmark$: all queries accepted; $\times$: at least one rejected; partial: median accepted, tail rejected.}
\label{tab:sweep}\vspace{0.25cm}
\setlength{\tabcolsep}{4.5pt}
\resizebox{\textwidth}{!}{%
\begin{tabular}{@{}lccc@{}}
\toprule
Experiment & $T=10^{-4}$ & $T=10^{-3}$ & $T=10^{-2}$ \\
\midrule
Disk impedance sweep (32 impedances) & $\checkmark$ (worst $4.9\times10^{-5}$) & $\checkmark$ & $\checkmark$ \\
Disk near-boundary $r_0=0.95$ / $0.98$ & $\times/\times$ & $\times/\times$ & $\checkmark/\times$ \\
Ellipse direct solves & $\checkmark$ & $\checkmark$ & $\checkmark$ \\
Ellipse cavity screening (14 candidates, \S\ref{sec:cavityscreen}) & partial ($6/14$) & $\checkmark$ & $\checkmark$ \\
Flower direct solves & $\times$ (min $1.3\times10^{-4}$) & $\checkmark$ & $\checkmark$ \\
Annulus direct solves & $\checkmark$ & $\checkmark$ & $\checkmark$ \\
Square / L-shape stress tests & $\times$ & $\times$ & $\times$ \\
Asymmetric cascade $d=0.21/0.10/0.02$ & $\checkmark/\times/\times$ & $\checkmark/\checkmark/\times$ & $\checkmark/\checkmark/\times$ \\
Asymmetric learned surrogate (64 sources) & $\times$ & $\times$ (min $5.2\times10^{-3}$) & partial \\
\bottomrule
\end{tabular}}
\end{table}

\subsection{Ordinary-accuracy amortization over operator parameters}\label{sec:kalpha}
\subsubsection{\texorpdfstring{Joint $(k,\alpha)$ amortization on the unit disk}{Joint (k,alpha) amortization on the unit disk}}
The high-precision experiments below fit a structured surrogate separately at each impedance. As a complementary operating point, we also tested a density network that takes $(z,x_0,k,\alpha)$ as input and is amortized jointly over source locations and both operator parameters on the unit disk. The network is a $4\times256$ MLP with Fourier features. Training uses $N_q=64$, batch size 32, $\lambda=10^{-6}$, sources sampled uniformly by area in $r_0\le0.6$, $k\sim U[0.6,2.6]$, log-uniform $\alpha\in[0.05,5]$, and $20{,}000$ steps; operators are preassembled on a 61-point $k$ grid and linearly interpolated during training. Samples within $\pm0.03$ of the nine explicit resonance curves in the training domain are excluded. Evaluation uses 15 random sources per tabulated parameter point from the same source distribution and a refined direct solve at $N=512$, on a $25\times120$ radial--angular grid with $0.03\le r\le0.90$ and points with $|x-x_0|\le0.03$ excluded. The reported numbers are measured relative $L^2$ field errors, not residual--stability estimates. The density-norm column uses the mean norm of 15 direct solves at $N_q=64$, assembled at the evaluation wavenumber, normalized by the corresponding mean at $(k,\alpha)=(1,1)$.

Excluding the two resonance-skirt rows and the extrapolation row, the remaining ten grid points have mean field errors between $4.53\times10^{-2}$ and $1.37\times10^{-1}$: this is an ordinary-accuracy operating point, not the high-precision regime studied below. The two resonance-skirt points, $(0.7,0.3)$ and $(2.5,1.0)$, also carry the largest density-norm amplifications, $6.45$ and $3.75$. The moderately amplified point $(1.5,2.0)$ remains at $1.37\times10^{-1}$, however, so the norm probe is a useful warning diagnostic rather than a sharp universal pass/fail threshold.

The extrapolation point $(3.2,1.0)$ has a much larger error, $1.46$, even though its normalized direct-solve density norm is only $1.03$. Thus $k$-extrapolation failure need not be accompanied by the density-norm signature seen at the two resonance-skirt points. These 13 sparse parameter points do not establish uniform accuracy over the $(k,\alpha)$ plane, and none of the errors in this subsection is a residual--stability estimate. The resonance-band exclusions used during training are certified in the $k$-direction by the law $\norm{\Acal^{-1}}\le C(\Omega,\alpha)/\mathrm{dist}(k^{2},\sigma(-\Delta_{\alpha}))$ \cite{Ainv}; on general domains, certified spectral distances yield stability upper bounds that close the $\eta_{\mathrm{stab}}$ term of the discrete-to-continuous bridge.

\begin{center}\small
\setlength{\tabcolsep}{8pt}
\begin{tabular}{@{}cccl@{}}
\toprule
$(k,\alpha)$ & $\|\sigma_{\mathrm{dir}}\|/\|\sigma_{(1,1)}\|$ & mean field error & diagnosis \\
\midrule
$(0.7,0.3)$ & $6.45$ & $4.86\times10^{-1}$ & resonance skirt \\
$(0.7,1.0)$ & $0.42$ & $7.47\times10^{-2}$ & in-domain \\
$(0.7,3.0)$ & $0.32$ & $6.10\times10^{-2}$ & in-domain \\
$(1.0,0.1)$ & $0.87$ & $7.64\times10^{-2}$ & in-domain \\
$(1.0,1.0)$ & $1.00$ & $8.61\times10^{-2}$ & in-domain \\
$(1.0,4.0)$ & $0.41$ & $4.53\times10^{-2}$ & in-domain \\
$(1.5,0.5)$ & $0.76$ & $5.17\times10^{-2}$ & in-domain \\
$(1.5,2.0)$ & $2.99$ & $1.37\times10^{-1}$ & in-domain, amplified \\
$(2.0,0.7)$ & $1.07$ & $1.03\times10^{-1}$ & in-domain \\
$(2.0,1.5)$ & $0.79$ & $1.16\times10^{-1}$ & in-domain \\
$(2.5,1.0)$ & $3.75$ & $5.06\times10^{-1}$ & resonance skirt \\
$(2.5,3.0)$ & $0.88$ & $1.14\times10^{-1}$ & in-domain \\
$(3.2,1.0)$ & $1.03$ & $1.46$ & $k$-extrapolation \\
\bottomrule
\end{tabular}
\end{center}

\subsubsection{Like-for-like baseline at matched configuration}\label{sec:likeforlike}
To separate representation from network capacity, we use a matched ResNet baseline with the BI-GreenNet backbone and optimization configuration: ResNet $8\times100$, Fourier scales $\{1,2,4,8\}$, $N_q=64$, Adam $10^{-3}$, batch $32$, Tikhonov $10^{-6}$, and $20{,}415$ steps. The Robin problem and single-layer output representation are those of the present method, and the corrected coincident quadrature replaces the baseline's corrected-trapezoid rule; thus this is a matched-backbone comparison, not a reproduction of the complete BI-GreenNet algorithm. All errors use the same $240\times50$ polar evaluation grid and $N=512$ reference.

The baseline field errors are $9.5\times10^{-4}$ ($r_0=0.15$), $1.8\times10^{-3}$ ($r_0=0.30$), and $7.4\times10^{-4}$ ($r_0=0.45$), reached within the first $400$ steps. For context, the original \cite{BiGreenNet} benchmarks report $4.7\times10^{-3}$--$8.6\times10^{-3}$ on different problem settings, so those published numbers are not a direct performance comparison. The hybrid surrogate stays below $5.7\times10^{-7}$ for $r_0\le0.9$. The three-to-four-order gap in this matched disk configuration quantifies the value of the modal representation, analytic tail, and corrected coincident quadrature; it is not a cross-geometry ranking.

\subsection{Discrete acceptance: disk surrogate and direct-solve geometry validation}\label{sec:certified}

The discrete acceptance estimate $\norm{\Acal_h^{-1}}_2\varepsilon_{\mathrm{res}}$ combines a fine-grid residual with an inverse-norm estimate. All computations in this subsection are on the unit disk at fixed $k=1$. Unless explicitly stated, ``bound'' in the experimental tables means the leading discrete density estimate, without the continuous correction terms of Theorem~\ref{thm:conv}. The $g$-/$\sigma$-normalization convention of \S\ref{sec:protocol} applies to every table in this subsection.

\subsubsection{Corrected coincident quadrature with exact diagonal entries}\label{sec:specquad}
The Nystr\"om format of \S\ref{sec:nystrom} is upgraded from midpoint to \emph{coincident} collocation, with the two exact diagonal corrections---the logarithmic single-layer diagonal $d_{\mathrm{eff}}=R/N$ (derived from $\prod_{j=1}^{N-1}2\sin\tfrac{\pi j}{N}=N$; Appendix~\ref{app:nystrom}(iii)) and the bounded normal-kernel diagonal $-\kappa/(4\pi)$. The format is verified by the method of manufactured solutions (the free-space field of an exterior source $|\hat x|=4$, for which the density is known in closed form) and against the analytic series of \S\ref{sec:benchmark}:

\begin{center}\small
\begin{tabular}{@{}lcc@{}}
\toprule
$N$ & MMS: $\max|u-u_{\mathrm{ex}}|$ & floor: rel $L^2$ of $G$ vs $N=2048$ \\
\midrule
$64$  & $1.8\times10^{-6}$ & $4.8\times10^{-4}$ \\
$128$ & --- & $1.0\times10^{-4}$ \\
$256$ & --- & $1.3\times10^{-5}$ \\
$512$ & $3.6\times10^{-9}$ & $4.6\times10^{-7}$ \\
\midrule
\multicolumn{3}{l}{series ($N_t=40$) vs $N=2048$ direct solve: rel $L^2 = 3.3\times10^{-10}$} \\
\bottomrule
\end{tabular}
\end{center}

\begin{table}[H]
\centering\small
\caption{Independent quadrature-refinement anchors for the diagonal-corrected disk implementation at $k=\alpha=R=1$. The same fixed continuous residual is re-evaluated on successive grids; no surrogate is refitted. These observations calibrate consistency only and do not enclose the continuous residual correction.}
\label{tab:quadrefine}\vspace{0.25cm}
\setlength{\tabcolsep}{5pt}
\resizebox{\textwidth}{!}{%
\begin{tabular}{@{}clll@{}}
\toprule
$r_0$ & refinement range & observed error behaviour & calibrated quantity \\
\midrule
$0.45$ & $N=64,128,256,512,1024$ & ratio $\approx8$ per halving of $h$ & $C_h\approx1.26465\times10^{-2}$ in $C_hh^3$, $h=2\pi/N$; drift $<0.08\%$ \\
$0.90$ & $N=64,\ldots,2048$ & saturation of the refinement curve & $3.39674\times10^{-7}$ at $N=2048$ \\
\bottomrule
\end{tabular}}
\end{table}

At $N_{\mathrm{eval}}=512$, the measured relative $L^2$ field discrepancy against the refined solve is $4.6\times10^{-7}$, below the $10^{-4}$ acceptance target. This refinement study supports the discretization choice but does not replace the uncomputed continuous quadrature and interpolation corrections. Read against Theorem~\ref{thm:conv}(3), the tabulated residuals of \S\ref{sec:hybrid} sit at this floor: on the disk hybrid operating point the quadrature channel leads the residual budget, with the tail channel closed explicitly and the learned-mode channel fitted below it (Remark~\ref{rem:channels}).

\subsubsection{A hybrid modal surrogate}\label{sec:hybrid}
\noindent\emph{Status of the construction.} Per the definition of \S\ref{sec:surrogate}, the hybrid model \emph{is} a surrogate: the learned component carries the impedance-dependent quotient, while the analytic Graf tail carries the disk geometry; Table~\ref{tab:positioning} lists it side by side with the fully learned rows for a like-for-like reading.

As $r_0\to1$, the disk density develops an increasingly slowly decaying modal tail and approaches a boundary distribution. On the disk, the Graf addition theorem supplies this tail analytically; the quotient $q_n$ and the resulting tail formula are disk-specific, not geometry-independent. The surrogate therefore splits: modes $|n|\le48$ are learned using the carrier $J_n(kr_0)$ and Chebyshev quotients in $r_0$, fitted by standardized ridge regression; modes $49\le|n|\le220$ are appended from the analytic disk quotient $q_n$ (Appendix~\ref{app:tail}). The omitted modes $|n|>220$ are bounded by Lemma~\ref{lem:tail} with $M=220$; this closes the disk modal-tail term for $\bar\rho<1$, but not the corresponding omitted-mode term for the all-learned flower surrogate or for near-boundary sources. The fine-grid residual at $N_{\mathrm{eval}}=512$ supplies $\varepsilon_{\mathrm{res}}$:
\begin{center}\small
\setlength{\tabcolsep}{4pt}
\begin{tabular}{@{}ccccc@{}}
\toprule
$r_0$ & residual $\varepsilon_{\mathrm{res}}$ & discrete density bound $\norm{\Acal^{-1}}\varepsilon_{\mathrm{res}}$ & meas. density err. & meas. field err. \\
\midrule
$0.10$--$0.80$ & $\le 3.3\times10^{-7}$ & $\le 8.4\times10^{-7}$ & $\le 2.7\times10^{-7}$ & $\le 5.4\times10^{-7}$ \\
$0.90$ & $1.7\times10^{-6}$ & $4.3\times10^{-6}$ & $1.6\times10^{-6}$ & $5.7\times10^{-7}$ \\
$0.95$ & $3.9\times10^{-3}$ & $9.8\times10^{-3}$ (above tolerance) & $3.8\times10^{-3}$ & $1.5\times10^{-6}$ \\
$0.98$ & $6.9\times10^{-2}$ & $1.8\times10^{-1}$ (rejected) & $7.1\times10^{-2}$ & $1.2\times10^{-5}$ \\
\bottomrule
\end{tabular}
\end{center}

Because the exact coefficients are available mode by mode, the error of any computed density resolves into three readable blocks: the learned modes $|n|\le48$, each of which must be learned to the target accuracy individually; the analytic tail $49\le|n|\le220$, whose coefficients are each known to roundoff at essentially no cost; and the omitted tail $|n|>220$, which is never learned but bounded as a single explicit sum ($9.5\times10^{-11}$ at $r_0{=}0.9$). Comparing the learned $\hat c_n$ against the closed form thereby localizes an observed error before any remedy is attempted: this reading diagnosed the flower surrogate's floor as a representation/training effect (mode-cutoff sweep, \S\ref{sec:all-learned}) and the asymmetric network's floor as the $|m|\le32$ bandwidth (truncated-direct comparison, \S\ref{sec:asym}).

\medskip
The same mode-wise diagnostic separates the \emph{estimation} channel from the approximation floor in the stochastic noise model of Proposition~\ref{prop:fitrate}. Demodulating through the carrier turns relative data noise $\varepsilon_d$ into mode- and sample-dependent noise $\sigma_{n,i}=C\varepsilon_d/|J_n(k\rho_i)|$: noise is inflated where the carrier is small. A fixed carrier threshold $\eta_c$ defines the identifiable band, whose aggregate prediction error obeys the $N_{\mathrm{tr}}^{-1/2}$ rate of Proposition~\ref{prop:fitrate}. Theorem~\ref{thm:continuousfit} lifts the general sub-Gaussian statement to every source radius in $I$ through the explicit design factor $\mathfrak L_X$; Theorem~\ref{thm:gaussianfit} gives the sharper joint chi-square bound only for the source-$L^2$ error. In the evanescent high-order range, unidentifiable modes can instead be supplied by the analytic tail. A low-order carrier zero is different: its fit error remains in $\delta_{\mathrm{blind}}$ and cannot be assigned to that tail. The acceptance residual is a linear image of the density error, $r_\theta=\Acal_h\delta\sigma$, diagonalized by $\lambda_n$; direct symbol computation on the identifiable estimation-error family gives the spectral gains $0.72$--$0.74$ reported in our diagnostics. If the residual's own consistency floor exceeds the estimation error, that response must be de-floored before it becomes visible. These stochastic rates close only the identifiable-band channel, with the blind-mode floor stated separately; the deployed per-query residual remains a distinct deterministic diagnostic.

At $r_0=0.45$, the discrete estimate $5.3\times10^{-7}$ is about three times the measured density error $1.5\times10^{-7}$. For $r_0\le0.9$ (at $\alpha=1$), the measured field errors are below $6\times10^{-7}$. At $r_0=0.95$ and $0.98$, the density estimates exceed the $10^{-4}$ tolerance even though the measured field errors remain much smaller; these queries are therefore not accepted at that tolerance. Near-boundary acceptance requires spike-aware representation and near-singular quadrature. As a residual diagnostic, dropping the single mode $n=-49$ ($|c_{-49}|=8.8\times10^{-4}$) raised the residual from $10^{-6}$ to $1.5\times10^{-3}$ in our implementation.

\medskip\noindent\textbf{Field-level reading.} The absolute field-level estimate \eqref{eq:field-est} converts the absolute fine-grid residual into application-facing numbers. On the evaluation grid ($0.05\le r\le0.97$, $50\times240$ polar points, trapezoid boundary quadrature), the mean kernel factor is $0.502$. The grid-averaged estimates are $1.52\times10^{-7}$ at $r_0=0.45$ and $1.44\times10^{-6}$ at $r_0=0.9$ (maxima $1.67\times10^{-7}$ and $1.58\times10^{-6}$). The rejected rows report grid-averaged values $4.27\times10^{-3}$ ($r_0=0.95$) and $1.18\times10^{-1}$ ($r_0=0.98$). These are absolute pointwise-field estimates; the measured relative $L^2$ field errors above are validation data in a different metric and are not used for a cross-metric pass/fail test.

\subsubsection{Acceptance at sampled impedances}\label{sec:accept}
The acceptance criterion $\norm{\Acal_h^{-1}}_2\varepsilon_{\mathrm{res}}\le T$ with $T=10^{-4}$ is evaluated at 32 sampled impedances in $[0.05,7]$, excluding the resonance band used in training. One surrogate is fitted per sampled $\alpha$, and the maximum residual over $r_0\in\{0.2,0.45,0.9\}$ is reported. On the disk, the tabulated inverse-norm factor is the finite-mode estimate $\big(\min_{|n|\le80}|\lambda_n|\big)^{-1}$; by Lemma~\ref{lem:tail} the omitted large-order modes do not contain the stability minimum for sources with $\bar\rho<1$, and the density tail is explicitly bounded:

\begin{center}\small
\begin{tabular}{@{}cccc@{}}
\toprule
$\alpha$ & $\norm{\Acal^{-1}}$ & $\varepsilon_{\mathrm{res}}$ (max) & discrete density bound \\
\midrule
$0.05$ & $2.06$ & 4.17e-07 & 8.59e-07 \\
$0.10$ & $2.27$ & 4.20e-07 & 9.55e-07 \\
$0.30$ & $3.93$ & 6.09e-07 & 2.39e-06 \\
$0.50$ & $14.39$ & 9.37e-07 & 1.35e-05 \\
$0.60$ & $43.34$ & 1.12e-06 & 4.85e-05 \\
$0.90$ & $3.32$ & 1.56e-06 & 5.18e-06 \\
$1.51$ & $1.96$ & 2.07e-06 & 4.05e-06 \\
$2.12$ & $1.95$ & 2.13e-06 & 4.15e-06 \\
$3.04$ & $1.93$ & 1.98e-06 & 3.82e-06 \\
$4.87$ & $1.89$ & 1.62e-06 & 3.06e-06 \\
$7.00$ & $1.84$ & 1.31e-06 & 2.40e-06 \\
\bottomrule
\end{tabular}
\end{center}

\textbf{All 32 sampled impedances pass the discrete test.} The worst bound, $4.9\times10^{-5}$, occurs at $\alpha=0.60$, next to the excluded band, where $\norm{\Acal_h^{-1}}_2=43.3$. For the sampled values in $[2,5]$, the bounds are at most $4.2\times10^{-6}$. The table lists a representative subset. These statements concern the tested impedances and three source radii; they are not a uniform proof over the continuum $\alpha\in[0.05,7]$ or over all source locations. \emph{Clarification:} the excluded training band $(0.555,0.60)$ of \S\ref{sec:resonance} is open at $0.60$; the $\alpha=0.60$ surrogate is fitted on samples outside the band, which explains its enlarged stability factor and the worst acceptance bound of the sweep.

\medskip\noindent\emph{Finite-mode status.} Direct evaluation gives $\min_{81\le|n|\le142}|\lambda_n(\alpha)|\ge0.50$ in double precision and $\min_{81\le|n|\le400}|\lambda_n(\alpha)|\ge0.50$ at $40$-digit precision (mpmath) at every sampled impedance; the overall scanned minimum is $0.500063$. Lemma~\ref{lem:tail} supplies the large-order asymptotic $\lambda_n\to\tfrac12$ and a computable threshold beyond which the tail cannot contain the stability minimum. Since the smallest low-mode value is $1/43.34\approx2.3\times10^{-2}$, the scanned tail is separated from it by a wide margin. We therefore use the factor as a high-precision finite-mode stability estimate supported by an analytic tail argument, not as an interval-arithmetic certificate over every finite Bessel evaluation or over a continuum of impedances. For the larger impedances whose tabulated low-mode factors are $1.84$--$1.96$, we retain those finite low-mode values rather than replacing them by the floor-corrected value of approximately $2.00$ associated with the scanned tail minimum $0.500063$; this convention differs by at most about $2\%$ from the floor-corrected reporting in the companion theory and changes no acceptance verdict. The check is released as the \texttt{stab-tail} mode of the reproduction file, and the tail table accompanies the reproduction package.

The degradation pattern of the table is dictated by the exact constant of the two-sided law: the worst bound occurs at the band edge $\alpha=0.60$, where $\norm{\Acal^{-1}}=43.3$, while in the benign strip $\alpha\in[2,5]$ the stability factor attains its global minimum $\approx1.9$ \cite{Ainv}.

\subsubsection{Pole-aware features near a resonance}\label{sec:polefact}
The radial modal quotient $q_0(\alpha)$ has a pole at $\alpha_*=J_1(1)/J_0(1)\approx0.57508$. We compare a smooth quotient model with one augmented by $1/(\alpha-\alpha_*)$. Both are trained on $[0.50,0.63]\setminus(0.545,0.605)$ and evaluated at the following near-pole points:

\begin{center}\small
\begin{tabular}{@{}ccc@{}}
\toprule
$\alpha$ (test) & model A: smooth features only & model B: $+\,1/(\alpha-\alpha_*)$ \\
\midrule
$0.5300$ & 8.0e-02 & 1.3e-09 \\
$0.5600$ & 7.7e-01 & 9.7e-10 \\
$0.5650$ & 8.7e-01 & 1.5e-09 \\
$0.5751$ & 1.0e+00 & 2.6e-09 \\
$0.6200$ & 8.3e-02 & 8.0e-10 \\
\bottomrule
\end{tabular}
\end{center}

The entries are relative errors for the density coefficient. The value $0.5751$ is a rounded test point close to, but not equal to, $\alpha_*$; the large coefficient there remains finite. Thus the experiment shows that an explicit pole feature can improve approximation in a punctured neighborhood of a resonance. It does not establish existence of a classical Green's function at $\alpha_*$, nor does it provide a rigorous field or density bound at the pole.

\subsubsection{Direct-solve acceptance without disk symmetry: an ellipse}\label{sec:ellipse}
Everything in this example so far used the disk's symmetry. The discrete acceptance estimate is not tied to it: the corrected coincident Nystr\"om format holds on any smooth boundary (coincident collocation at arclength-equispaced nodes, with the same two diagonal corrections---the logarithmic diagonal $d_{\mathrm{eff}}=w/(2\pi)$ and the curvature diagonal $-\kappa/(4\pi)$), the fine-grid residual is evaluated on a finer quadrature by trigonometric interpolation of the density, and the stability factor is available as the discrete inverse norm $\norm{\Acal_h^{-1}}_2$---no closed-form symbol needed. As a demonstration, we apply the acceptance estimate to the \emph{direct solve} (no amortization, no symmetry) on the ellipse $a=1.3$, $b=0.7$ at $k=1$, $N=256$, against an $N_{\mathrm{eval}}=2048$ quadrature:

\begin{center}\footnotesize
\resizebox{\textwidth}{!}{%
\begin{tabular}{@{}ccccccc@{}}
\toprule
source $x_0$ & $\alpha$ & $\norm{\Acal_h^{-1}}_2$ & fine-grid residual & discrete density bound & meas.\ density err. & meas.\ field err. \\
\midrule
$(0.30,0.10)$ & $1.0$ & $2.20$ & $3.7\times10^{-7}$ & $8.1\times10^{-7}$ & $2.9\times10^{-7}$ & $1.3\times10^{-7}$ \\
$(-0.55,0.25)$ & $1.0$ & $2.20$ & $5.4\times10^{-7}$ & $1.2\times10^{-6}$ & $4.6\times10^{-7}$ & $1.4\times10^{-7}$ \\
$(0.30,0.10)$ & $2.0$ & $1.96$ & $1.9\times10^{-7}$ & $3.8\times10^{-7}$ & $3.1\times10^{-7}$ & $8.1\times10^{-8}$ \\
$(-0.55,0.25)$ & $2.0$ & $1.96$ & $3.9\times10^{-7}$ & $7.7\times10^{-7}$ & $8.6\times10^{-7}$ & $8.1\times10^{-8}$ \\
$(0.30,0.10)$ & $5.0$ & $1.90$ & $1.7\times10^{-7}$ & $3.1\times10^{-7}$ & $3.6\times10^{-7}$ & $1.1\times10^{-7}$ \\
$(-0.55,0.25)$ & $5.0$ & $1.90$ & $1.9\times10^{-7}$ & $3.6\times10^{-7}$ & $8.0\times10^{-7}$ & $9.2\times10^{-8}$ \\
\bottomrule
\end{tabular}
}
\end{center}

For these six direct-solve queries, the discrete bounds are $3.1\times10^{-7}$--$1.2\times10^{-6}$ and have the same order as the measured refined-solve discrepancies. Individual measured errors can slightly exceed the leading discrete estimate, so the table is not a uniform enclosure proof. The inverse norm is computed from the discrete operator offline and does not require the disk Fourier symbol. This example tests the acceptance estimate on a direct solve without symmetry; it is not a learned-surrogate experiment and does not establish geometry-independent stability constants. Near-boundary sources still require appropriate near-singular and spike-aware quadrature.

\subsubsection{Screening of an elliptical cavity family with per-candidate acceptance estimates}\label{sec:cavityscreen}
\noindent\emph{Application-driven experiment.} The ellipse experiment of \S\ref{sec:ellipse} is a single point of the design question raised in \S\ref{sec:intro}: in elliptical-cavity design one semi-axis is typically fixed by the target frequency and the other is the shape variable, so the designer screens a one-parameter family of geometries against excitation and wall-impedance requirements. The acceptance instrument supports this screening directly, at the per-candidate cost of a direct solve, with an acceptance estimate attached to each answer. We demonstrate it on the family
\[
\Omega_{a}=\bigl\{(x,y): x^{2}/a^{2}+y^{2}\le 1\bigr\},\qquad a\in[1.02,1.80],
\]
at $k=1$, $\alpha=1$, with source $x_{0}=0.3+0.1\ii$ and the query quantity $|G(x_{t},x_{0})|$ at $x_{t}=-0.45+0.05\ii$. Every candidate is processed by the protocol of \S\ref{sec:protocol}: coincident Nystr\"om solve at $N=64$, fine-grid residual at $N=256$, discrete inverse norm, and the field-level estimate \eqref{eq:field-est} (sub-second per candidate at this resolution, single core).

\begin{table}[H]
\centering\small
\caption{Screening of the fixed-height ellipse family with per-candidate acceptance estimates $\Omega_a$ at $k=1$, $\alpha=1$ (coarse $N=64$, fine residual $N=256$, reference $N=512$). The two field columns are absolute errors and estimates at $x_t$; the density columns retain the $g$-normalization of \S\ref{sec:certified}.}
\label{tab:cavity}\vspace{0.25cm}
\setlength{\tabcolsep}{2.5pt}
\resizebox{\textwidth}{!}{%
\begin{tabular}{@{}cccccccc@{}}
\toprule
$a$ & $\norm{\Acal_h^{-1}}_2$ & fine-grid residual & discrete density bound & abs.\ field err. & abs.\ field-level est. & $|G(x_t,x_0)|$ & verdict \\
\midrule
$1.02$ & $2.59$ & $1.2\times10^{-5}$ & $3.1\times10^{-5}$ & $8.1\times10^{-6}$ & $8.8\times10^{-6}$ & $0.666$ & pass $10^{-4}$ \\
$1.14$ & $2.90$ & $1.6\times10^{-5}$ & $4.5\times10^{-5}$ & $1.2\times10^{-5}$ & $1.3\times10^{-5}$ & $0.757$ & pass $10^{-4}$ \\
$1.26$ & $3.26$ & $2.1\times10^{-5}$ & $6.8\times10^{-5}$ & $1.7\times10^{-5}$ & $2.0\times10^{-5}$ & $0.859$ & pass $10^{-4}$ \\
$1.38$ & $3.70$ & $2.8\times10^{-5}$ & $1.03\times10^{-4}$ & $2.5\times10^{-5}$ & $3.0\times10^{-5}$ & $0.974$ & pass $10^{-3}$ \\
$1.50$ & $4.21$ & $3.7\times10^{-5}$ & $1.6\times10^{-4}$ & $3.7\times10^{-5}$ & $4.7\times10^{-5}$ & $1.106$ & pass $10^{-3}$ \\
$1.62$ & $4.84$ & $5.1\times10^{-5}$ & $2.4\times10^{-4}$ & $5.5\times10^{-5}$ & $7.3\times10^{-5}$ & $1.259$ & pass $10^{-3}$ \\
$1.74$ & $5.59$ & $6.9\times10^{-5}$ & $3.8\times10^{-4}$ & $8.2\times10^{-5}$ & $1.2\times10^{-4}$ & $1.441$ & pass $10^{-3}$ \\
$1.80$ & $6.04$ & $8.0\times10^{-5}$ & $4.8\times10^{-4}$ & $1.0\times10^{-4}$ & $1.5\times10^{-4}$ & $1.545$ & pass $10^{-3}$ \\
\bottomrule
\end{tabular}}
\end{table}

Two facts in Table~\ref{tab:cavity} carry the design content. \emph{First}, the query quantity increases monotonically over this window, and the accepted maximizer $a=1.80$ is separated from the runner-up by $0.104$ --- nearly three orders of magnitude above the maximum absolute field-level uncertainty $1.5\times10^{-4}$ over the sweep: the ranking, not only the values, is separated at the resolution of the reported discrete field-level estimates. This is not yet an interval-certified continuous ranking proof. Candidates whose gap fell within a few error bars would be flagged as indistinguishable at the stated tolerance, which is the intended use of the instrument in screening. \emph{Second}, the verdicts degrade gradually with elongation --- six of the fourteen swept candidates pass $10^{-4}$,\footnote{The six are $a\in\{1.02,1.08,1.14,1.20,1.26,1.32\}$; the remaining eight pass $10^{-3}$, the largest bound being $4.8\times10^{-4}$ at $a=1.80$.} the rest pass $10^{-3}$ --- tracking the discretization adequacy of the fixed $N=64$ mesh rather than any physical defect: without the acceptance estimate, the numbers for $a=1.02$ and $a=1.80$ would look equally trustworthy. The density-level effectivity stays between $2.5$ and $7$ throughout, an order of magnitude tighter than the learned surrogates of \S\ref{sec:asym}: direct solves at this accuracy level operate in the decision-grade band of \S\ref{sec:robingap}.

\begin{figure}[t]
\centering
\includegraphics[width=0.9\textwidth]{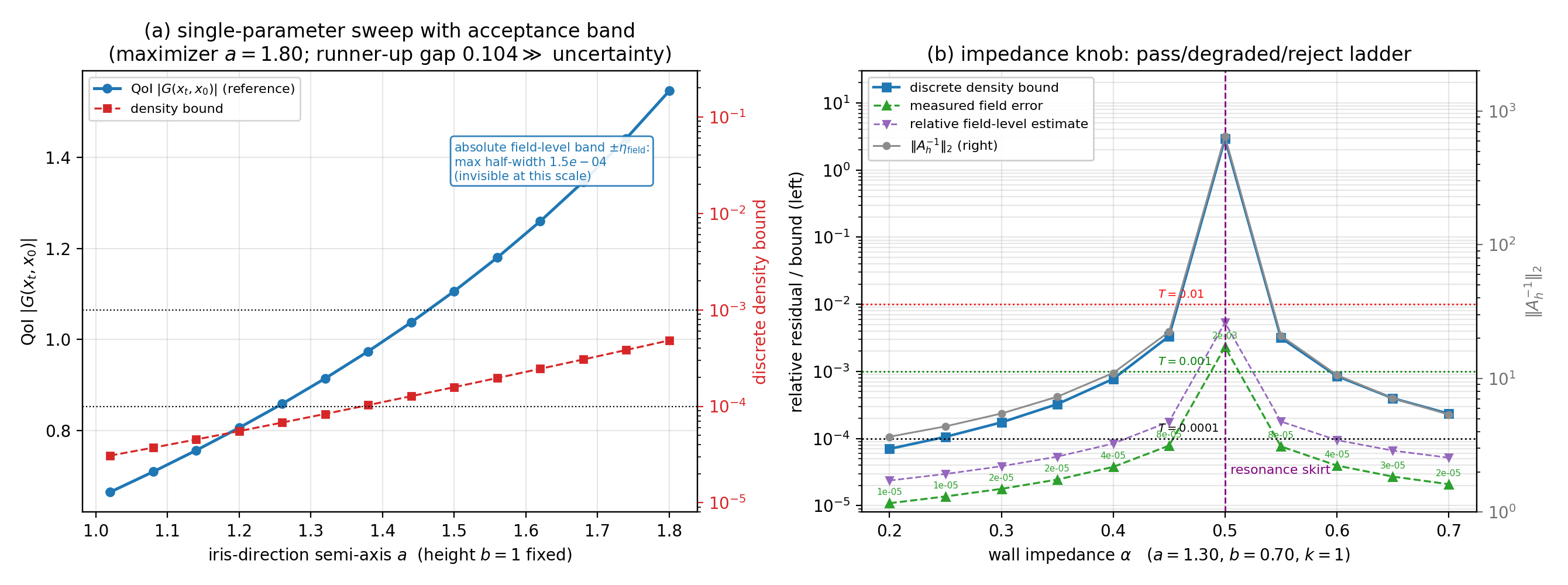}
\caption{Screening of the ellipse family with per-candidate acceptance estimates. \textbf{(a)} Single-parameter sweep with fixed height $b=1$: query quantity with the field-level band $\pm\eta_{\mathrm{field}}$ of \eqref{eq:field-est} (right axis: discrete density bound). \textbf{(b)} Impedance sweep at $(a,b)=(1.30,0.70)$ across the resonance skirt near $\alpha=0.50$: inverse norm and density bound, together with the measured relative pointwise field error and the equally normalized estimate $\eta_{\mathrm{field}}/|G^*|$; the reference tolerances apply to the density estimate, not to these two relative field curves.}
\label{fig:cavity}
\end{figure}

\medskip\noindent\textbf{The impedance knob and the resonance skirt.} Fixing $(a,b)=(1.30,0.70)$ and sweeping $\alpha\in[0.20,0.70]$ crosses a Robin resonance near $\alpha=0.50$ (where $k^{2}$ meets the Robin spectrum), and the instrument reproduces the full pass/degraded/reject ladder of \S\ref{sec:hybrid} without disk symmetry: $\norm{\Acal_h^{-1}}_2$ rises from $3.7$ to $645.7$, the density bound to $2.9$ at $\alpha=0.50$ (rejected at any stated tolerance), with $\alpha=0.45$ and $0.55$ degraded ($3.3\times10^{-3}$) and all remaining candidates passing $10^{-3}$ (Figure~\ref{fig:cavity}(b)). The measured field error stays at $10^{-5}$--$10^{-3}$ throughout: the estimate tracks the spectral defect class, not the approximation error --- precisely the behaviour the density-norm screening probe of \S\ref{sec:resonance} cannot be trusted to detect on its own. The query quantity itself peaks at $138.8$ at $\alpha=0.50$; this amplification is physical (the Green function is large near resonance), but the estimate states that this candidate's field is not controlled at the stated tolerance, and refinement cannot repair the verdict arbitrarily close to the pole.

\medskip\noindent\emph{Scope.} This is a direct-solve screening demonstration: no amortized geometry network is involved, the query quantity is a single field value, and the demonstration concerns one geometry family at one wavenumber. Design functionals of engineering interest (eigenfrequency placement, peak-field ratios) are $J(G)$-type quantities, to which the dual-weighted bound of Proposition~\ref{prop:dual} applies verbatim at one adjoint solve per quantity; extending the sweep to a learned-across-shapes surrogate with per-query acceptance is the cross-geometry amortization identified as future work in \S\ref{sec:conclusion}. All sweeps are reproduced by the \texttt{ellipse} mode of the reproduction file (\texttt{slgreennet\_all\_in\_one.py}).

\medskip\noindent\textbf{Driven-response discrete estimates are essentially sharp in this regime.} Replacing the point query by a driven-response quantity $J(G)=\langle G,\psi\rangle$ over a target patch ($\psi$ the normalized indicator of the disk of radius $0.08$ about $x_{t}$, $613$ quadrature points), the dual-weighted bound of Proposition~\ref{prop:dual} is nearly saturated: over the fourteen candidates the dual bound exceeds the measured functional error by a factor of only $1.0$--$1.3$, and the standard stability-times-kernel bound exceeds the dual bound by the benign-regime factor $1.06$--$1.10$ familiar from \S\ref{sec:asym} --- for direct solves, functional screening carries essentially no conservatism, and the estimate is a property of the candidate's own solve rather than of the family.

\medskip\noindent\textbf{Order-of-magnitude cavity comparison.} For context only, the TESLA nine-cell cavity has a $1.3$ GHz accelerating mode and half-cell/equator-arc dimensions $57.7/42.0$ mm \cite{TESLATDR}; its aspect ratio is of the same order as the scanned ellipse family. The computed first Robin eigen-wavenumber satisfies $k^*D\in[2.17,2.50]$, and $(a/k^*)\,dk^*/da\approx-0.21$ to $-0.28$. Thus a $0.1\%$ shape change corresponds to a relative frequency change of order $2.5\times10^{-4}$. This is not a quantitative TESLA simulation: the physical problem is an axisymmetric TM eigenmode with near-Dirichlet superconducting walls, whereas the impedance skirt above concerns absorber-wall physics. An axisymmetric Leontovich formulation is the required quantitative extension.

\subsubsection{Direct-solve acceptance on a non-convex flower}\label{sec:flower}
On the smooth non-convex flower $r(t)=\tfrac{9}{20}-\tfrac{\cos 5t}{9}$, also used by BINet \cite{BINet}, we compute $\norm{\Acal_h^{-1}}_2=2.84$ at $\alpha=1$. MMS errors decrease from $8.3\times10^{-6}$ to $2.2\times10^{-7}$ as $N$ increases from $128$ to $512$. For the three direct-solve queries in Table~\ref{tab:flower}, the discrete bounds at $N=512$ are $1.3\times10^{-4}$--$3.8\times10^{-4}$, with measured effectivity approximately $2$--$3$.

\begin{table}[H]
\centering\small
\caption{Discrete acceptance estimates for direct Nystr\"om solves on the flower (coarse $N=512$, reference $N=1024$, $\alpha=1$).}
\label{tab:flower}\vspace{0.3cm}
\resizebox{\textwidth}{!}{%
\begin{tabular}{@{}cccccc@{}}
\toprule
source $x_0$ & $\norm{\Acal_h^{-1}}_2$ & fine-grid residual & discrete density bound & meas.\ density err.\ & meas.\ field err.\ \\
\midrule
$(0.00,0.00)$ & $2.84$ & $4.4\times10^{-5}$ & $1.3\times10^{-4}$ & $5.4\times10^{-5}$ & $5.6\times10^{-6}$ \\
$(0.20,0.00)$ & $2.84$ & $1.3\times10^{-4}$ & $3.8\times10^{-4}$ & $1.6\times10^{-4}$ & $9.8\times10^{-6}$ \\
$(-0.15,-0.15)$ & $2.84$ & $5.9\times10^{-5}$ & $1.7\times10^{-4}$ & $6.8\times10^{-5}$ & $6.2\times10^{-6}$ \\
\bottomrule
\end{tabular}
}
\end{table}

\begin{figure}[t]
\centering
\includegraphics[width=0.75\textwidth]{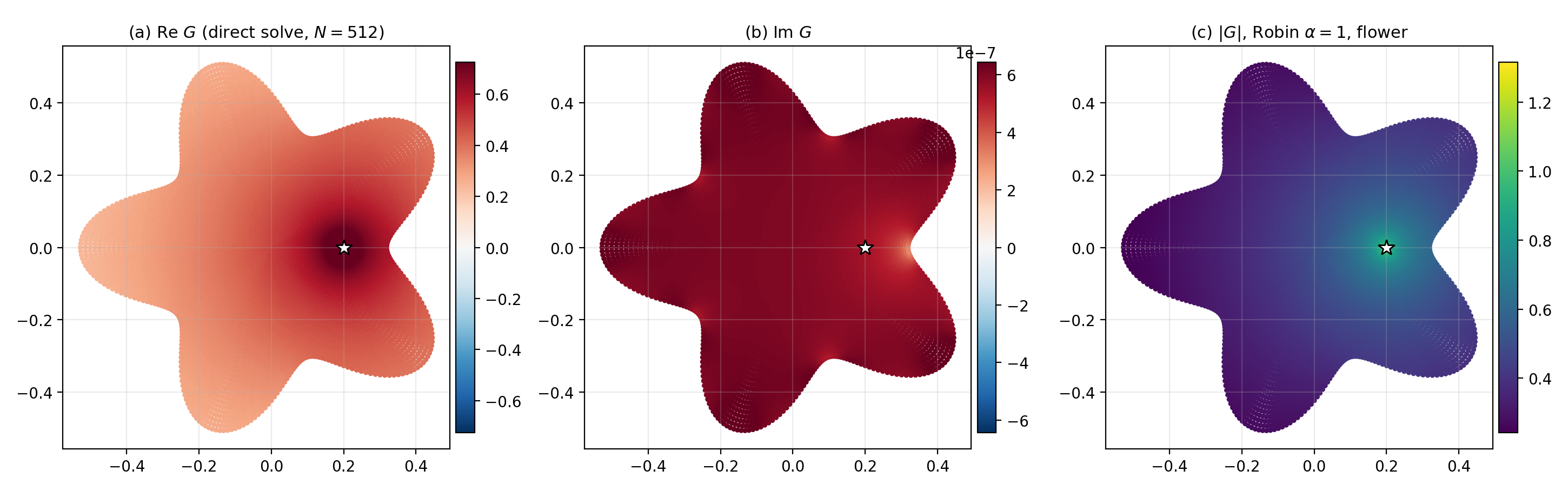}
\caption{Robin--Helmholtz Green's function on the flower at $\alpha=1$, computed by a direct solve at $N=512$: \textbf{(a)} $\operatorname{Re}G$; \textbf{(b)} $\operatorname{Im}G$; \textbf{(c)} $|G|$. The star marks $x_0=(0.2,0)$.}
\label{fig:flower}
\end{figure}

\subsubsection{High-frequency stress test: the disk at $k=185$}\label{sec:highfreq}
The experiments above fix $k=1$. Boundary integral methods do not suffer the FEM pollution effect, so the natural stress test of the \emph{instrument} (as opposed to a learned surrogate) is the same $10$\,kHz-scale acoustic regime adopted by the learning-based high-frequency reference \cite{ChenChengJCP2025}: we run the unit-disk Robin problem at $k=185$ with $\alpha=1$ and interior source $x_0=0.5$, coarse solves at $N$ and honest residuals on a fixed $N=4096$ grid.

Stability first: the Plancherel-exact symbol scan of Proposition~\ref{prop:stability} gives $\min_n|\lambda_n|=0.0301$, attained at the turning point $n^{*}=108$ --- a high-frequency turning-point minimum, not a resonance-skirt event --- whence $\|\Acal^{-1}\|_2=33.255$ exactly, with no impedance-curve exclusion required at this $(k,\alpha)$; the disk happens to sit in a wide spacing gap at this wavenumber, and a $0.5$-step scan of the exact symbol over $k\in[180,190]$ shows $\|\Acal^{-1}\|_2$ ranging $33$--$842$ there --- knowable only through Proposition~\ref{prop:stability}. The acceptance-estimate--$N$ curve:
\begin{center}\small
\begin{tabular}{@{}ccccc@{}}
\toprule
coarse $N$ & resid $\varepsilon_{\mathrm{res}}$ & bound $\|\Acal^{-1}\|\,\varepsilon_{\mathrm{res}}$ & meas.\ density err.\ & meas.\ field err.\ \\
\midrule
$512$  & $1.8\times10^{-3}$ & $6.1\times10^{-2}$ & $5.3\times10^{-3}$ & $5.5\times10^{-3}$ \\
$1024$ & $2.0\times10^{-4}$ & $6.5\times10^{-3}$ & $5.7\times10^{-4}$ & $5.9\times10^{-4}$ \\
$2048$ & $2.1\times10^{-5}$ & $7.0\times10^{-4}$ & $6.1\times10^{-5}$ & $6.3\times10^{-5}$ \\
\bottomrule
\end{tabular}
\end{center}
The bound improves by a factor $\approx9.3$ per doubling of $N$: the same $h^{3}$ consistency law calibrated in Table~\ref{tab:quadrefine}, so extrapolation to $N=4096$ places the bound in the $10^{-4}$ band. At $N=2048$ ($11$ nodes per wavelength) the query therefore passes the $10^{-3}$ acceptance with a measured field error of $6.3\times10^{-5}$.

Two remarks qualify the reading. \emph{MMS field metric at high frequency.} For the exterior-source MMS check at this frequency the relative field-error metric degenerates: the exact field is $\Phi_k(\hat x,\cdot)$ itself, the single-layer correction vanishes identically in $\Omega$, and the free-space magnitude $|H_0^{(1)}(k|\cdot-\hat x|)|/4\sim 8\times10^{-3}$ amplifies boundary-quadrature noise by two orders. Density-level quantities remain meaningful and are the honest high-frequency diagnostics (residual $7.8\times10^{-6}$, bound $2.6\times10^{-4}$, measured density error $2.6\times10^{-5}$). \emph{Frequency-matched comparison.} The learning-based reference at $k\approx184.75$ on a flower-shaped Dirichlet problem reconstructs a boundary-input solution operator with error-vector $2$-norm $1.06\times10^{-3}$ (maximum $7.5\times10^{-5}$) and no per-query error information \cite{ChenChengJCP2025}; the present point-source Robin Green queries --- a strictly harder data class --- attain $6.3\times10^{-5}$ measured field error with the acceptance bound $7.0\times10^{-4}$ attached to the same query; and on a qualitatively matched flower-shaped domain the present protocol reaches the reference's $10^{-3}$ field level at $9.4$ nodes per wavelength ($N{=}1024$) and surpasses it by an order of magnitude at $N{=}2048$. The problems and geometries differ, so this is a matched-\emph{frequency} instrument comparison, not a same-problem ranking.

\emph{Non-convex geometry.} The sweep was repeated on the non-convex flower $r(t)=\tfrac{9}{20}-\tfrac{\cos 5t}{9}$ of \S\ref{sec:flower} (perimeter $L\approx3.72$, interior source $x_0=0.2$), where the discrete inverse norm rises to $\norm{\Acal_h^{-1}}_2=107.3$ --- a moderately conditioned but non-resonant operator --- and the fine residuals are evaluated on $N=4096$:
\begin{center}\small
\begin{tabular}{@{}ccccc@{}}
\toprule
coarse $N$ & resid $\varepsilon_{\mathrm{res}}$ & bound $\|\Acal^{-1}\|\,\varepsilon_{\mathrm{res}}$ & meas.\ density err.\ & meas.\ field err.\ \\
\midrule
$512$  & $2.3\times10^{-3}$ & $2.5\times10^{-1}$ & $1.1\times10^{-2}$ & $9.9\times10^{-3}$ \\
$1024$ & $2.3\times10^{-4}$ & $2.5\times10^{-2}$ & $1.1\times10^{-3}$ & $1.0\times10^{-3}$ \\
$2048$ & $3.0\times10^{-5}$ & $3.2\times10^{-3}$ & $1.3\times10^{-4}$ & $1.2\times10^{-4}$ \\
\bottomrule
\end{tabular}
\end{center}
The bound improves by a factor $\approx7.8$ per doubling --- the $h^{3}$ law with a non-convex constant --- so the flower requires roughly one doubling more than the disk at the same frequency. The measured field error at $N{=}1024$ ($9.4$ nodes per wavelength) matches the $10^{-3}$-level headline of the learning-based reference on its own flower-shaped domain at the same frequency, and $N{=}2048$ improves it by an order of magnitude; the acceptance bound is attached at every level. The two flowers differ in petal count and size ($5$- versus $6$-petal), so the geometry match is qualitative.

\emph{Dense-resonance refusal on the asymmetric bean.} The same sweep on the asymmetric bean $r(t)=1+0.42\cos t-0.22\sin 2t$ of \S\ref{sec:asym} (perimeter $\approx6.86$, roughly twice the flower's) exposes the regime the instrument is designed to police: at $k=185$ the discrete inverse norm is $\norm{\Acal_h^{-1}}_2\approx971$ at $\alpha=1$, stable under refinement ($N=512,1024,2048$; the $N{=}4096$ operator agrees) --- a genuine spectral fact, not a discretization artifact --- and it remains $O(10^{2})$--$O(10^{3})$ at every sampled $\alpha\in[0.3,3]$. This is the high-frequency dense-resonance regime: Weyl counting places $\approx\!|\Omega|k^{2}/4\pi\sim 10^{4}$ Robin eigenvalues below $k^{2}=185^{2}$, so the impedance resonance curves crowd every $\alpha$ window. A wavenumber scan ($k\in[180,190]$, step $0.5$) confirms that the bean is not exceptional: $\|\Acal_h^{-1}\|_2$ fluctuates between $126$ and $1360$ (median $472$), the disk between $33$ and $842$ (median $226$) --- the same $10^{2}$--$10^{3}$ band, with the bean a mild factor $\approx\!3$ worse on median. The disk's benign $33$ at $k{=}185$ is a favorable draw made knowable by the explicit symbol; integrability does not sparsify the spectrum, it makes it auditable. The acceptance bound then exceeds $2.5\times10^{-1}$ at $N\le2048$ and no stated tolerance is met; crucially, the refined reference inherits the same conditioning, so \emph{no} solver output at this operating point is trustworthy at precision below $\sim\!10^{-3}$ without a certificate. An instrument-free method cannot distinguish this regime from a benign one at the solver level; the same spectral factor also degrades the a priori theory of the learning-based reference, whose stability constant carries $k^{3}/\mathrm{dist}$ \cite{ChenChengJCP2025}. The verdict here is refusal --- the designed behaviour of the last row of Table~\ref{tab:ladder} --- and it is reached on an asymmetric domain at high frequency without any symmetry to lean on.

\emph{Scope.} These are direct solves; a continuously certified \emph{learned} surrogate at high $k$ remains unreported (Table~\ref{tab:ladder}).

\subsubsection{Direct-solve stress tests on polygonal boundaries}\label{sec:polygon}
We next test the acceptance instrument outside the smooth-boundary setting, using the square $[-1,1]^2$ and an L-shaped domain with a reentrant corner. These are direct panel-Nystr\"om solves, not learned-surrogate results. Graded panels are concentrated toward the corners, and the residual is evaluated after panel-wise interpolation to a finer mesh. On the square, the manufactured-solution error decreases only algebraically, from $3.9\times10^{-5}$ to $9.2\times10^{-6}$ as the number of panels per edge increases from $4$ to $32$.

\begin{center}\small
\setlength{\tabcolsep}{4pt}
\begin{tabular}{@{}llcccc@{}}
\toprule
Domain & $\alpha$ & $\norm{\Acal_h^{-1}}_2$ & fine-grid residual & discrete density estimate & measured field error \\
\midrule
Square & $1$ & $15.4$ & $1.1$--$1.2\times10^{-1}$ & $1.8$--$1.9$ & $4.7$--$5.0\times10^{-6}$ \\
Square & $2$ & $7.9$ & $2.6$--$3.0\times10^{-2}$ & $0.21$--$0.23$ & $4.3$--$4.7\times10^{-6}$ \\
L-shape & $1$ & $12.53$ & $9.1\times10^{-2}$ & $1.14$ & $4.8$--$4.9\times10^{-6}$ \\
L-shape & $2$ & $10.24$ & $3.2$--$3.3\times10^{-2}$ & $0.32$--$0.34$ & $4.6$--$5.0\times10^{-6}$ \\
\bottomrule
\end{tabular}
\end{center}

The ranges cover three tested sources at each impedance. None of these queries meets a high-precision density tolerance, so the appropriate verdict is rejection rather than acceptance. The residual also localizes the defect: on the square, the eight corner-adjacent panels rise above the median panel residual by a factor of up to $7.7\times10^{4}$; on the L-shape, the two panels meeting at the reentrant corner are the largest. The small measured field discrepancies do not turn the large density estimates into rigorous field bounds. Rather, these examples separate density-level rejection from field-level accuracy and show that the instrument can identify corner-discretization defects even when the smooth-density theory no longer applies.

\medskip
A focused diagnosis supports this reading. For the L-shape at $(\alpha,p)=(1,24)$ the fine-grid residual is corner-dominated: the two panels adjacent to the reentrant corner carry $45.6\%$ of its mass, and the three largest inter-panel jumps (at the corners) exceed the background by three orders of magnitude. Three natural repairs were tested and quantitatively excluded. \textbf{(a)}~Post-hoc singularity subtraction with the two-term Grisvard span $\{\chi(r)\,r^{\lambda-1}:\lambda\in\{\tfrac23,\tfrac43\}\}$ --- coefficients fitted either from the coarse solution or optimally against the fine-grid residual --- reduces the bound by at most $20\%$ ($1.14\to0.85$), the map from this span to the residual having only two active directions. \textbf{(b)}~Enforcing $C^{0}$ continuity of the panel interpolant across breakpoints is incompatible with the unbounded corner density and degrades the solve. \textbf{(c)}~Mesh refinement alone follows $\mathrm{bound}\asymp p^{-1.3}$, so reaching $10^{-3}$ would require $\sim\!5000$ panels per edge. The viable route is a corner-adapted \emph{coarse} representation --- singularity-adapted degrees of freedom in place of the corner panels, in the spirit of the corner singular samples of the approximation-based learning line \cite{ChenChengJCP2025,ChenChengCorner} --- which is left as future work. This campaign is a stress test of the instrument, not a validated corner surrogate.

\subsection{All-learned surrogates on general smooth boundaries}\label{sec:all-learned}
The fixed-geometry workflow is summarized below. Screening is an empirical preprocessing step; acceptance is determined only by the discrete residual--stability estimate.

\medskip\noindent\textbf{Algorithm 1: fixed-geometry surrogate construction and query acceptance}
\begin{enumerate}[leftmargin=1.8em,itemsep=1pt]
\item \textbf{Discretize and verify.} Assemble $\Acal_h$ with the jump, logarithmic-diagonal, and curvature corrections, and gate the assembly by an MMS or independent benchmark check.
\item \textbf{Screen parameters (optional).} Use a small direct solve to flag exceptionally large density norms near resonance skirts. This is a heuristic screen, not an error bound.
\item \textbf{Generate offline data.} Compute direct-solve densities at training sources in a benign interior region, reusing the operator factorization at each fixed parameter value.
\item \textbf{Fit the surrogate.} Approximate the density by a structured modal expansion or a density network at fixed geometry.
\item \textbf{Accept or reject online.} Upsample the predicted density, evaluate the fine-grid residual and the discrete inverse-norm estimate, and accept only when the resulting discrete density estimate meets the prescribed tolerance.
\end{enumerate}

\noindent\textbf{Fallback.} A query that fails the acceptance test is answered by a direct Nystr\"om solve at the fine grid, for which the same residual--stability estimate is available; the surrogate is bypassed for that query. A rejected query is therefore always routed to the direct baseline. This fallback statement does not assert that an accepted surrogate is more accurate than a direct solve; acceptance only records that the stated discrete criterion has been met.

\medskip\noindent\textbf{Secondary flower-boundary check.} An all-learned modal experiment on the flower boundary gave measured field errors at most $4.4\times10^{-4}$, while its discrete density estimates ranged from $4.6\times10^{-2}$ to $1.8\times10^{-1}$. At $\alpha=1$, increasing the boundary-mode cutoff from $M=40$ to $56$ reduced the worst estimate only from $1.5\times10^{-1}$ to $1.2\times10^{-1}$; together with the truncated-direct comparison below, this identifies representation and fitting, rather than the inverse-norm factor, as the residual-limiting mechanism. At the checked query, an independent $P1$ FEM calculation agreed with the fine BEM reference to $2.5\times10^{-5}$ and with the learned surrogate to $5.4\times10^{-5}$. This is an independent field comparison, not a continuous error bound. The main general-geometry acceptance study is the asymmetric experiment below.

\subsubsection{Asymmetric smooth boundary: a source-point neural surrogate}\label{sec:asym}
We finally test a neural source-point surrogate on a smooth asymmetric boundary
without using rotational or reflection symmetry:
\[
z_b(t)=r(t)e^{\ii t},\qquad r(t)=1+0.42\cos t-0.22\sin(2t),\qquad 0\le t<2\pi .
\]
The setting is fixed-geometry and fixed-parameter ($k=1$, $\alpha=1$). The network input is the source coordinate with multiscale Fourier features; the output is the $65$ complex Fourier coefficients of the whole boundary density for $|m|\le32$. These coefficients are taken with respect to the equi-arclength boundary sampling coordinate, not the polar angle $t$ in $z_b(t)$; the same coordinate is used in truncation, synthesis, and upsampling.
Training uses $1600$ direct solves at $N=192$ ($160$ validation sources) and the combined coefficient-regression and coarse-operator residual objective of \S\ref{sec:surrogate}, with residual weight $0.1$; testing
uses $64$ independent sources in $0.10\le\rho\le0.58$, with minimum distance
$0.217$ to the boundary, so this is a safe-interior-source experiment rather than
a near-boundary certificate. The distance-resolved cascade at the end of this subsection supplies the near-boundary counterpart on the same boundary. Field errors are measured on $0.04\le\rho\le0.90$
with $|x-x_0|\le0.05$ excluded.

\begin{center}\small
\begin{tabular}{@{}lcccc@{}}
\toprule
Quantity & min & median & $90\%$ quantile & max \\
\midrule
full-field relative $L^2$ error &
$3.51\times10^{-5}$ & $1.65\times10^{-4}$ & $2.87\times10^{-4}$ & $1.29\times10^{-3}$ \\
density relative $L^2$ error &
$1.37\times10^{-3}$ & $1.99\times10^{-3}$ & $4.13\times10^{-3}$ & $9.00\times10^{-3}$ \\
discrete residual--stability estimate &
$5.22\times10^{-3}$ & $7.25\times10^{-3}$ & $1.41\times10^{-2}$ & $3.00\times10^{-2}$ \\
\bottomrule
\end{tabular}
\end{center}

The acceptance estimates in this experiment use the same relative residual $\norm{r}_2/\norm{g}_2$ as the other tables. The conclusion is unchanged:
all $64$ unseen sources exceed the strict $T=10^{-4}$ discrete density tolerance.
Measured full-field errors are far smaller: $63/64$ sources are below
$5\times10^{-4}$, with a single $1.286\times10^{-3}$ tail point, so we do not claim
uniform $10^{-4}$ field accuracy. With the stability factor
$\|\Acal_{192}^{-1}\|_2\approx2.674$ computed on the coarse grid (see
\S\ref{sec:protocol} for the coarse/fine-grid convention), this again separates field
accuracy from conservatism of the leading discrete density estimate; it is not a
continuous field certificate. The released implementation also evaluates \eqref{eq:field-est} with the absolute weighted residual and writes its grid mean and maximum for every source. The density-level summary above does not specify an absolute field tolerance, so we do not infer a $T_{\mathrm{field}}$ verdict by comparing it with the relative $L^2$ field errors; absolute field-acceptance statistics require the regenerated per-source field output rather than the $g$-normalized density estimate alone.

\medskip\noindent\textbf{Dual-weighted discrete estimate (Proposition~\ref{prop:dual}).} For a fixed target region --- the annular sector $0.25\le\rho\le0.55$, $|\arg z|\le\pi/4$ --- one adjoint solve at $N{=}512$ gives $\|z_\psi\|_2=2.087$ against $\|\Acal_h^{-1}\|_2\,\|\Scal_k^{*}\psi\|_2=2.207$: in this benign regime ($\|\Acal_h^{-1}\|_2\approx2.67$, far from resonance) the adjoint norm nearly saturates the stability-times-kernel product (ratio $1.06$), so dual weighting tightens the functional estimate only modestly here; its value is direction-selectivity, which grows as the spectrum is approached. For all $64$ unseen sources the measured functional error lies below the computed quantity $\|r_\theta\|_2\,\|z_\psi\|_2$, with median effectivity $22.8$ versus $24.1$ for the standard density estimate. This is an observed discrete effectivity audit, not an interval-enclosed continuous bound; \S\ref{sec:discussion} decomposes the remaining conservatism quantitatively. The released \texttt{dual-demo} mode reproduces this table. Truncating the direct $N=192$ density to the same $65$ modes without
the network gives field errors of median $1.89\times10^{-6}$ and maximum
$8.80\times10^{-6}$, while refining the reference from $N=512$ to $N=768$ changes
the direct field by only $1.6$--$1.9\times10^{-8}$; the dominant error is therefore
network approximation rather than the boundary-integral solver, with the fixed
$|m|\le32$ bandwidth setting a residual floor.

\begin{figure}[t]
\centering
\includegraphics[width=0.82\textwidth]{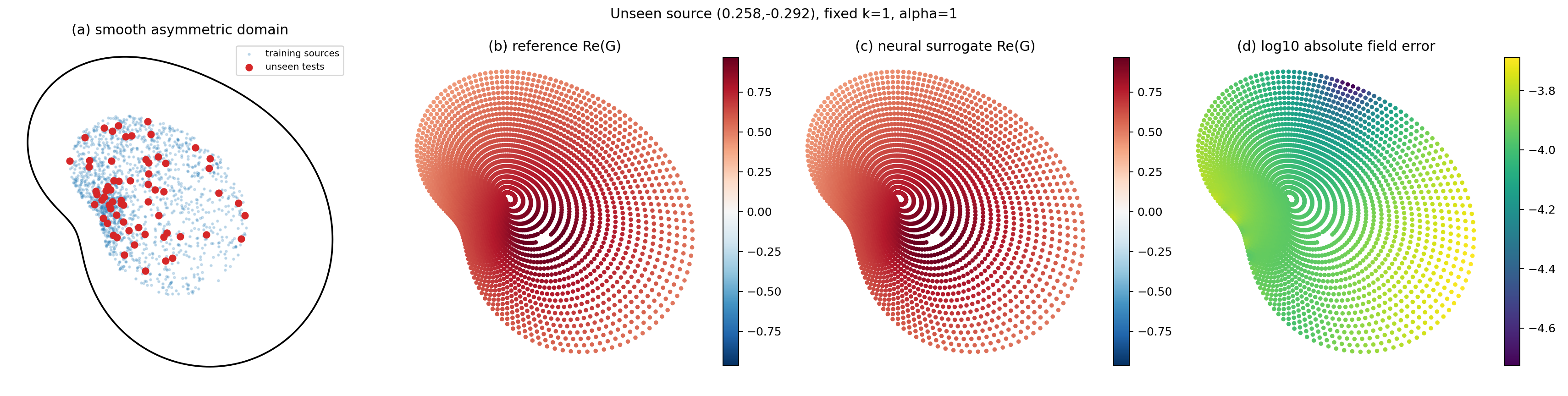}
\caption{Asymmetric smooth-boundary source-point surrogate: the network predicts
the full boundary-density Fourier coefficients for each new source and
$G=\Phi_k+\Scal\widehat\sigma_{x_0}$ is evaluated by single-layer quadrature.}
\label{fig:asym}
\end{figure}

\begin{figure}[t]
\centering
\includegraphics[width=0.7\textwidth]{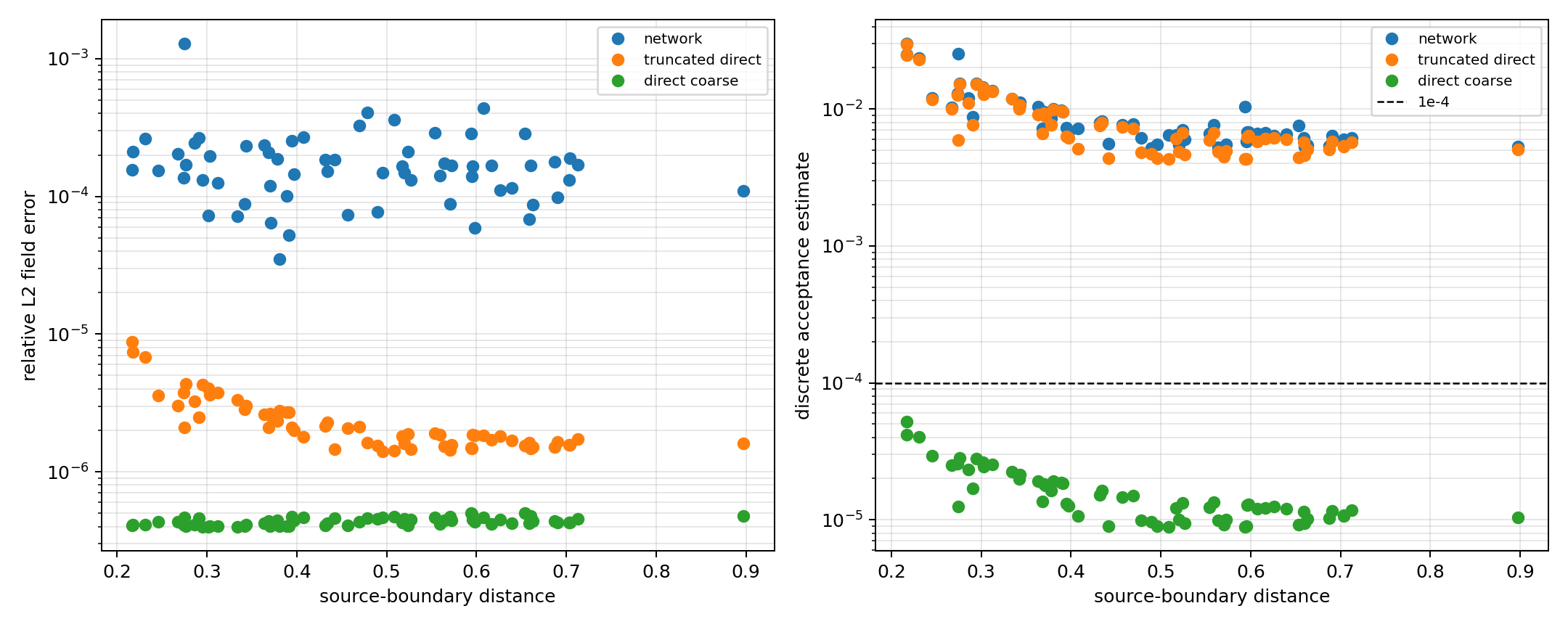}
\caption{Asymmetric smooth boundary: full-field accuracy and per-query discrete
residual--stability estimates for $64$ unseen sources.}
\label{fig:asymaccept}
\end{figure}

\subsubsection{Distance-resolved acceptance cascade on the same boundary}
\begin{figure}[t]
\centering
\includegraphics[width=0.68\textwidth]{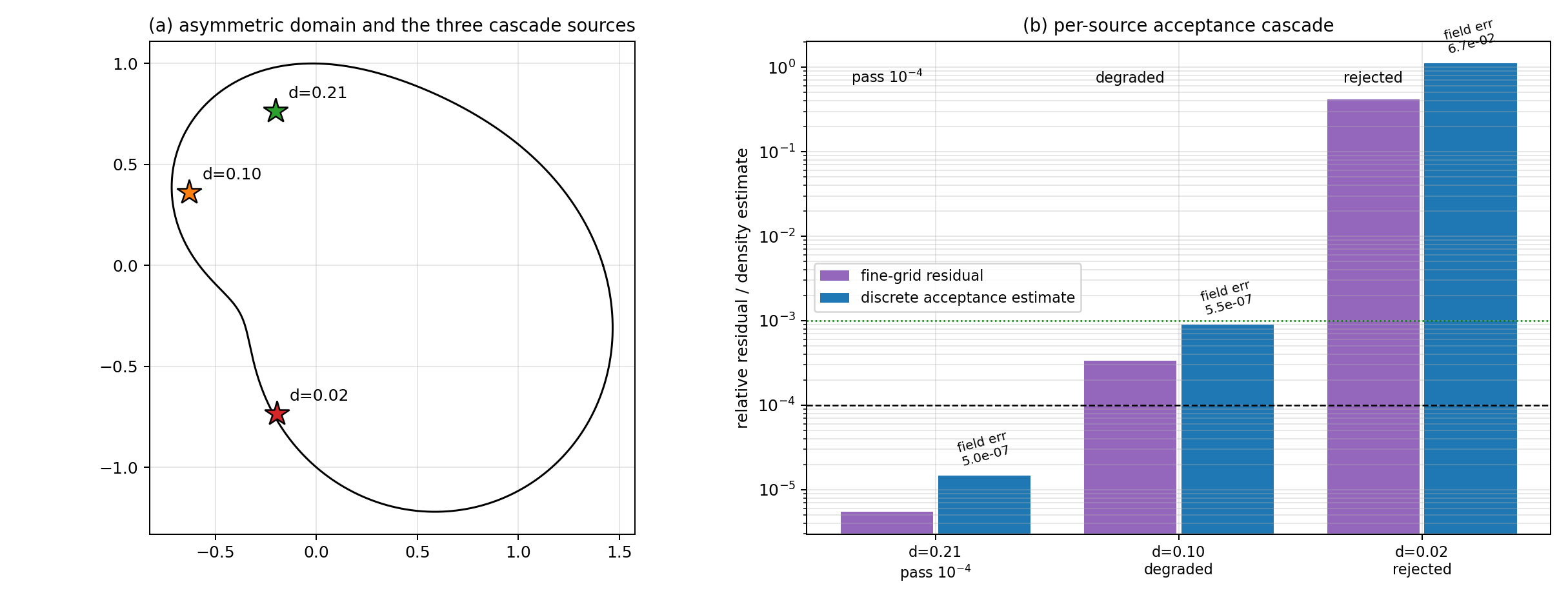}
\caption{\textbf{(a)} Asymmetric boundary and the three cascade sources (green: $d=0.21$, pass $10^{-4}$; orange: $d=0.10$, degraded; red: $d=0.02<h$, rejected). \textbf{(b)} Relative fine-grid residual, discrete density acceptance estimate, and measured field error per source; the instrument reproduces the pass/degraded/reject hierarchy of \S\ref{sec:hybrid} without disk symmetry.}
\label{fig:asymcascade}
\end{figure}

To check that the source--boundary distance hierarchy observed on the disk in \S\ref{sec:hybrid} is reproduced by the acceptance instrument on a boundary without symmetry, we add three direct-solve queries at decreasing distance from the asymmetric boundary (fallback path of Algorithm~1; the trained network is bypassed, since these sources lie outside its training region). Table~\ref{tab:asymcascade} and Figure~\ref{fig:asymcascade} summarize the results. \emph{Protocol note.} The three cascade sources were placed by hand during the writing of Table~\ref{tab:asymcascade} and regenerated at three separated convex-arc boundary points for Figure~\ref{fig:asymcascade} (the regeneration protocol is documented in the reproduction guide); the two protocols differ in detail --- fine-grid residuals $6.8\times10^{-6}$/$8.6\times10^{-5}$/$2.4\times10^{-1}$ (Table~\ref{tab:asymcascade}) versus $5.5\times10^{-6}$/$3.4\times10^{-4}$/$4.1\times10^{-1}$ (regenerated), field errors $5.0\times10^{-7}$/$5.3\times10^{-7}$/$2.1\times10^{-2}$ versus $5.0\times10^{-7}$/$5.5\times10^{-7}$/$6.7\times10^{-2}$ --- while the pass/degraded/rejected verdicts coincide in both protocols, with the regenerated $d{=}0.10$ point now in the degraded band $(10^{-4},10^{-3}]$ exactly as in Table~\ref{tab:asymcascade}.

\begin{table}[H]
\centering\small
\caption{Direct-solve acceptance cascade on the asymmetric boundary (coarse $N=192$, $h\approx0.036$; fine residual at $N=512$). The $d=0.21$ row lies in the safe regime of the $64$-source test; the $d=0.02$ density spike is unresolved at the coarse mesh.}
\label{tab:asymcascade}\vspace{0.25cm}
\begin{tabular}{@{}ccccc@{}}
\toprule
$d=\mathrm{dist}(x_0,\partial\Omega)$ & fine-grid residual & discrete density bound & meas.\ field err.\ & verdict \\
\midrule
$0.21$ & $6.8\times10^{-6}$ & $1.8\times10^{-5}$ & $5.0\times10^{-7}$ & pass $10^{-4}$ \\
$0.10$ & $8.6\times10^{-5}$ & $2.3\times10^{-4}$ & $5.3\times10^{-7}$ & above $10^{-4}$ \\
$0.02<h$ & $2.4\times10^{-1}$ & $6.5\times10^{-1}$ & $2.1\times10^{-2}$ & \textbf{rejected} \\
\bottomrule
\end{tabular}
\end{table}

\subsection{Multiply connected domains: annulus and crescent}\label{sec:annulus}
On multiply connected boundaries, the operator becomes a block matrix over boundary components. For the annulus, the inner and outer traces are $-\partial_r+\alpha_1$ and $\partial_r+\alpha_2$, with the inner curvature-sign convention of Appendix~\ref{app:nystrom}(vi). The modal benchmark of \S\ref{sec:benchmark} supplies an independent check. The outer-circle block agrees with the disk format to $3\times10^{-16}$, the cross blocks agree with direct quadrature to roundoff, and the interior normal derivative shows decreasing MMS errors from $1.7\times10^{-5}$ to $2.0\times10^{-7}$ as $N_1$ increases from $32$ to $128$. A small pilot density network attained only $1.2\times10^{-1}$ relative $L^2$ error; it is not part of the high-precision surrogate claim.

We apply the discrete acceptance estimate to direct solves on the annulus $R_1=1$, $R_2=2$ at $k=1$, $\alpha_1=1$, $\alpha_2=2$. Coarse solves use $(N_1,N_2)=(64,128)$ and fine residuals use $(128,256)$.

\begin{center}\small
\begin{tabular}{@{}cccccc@{}}
\toprule
source $x_0$ & $\norm{\Acal_h^{-1}}_2$ & fine-grid residual & discrete bound & meas.\ density err.\ & meas.\ field err.\ \\
\midrule
$(1.50,0.00)$ & $1.92$ & $1.4\times10^{-5}$ & $2.7\times10^{-5}$ & $3.9\times10^{-5}$ & $3.9\times10^{-6}$ \\
$(-1.20,0.90)$ & $1.92$ & $1.4\times10^{-5}$ & $2.7\times10^{-5}$ & $3.9\times10^{-5}$ & $3.9\times10^{-6}$ \\
$(0.00,-1.70)$ & $1.92$ & $3.7\times10^{-5}$ & $7.0\times10^{-5}$ & $1.0\times10^{-4}$ & $4.6\times10^{-6}$ \\
\bottomrule
\end{tabular}
\end{center}

The discrete estimates and measured density errors have the same order, but the measured errors can exceed the leading discrete bound. These rows therefore must not be described as uniformly enclosing the refined-solve error: residual upsampling, norm conversion, and discrete-to-continuous corrections remain. The measured field errors are smaller, at most $4.6\times10^{-6}$.

\begin{figure}[t]
\centering
\includegraphics[width=0.85\textwidth]{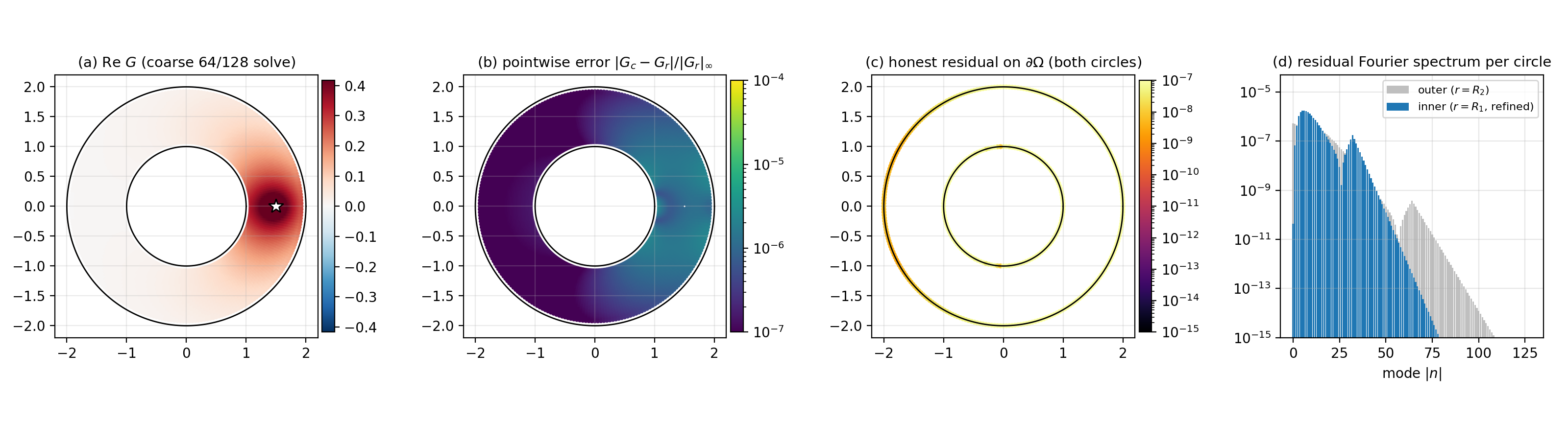}
\caption{Direct-solve diagnostics on the annulus ($R_1=1$, $R_2=2$, $\alpha_1=1$, $\alpha_2=2$, $k=1$): \textbf{(a)} $\operatorname{Re}G$ for the coarse solve; \textbf{(b)} normalized pointwise field error against the fine reference; \textbf{(c)} fine-grid residual on both boundary components; \textbf{(d)} residual Fourier spectra.}
\label{fig:annulus}
\end{figure}

\textbf{Crescent stress test.} Offsetting the inner circle to $0.7$ with $R_1=1.2$ and $R_2=2$ produces a thin sliver. MMS errors decrease from $3.6\times10^{-7}$ to $8.3\times10^{-9}$ as $N_1$ increases from $32$ to $128$. For direct-solve queries at $(-0.90,0.10)$ and $(-0.40,0.90)$, the discrete bounds are $5.1\times10^{-4}$ and $5.0\times10^{-3}$, and the measured field errors are $6.9\times10^{-6}$ and $3.1\times10^{-5}$. A source $0.06$ from the outer arc in the thin sliver gives a bound of $0.93$ and is rejected at the stated tolerances. These examples again concern direct solves, not an amortized geometry network.

\begin{figure}[t]
\centering
\includegraphics[width=0.85\textwidth]{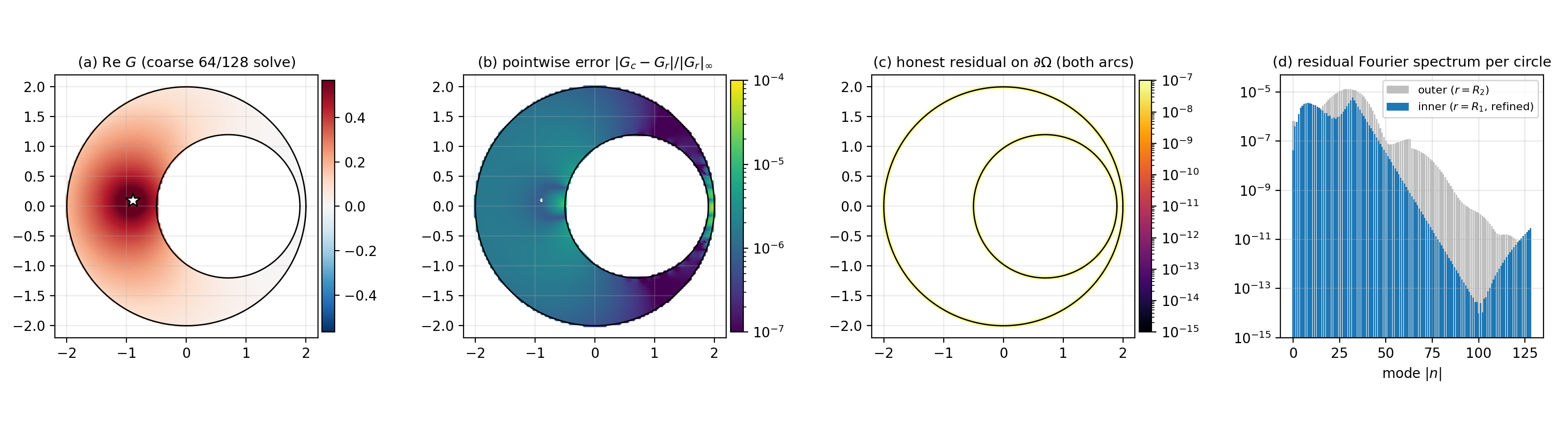}
\caption{Direct-solve diagnostics on the crescent at source $(-0.90,0.10)$: \textbf{(a)} $\operatorname{Re}G$; \textbf{(b)} normalized pointwise field error (the relative $L^2$ field error is $6.9\times10^{-6}$); \textbf{(c)} fine-grid residual on both arcs; \textbf{(d)} residual Fourier spectra.}
\label{fig:crescent}
\end{figure}

\section{Discussion}\label{sec:discussion}

\medskip\noindent\emph{Division of labor with the companion works.} The explicit kernels, spectral branches, singular-perturbation theory, and quantitative inverse-norm laws used to calibrate this instrument are established in the companion papers \cite{RobinNeumann,yang2026balls,HelmholtzAnnulus,Ainv}; the general three-leg error budget and its certified numerical sample budgets are developed separately in \cite{ErrorBudget}. The present paper contributes the single-layer representation, the discretization conventions, the per-query acceptance workflow, and the experimental validation across eight geometry classes. It specializes the continuous-source estimation mechanism to the disk in Proposition~\ref{prop:fitrate} and Theorems~\ref{thm:continuousfit}--\ref{thm:gaussianfit}, but does not import numerical budgets proved for the two-layer model of \cite{ErrorBudget}. The companions are written to be read independently; the shared acceptance narrative is kept distinct from both the stochastic a-priori theory and the deterministic query-time estimate.

\subsection{Relation to existing work}
BINet and BI-GreenNet learn boundary-integral representations through layer potentials \cite{BINet,BiGreenNet}; the present Robin correction problem uses the single-layer ansatz to avoid the hypersingular normal trace of a double-layer representation. Neural Green's Functions and neural Green's operators address complementary cross-geometry or architecture-prior questions \cite{NGF,MelchersNGO}, whereas our learned results remain at fixed geometry. The deep generalized Green's function removes the Dirac singularity through a smoothed excitation \cite{DGGF}; our representation instead subtracts the exact fundamental solution and learns the regular boundary correction.

The kernel-decomposition analysis of \cite{DuKernelDec} likewise subtracts the singular part before learning, but addresses a priori approximation rates for self-adjoint elliptic Dirichlet problems within an RKHS/optimal-recovery framework \cite{LinDuOptRecovery,LinDuOptRecHS}. Its parameter-depth degradation complements the accuracy ladder observed here. Our contribution combines query-time acceptance for point-source Robin--Helmholtz surrogates with an explicit continuous-source estimation analysis on the disk. Proposition~\ref{prop:fitrate} analyzes the identifiable-band fixed-design channel, Theorem~\ref{thm:continuousfit} supplies pointwise source-interval lifting, and Theorem~\ref{thm:gaussianfit} supplies only an integrated-source Gaussian statement. Deterministic discretization remains in the residual decomposition, and the blind-mode floor is explicit.

\medskip\noindent\emph{Trust signals versus priced acceptance.} The closest 2026 relatives of the
acceptance instrument are the calibrated physics-residual trust scores of \cite{SelfAudit}, the
conformal coverage bands for operator learning of \cite{UQNO,ConfDeepONet,PICP}, and the
ensemble-disagreement validity check with solver routing of \cite{CLAEM}. They share our central
conviction that a surrogate must decline queries it cannot support. The difference is the currency of
the verdict: a conformal band covers the truth with a stated probability \emph{over the calibration
distribution}, and an ensemble-disagreement score is informative when the failure mechanism resembles
extrapolation; neither delivers a deterministic statement about \emph{this} query's error. The
residual--stability product is deterministic and query-specific. On one fixed discretization it is
the exact worst-case discrete error bound, with the inverse norm of the discrete Robin operator as
its stability input --- Plancherel-exact on the disk and a dense inverse norm elsewhere. The
reported coarse-stability/fine-residual implementation is instead a discrete estimate; its
conservatism is itself measurable (the effectivity columns of \S\ref{sec:results}). When the failure mode is spectral rather than statistical --- a resonance skirt
(\S\ref{sec:cavityscreen}) --- calibration-set and ensemble signals have no direct purchase, while
the stability factor diverges through the exact two-sided law of Appendix~\ref{app:ainv}.

\subsection{Scope of the discrete acceptance estimate}
For an exact continuous operator and residual, error--residual equivalence is two-sided. The experiments instead use a discrete inverse norm and a fine-grid residual; hence the reported quantity is a computable discrete estimate, not an interval-enclosed continuous bound. Direct solves also require the same distinction, and measured errors can locally exceed the leading estimate.

A continuous certificate requires four additions: computable residual/interpolation and omitted-mode enclosures, a stability upper bound, and field-kernel quadrature control. On the disk, explicit modal ingredients are supplied by Lemma~\ref{lem:tail} and Proposition~\ref{prop:stability}. On general domains, the inverse-norm route of \cite{Ainv} can bound stability through certified spectral distances, while residual and tail enclosures require analyticity or singular-expansion information. None of these missing corrections is silently set to zero here.

For the exact continuous residual and stability constant, the remaining issue is effectivity rather than formal validity; for the computed coarse/fine-grid quantity, the correction terms above remain separate. The acceptance product amplifies all residual directions by the worst stability constant, although field errors depend on their directions. Direct solves on benign geometries have density effectivity close to one; learned general-domain surrogates can instead have small field errors but density bounds larger by orders of magnitude. The dual-weighted bound removes irrelevant directions for a fixed quantity of interest, but not the component of network error visible to that quantity. Accordingly, density acceptance and field-level estimates answer different tolerances and should not be compared as if they shared a normalization.

\begin{table}[H]
\centering\footnotesize
\caption{Accuracy ladder across parameter axes (one axis varied, the rest fixed; bands per \S\ref{sec:robingap}:
decision-grade $10^{-4}$, screening-grade $10^{-3}$, triage-grade $10^{-2}$).}
\label{tab:ladder}\vspace{0.25cm}
\setlength{\tabcolsep}{3.5pt}
\begin{tabular}{@{}>{\raggedright\arraybackslash}p{2.55cm}>{\raggedright\arraybackslash}p{2.55cm}>{\raggedright\arraybackslash}p{2.55cm}>{\raggedright\arraybackslash}p{1.95cm}>{\raggedright\arraybackslash}p{5.5cm}@{}}
\toprule
Axis varied & Setting & Method (\S) & Band & Reason \\
\midrule
impedance $\alpha$ & disk, $k{=}1$, sources $r_0\le0.9$ & hybrid modal $+$ Graf tail (\S\ref{sec:certified}) & decision & $\alpha$ enters as the zeroth-order term $\alpha\Scal$; explicit symbol; Plancherel-exact stability constant \\
source location $x_0$ & disk, $k{=}1$, $\alpha$ fixed & hybrid modal (\S\ref{sec:hybrid}) & decision ($r_0\le0.9$) & source enters only the right-hand side; analytic Graf tail; discrete acceptance test rejects $r_0=0.95,0.98$ \\
source location $x_0$ & general smooth boundary & all-learned network (\S\ref{sec:asym}) & screening & no analytic structure; measured field $1.65\times10^{-4}$ median, but density estimate conservative ($5.2\times10^{-3}$--$3.0\times10^{-2}$) \\
wavenumber $k$ & disk & joint $(k,\alpha)$ network (\S\ref{sec:kalpha}); direct solves (\S\ref{sec:highfreq}) & triage--screening & $k$ enters the kernel itself: oscillatory quotient dependence, conditioning degrades with $k$; direct-solve discrete acceptance estimate reported at $k{=}185$; continuously certified learned $k$-amortization not yet reported \\
geometry & fixed $k$, $\alpha$, source & no cross-geometry learned experiment & not yet assigned & only direct ellipse enumeration is reported; fixed-geometry flower/asymmetric results must not be relabelled as geometry generalization \\
\midrule
near-resonance / near-boundary & any axis & any method & rejected & spectral or spike defect; the discrete acceptance test rejects regardless of measured field accuracy (\S\ref{sec:hybrid}, \S\ref{sec:cavityscreen}) \\
\bottomrule
\end{tabular}
\end{table}

\subsection{Accuracy bands form a ladder across parameter axes}\label{sec:ladder}
The experiments of \S\ref{sec:results} admit a compact reading. When a single parameter axis is varied with the
others fixed (typically $k{=}1$, $\alpha{=}1$, interior sources), the accuracy band attained by a learned surrogate
tends to degrade with the depth at which the parameter enters the problem --- right-hand side ($x_0$),
zeroth-order operator term ($\alpha$), and kernel itself ($k$). Geometry is listed separately as an open axis, because no cross-geometry learned surrogate was run. Table~\ref{tab:ladder}
records the ladder with the supporting experiments; the bands are those of \S\ref{sec:robingap}. The ladder records today's evidence, not an impossibility statement: no continuously certified learned result is reported for the $k$- and geometry axes, and the near-resonance/near-boundary row is the discrete acceptance test's designed behaviour, not a failure of the method.

The ladder depends on an axis together with its representation prior: source amortization reaches decision grade on the disk, where the Graf tail is available, but only screening grade on the asymmetric boundary. Joint $(k,\alpha)$ amortization is harder than either axis alone. In deployment, low-dimensional geometry families can be enumerated directly, while source or impedance axes benefit more from amortization. High precision in the present experiments is bought primarily by analytic structure rather than network capacity; no cross-geometry learned accuracy band is claimed.

\subsection{Summary of numerical evidence and limitations}

The ordinary parameter-amortization experiment also exposes a structured difficulty map. On the disk, the modal quotients have low-order rational dependence on $\alpha$, whereas their $k$ dependence contains oscillatory Bessel and Hankel factors. This is consistent with the observed distinction between resonance-skirt failures, which carry conditioning signatures, and $k$-extrapolation failures, which need not do so. Cross-geometry amortization adds a still higher-dimensional shape dependence and remains untested. This ordering is a disk-based interpretation of the experiments, not a universal theorem for all geometries.

Several limitations delimit these results. High-precision learning is source-amortized at fixed $k=1$, with a separate fit per impedance; joint $(k,\alpha)$ learning is only ordinary accuracy, and no cross-geometry learned surrogate is reported. The Graf tail is disk-specific, and general-domain density bounds remain conservative. Near-boundary sources require specialized quadrature, while polygonal examples test rejection rather than validated corner surrogates. High-frequency direct solves are validated at $k{=}185$ (\S\ref{sec:highfreq}); high-frequency trapping regimes, learned high-frequency surrogates, the Laplace limit, exterior resonances, and three-dimensional surrogates remain outside the validated scope. The cavity demonstration is a two-dimensional direct-solve screen; its axisymmetric Leontovich extension remains open. Finally, continuous residual, stability, tail, and kernel-quadrature corrections must be enclosed before the discrete estimates become rigorous continuous certificates.

\section{Conclusion}\label{sec:conclusion}
We formulated Robin Green surrogates through a single-layer correction density, with the same boundary operator supplying the training residual and a per-query discrete residual--stability acceptance estimate. For the structured disk estimator, an explicit Chebyshev design factor lifts fixed-design estimation to the continuous source interval. General sub-Gaussian noise gives a pointwise source bound; proper complex Gaussian noise gives a sharper integrated-source bound, without any pointwise Gaussian overclaim. Blind modes and deterministic quadrature and tail channels remain separate from both statistical routes.

Numerically, learned low modes plus an analytic Graf tail pass the $10^{-4}$ discrete density test at all 32 sampled nonresonant impedances and three tested source radii. All-learned surrogates on general smooth boundaries attain smaller measured field errors than their conservative density estimates suggest, while near-boundary, corner-defective, and resonance-skirt queries expose genuine failure or rejection regimes. Direct-solve tests on eight geometry classes, a high-frequency ($k{=}185$) disk validation, an ellipse screening workflow, and dual-weighted output estimates demonstrate that the instrument is not restricted to the disk, to low frequencies, or to learned queries.

The contribution is thus a fixed-geometry Robin representation-and-acceptance workflow with a continuous-source estimation theory for the structured disk fit. It does not establish general method superiority, a training guarantee for deterministic direct-solve data or generic nonconvex networks, or a rigorous continuous field certificate for the reported coarse/fine-grid experiments. The remaining next steps are cross-geometry surrogates with query-time estimates, high-precision joint $(k,\alpha)$ amortization, and continuous certification beyond the present exact-modal ingredients by enclosing stability, residual/interpolation, tail, and kernel-quadrature corrections. Near-boundary and corner-aware representations remain equally important for extending the accepted query set.

\appendix

\section{Nystr\"om implementation notes}\label{app:nystrom}
The following items all come from actual implementation and item-by-item verification; we record them so that the method is reproducible.

\textbf{(i) Jump term $\tfrac12\sigma$ of the interior normal trace.} The principal-value integral gives $\Kcal'\sigma$; the interior Robin condition needs $\partial_n^{-}\Scal\sigma = \big(\tfrac12 I + \Kcal'\big)\sigma$. The index binding of the jump term $\tfrac12\sigma(x_i)$ must match the collocation scheme (coincident collocation: same index; midpoint collocation: the value of $\sigma$ at the midpoint).

\textbf{(ii) Sign convention of the kernel.} The Helmholtz normal kernel is
\begin{equation}\label{eq:normal-kernel}
\partial_n\Phi_k = -\frac{\ii k}{4}H_1^{(1)}(kd)\,\frac{(x-z)\cdot n}{d},
\end{equation}
the minus sign coming from $H_0^{(1)\prime} = -H_1^{(1)}$. Crucially, the operator matrix and the right-hand side must carry the \emph{same} sign: flipping both consistently yields a system that ``solves the wrong problem self-consistently'' (tiny least-squares residual, but the Robin problem for $-\partial_n+\alpha$).

\textbf{(iii) Diagonal correction.} Under coincident collocation, the $\partial_n\Phi$ kernel is bounded with diagonal limit $-\kappa/(4\pi)$ in 2D (and $-\tfrac{1}{4\pi R}$ on a sphere of radius $R$ in 3D); the $\Phi$ kernel is logarithmically singular, and for equispaced nodes its effective diagonal is $d_{\mathrm{eff}} = R/N$ (derived exactly from $\prod_{k=1}^{N-1}2\sin\tfrac{\pi k}{N}=N$: the off-diagonal logarithmic sums give $w[(N{-}1)\ln R + \ln N]$, and appending $w\ln(R/N)$ makes the discrete single-layer potential reproduce the logarithmic ($k=0$) identity $\mathcal{V}[1]=2\pi R\ln R$ exactly). An early implementation used $2\pi R/N$, introducing an $O(w)$ diagonal error; after correction, the direct solve at $k{=}1$, $N{=}64$ improves from $2\times10^{-5}$ to $3.4\times10^{-6}$.

\textbf{(iv) The 2D Laplace logarithmic-capacity obstruction.} On the unit disk, the constant density lies in the nullspace of the pure single-layer operator $\Scal$ at $k=0$ (since $\Scal[1]=R\ln R=0$ for $R=1$), so the single-layer representation can be singular or nonunique even when the underlying Robin boundary-value problem is well posed. This is a representation-specific obstruction, not a general assertion that every nonzero-mean right-hand side is incompatible. Solvability must be checked through the adjoint nullspace. A zero-mean or other gauge constraint can select a density only after compatibility is established; by itself it does not restore existence. The numerical experiments avoid this obstruction by using $k=1$; modified logarithmic kernels or constrained formulations are alternatives for Laplace problems.

\textbf{(v) Index binding of quadrature weights (general domains).} Under a nonuniform parameterization, weights must be bound to the \emph{integration} variable (columns): $A=(\mathcal{K}'+\alpha\mathcal{S})\operatorname{diag}(w)$. Binding them to rows instead is \emph{completely invisible on the disk} (weights constant there: the disk limit passes at machine precision and the Green's verification against the series also passes), while producing an $O(\varepsilon)$ systematic error on general domains that does not blow up by itself. A manufactured-solution (MMS) check pins it down: error stuck at $1.5\times10^{-2}$ with no convergence in $N$. The lesson, of the same type as (ii): \textbf{highly symmetric benchmarks can hide entire classes of bugs}; general geometry is the touchstone. The general-domain validations are gated by MMS or direct cross-checks.

\textbf{(vi) Curvature sign on multiply connected domains.} For inner boundary components, the diagonal limit of $\Kcal'$ is $+\kappa/(4\pi)$ (the domain's outward normal points toward the center of curvature), opposite in sign to the outer boundary's $-\kappa/(4\pi)$. Reusing the outer formula on inner components creates an $O(1)$ systematic error that is invisible in the Dirichlet limit ($\alpha$ large, $\Kcal'$ swamped by $\alpha\Scal$). Located by annulus-series cross-check and fixed ($4.9\times10^{-3}\to1.4\times10^{-5}$; verified against the annulus-series benchmark of \S\ref{sec:benchmark}).

\section{The analytic disk modal quotient}\label{app:tail}
On the unit disk the Graf addition theorem gives the exact modal quotient in closed form. Define
\[
r_m := \frac{H_{m+1}^{(1)}(kR)}{H_m^{(1)}(kR)},
\qquad
\mathcal{B}_\alpha[J_n](kR) = (\alpha + n/R)J_n(kR) - kJ_{n+1}(kR).
\]
We compute $r_m$ recursively from $r_0=H_1^{(1)}/H_0^{(1)}$ using $r_m=2m/(kR)-1/r_{m-1}$.
\[
q_n(\alpha) = -\frac{1}{2\pi R}\,\frac{(\alpha + n/R) - k\,\widehat{r}_n}{\mathcal{B}_\alpha[J_n](kR)},
\]
\[
c_n(r_0) = q_n(\alpha)\,J_n(kr_0).
\]
The upward recurrence is anchored against direct Hankel evaluations through $|n|=142$ and against its large-order asymptotic at $|n|=220$ in the released unit tests. Here $\widehat{r}_n = r_n$ for $n\ge0$ and $\widehat{r}_n = -1/r_{|n|-1}$ for $n<0$ (from $H_{-m}^{(1)} = (-1)^m H_m^{(1)}$). The appended tail $49\le|n|\le220$ used in \S\ref{sec:hybrid} is $c_n$ with $q_n$ from this formula; the quotient form was numerically stable throughout the tested mode range, whereas direct evaluation of the separate factors could underflow or overflow.

\medskip\noindent\textbf{Finite-mode computational status.} The stability factor is evaluated from $\lambda_n$ in Proposition~\ref{prop:stability} for $|n|\le80$ and the analytic tail is terminated at $|n|=220$; the omitted-mode sums are directly evaluable from the closed-form coefficients \eqref{eq:normalized-coeff} (evaluated values in \S\ref{sec:accept}), so no enclosure constants are needed for the disk tail. These formulas do not, however, supply a roundoff bound for the quotient recurrence, numerical values for all constants in Theorem~\ref{thm:conv}, or a modal-tail enclosure on general domains: accordingly, ``computable'' here means directly evaluable in the stated finite discretization, with the disk modal tail bounded and general-domain omitted modes unenclosed.

\section{Exact inverse-norm law on the disk}\label{app:ainv}
The calibration statements of \S\ref{sec:resonance} rely on disk facts only; we prove them here self-contained, so that the quantitative use made of them does not depend on unpublished material. Generalizations to arbitrary domains and to multiply connected blocks are the subject of the companion theory \cite{Ainv} and are not needed by the acceptance instrument.

\begin{proposition}[disk law]\label{prop:ainv}
Let $k>0$, $R>0$, and let $\alpha>0$ be bounded away from the resonance curves \eqref{eq:alpha-res}, with $\mathrm{dist}(\alpha,\Sigma(k))=\delta$ attained at a unique nearest mode $n^{*}$. Then:
\begin{enumerate}[leftmargin=1.6em,itemsep=1pt]
\item \textbf{exact pole formula.} $\Delta_n(\cdot,k)$ is affine with root $\alpha_n(k)$ and slope $J_n(kR)$, hence
\[
\big|\lambda_{n^{*}}(\alpha,k)\big|=\frac{\pi R}{2}\,\big|H_{n^{*}}^{(1)}(kR)\,J_{n^{*}}(kR)\big|\cdot\delta;
\]
\item \textbf{skirt exponent exactly $1$.} as $\alpha\to\alpha_{n^{*}}(k)$ at fixed $k$,
\[
\big\|\Acal_{k,\alpha}^{-1}\big\|\cdot\delta \;\longrightarrow\; \frac{2}{\pi R\,\big|H_{n^{*}}^{(1)}(kR)\,J_{n^{*}}(kR)\big|};
\]
\item \textbf{safe zone.} if $k$ is not a Neumann resonance and $0<\alpha\le\alpha_0(k,R)$ lies below the first resonance curve, then $\|\Acal_{k,\alpha}^{-1}\|\le C(k,R,\alpha_0)$;
\item \textbf{computable two-sided law.} for $\alpha$ in any closed bounded window $W$: the lower bound
\[
\frac{c(k,R)}{\delta}\;\le\;\big\|\Acal_{k,\alpha}^{-1}\big\|,\qquad c(k,R)=\frac{2}{\pi R\,\sup_n\big|H_n^{(1)}(kR)J_n(kR)\big|},
\]
holds for every $\alpha\in W$ (no regime restriction: $\|\Acal^{-1}\|=1/\min_n|\lambda_n|\ge1/|\lambda_{n^*}|$ and $|H_{n^*}^{(1)}J_{n^*}|\le\sup_n|H_n^{(1)}J_n|$); and there is a computable \emph{window gap} $d_0(k,R,W)>0$ --- the minimum of the distance from $W$ to the nearest resonance curve $\alpha_n(k)\notin W$ and of $\frac{\pi R}{2}m(k,W)\,\delta_{\mathrm{tail}}(k,W)$, where $m(k,W)$ is the minimum of $|H_n^{(1)}(kR)J_n(kR)|$ over the finite active set $\{n:\alpha_n(k)\in W\}$ and $\delta_{\mathrm{tail}}(k,W)>0$ is a uniform lower bound of $|\lambda_n|$ over the remaining modes --- such that for every $\alpha\in W$ with $\delta<d_0(k,R,W)$,
\[
\big\|\Acal_{k,\alpha}^{-1}\big\|\;\le\;\frac{C(k,R,W)}{\delta},\qquad C(k,R,W)=\frac{2}{\pi R\,m(k,W)}.
\]
Outside this regime, $\|\Acal_{k,\alpha}^{-1}\|\le\max\{C(k,R,W)/d_0(k,R,W),\,4\}$ holds trivially, so the resonance spikes are the only unbounded regime.
\end{enumerate}
All constants are computable from the explicit symbol: $c$ through $\sup_n|H_n^{(1)}(kR)J_n(kR)|$, and $C$ through the finite active set together with the uniform tail bound below.
\end{proposition}
\begin{proof}
Proposition~\ref{prop:stability} gives $\|\Acal^{-1}\|=1/\min_n|\lambda_n|$ with $\lambda_n=\frac{\ii\pi R}{2}H_n\Delta_n$, and $\Delta_n(\alpha)=(\alpha-\alpha_n(k))J_n(kR)$ by affinity, whence $|\lambda_n|=\frac{\pi R}{2}|H_nJ_n|\,|\alpha-\alpha_n(k)|$; (i) follows. For (ii), note first that $J_{n^{*}}(kR)\neq0$ on $\Sigma(k)$ by definition of the curves. Every other mode either has $\alpha_n(k)=\alpha_{n^{*}}(k)$ --- excluded by the uniqueness hypothesis, the general case adding the corresponding pole terms --- or satisfies $\liminf_{\alpha\to\alpha_{n^{*}}(k)}|\alpha-\alpha_n(k)|\ge|\alpha_{n^{*}}(k)-\alpha_n(k)|>0$, with $|H_nJ_n|>0$ ($H_n^{(1)}$ has no real zeros, and $J_n(kR)=0$ would make $\Delta_n$ constant and nonzero, so $\lambda_n\neq0$). Hence $\min_{n\neq n^{*}}|\lambda_n|$ stays bounded away from zero and the minimum is attained at $n^{*}$ for small $\delta$, proving (ii). For (iii) and the tail estimate underlying (iv), use that $\lambda_n\to\tfrac12$ in the precise sense of Lemma~\ref{lem:tail}: there is a computable $n_2$ with $|\lambda_n(\alpha)|\ge\tfrac14$ for all $|n|\ge n_2$ and all $\alpha$ in the window (for $|n|\ge n_2$, $\alpha_n(k)\le-2\alpha_{\max}$ and $|H_nJ_n|\,|\alpha-\alpha_n(k)|\to\tfrac12$ uniformly by the large-order asymptotics $|H_n(kR)J_n(kR)|=\frac{2}{\pi|n|}(1+o(1))$ and $|\alpha_n(k)|=\frac{|n|}{R}(1+o(1))$). The remaining $|n|<n_2$ modes are a finite set: those whose $\alpha_n(k)$ meets the window (the active set) contribute the factor $\delta$ through $\min_{\text{active}}|H_nJ_n|=m(k,W)$, giving the upper bound $2/(\pi R\,m(k,W)\,\delta)$, and the rest are bounded below by a positive constant on the compact window; the tail modes $|n|\ge n_2$ satisfy $|\lambda_n|\ge\tfrac14$. The lower bound in (iv) holds at every $\alpha\in W$ because $\|\Acal^{-1}\|\ge1/|\lambda_{n^*}|=2/(\pi R\,|H_{n^*}J_{n^*}|\,\delta)\ge c(k,R)/\delta$ with $c(k,R)=2/(\pi R\sup_n|H_n^{(1)}J_n|)$. The upper bound requires the window-gap regime $\delta<d_0(k,R,W)$: then the nearest mode is active (every inactive mode satisfies $|\alpha-\alpha_n(k)|\ge\mathrm{dist}(\alpha_n(k),W)\ge d_0>\delta$), the active modes are controlled by $m(k,W)$ as above, and the choice $d_0\le\tfrac{\pi R}{2}m(k,W)\delta_{\mathrm{tail}}$ upgrades the tail floor $|\lambda_n|\ge\delta_{\mathrm{tail}}$ to $\delta_{\mathrm{tail}}^{-1}\le 2/(\pi R\,m(k,W)\,\delta)$. Combining the three groups yields (iii) (no resonance below $\alpha_0$) and (iv) with the stated computable constants.\end{proof}

\section{Reproducibility map}\label{app:repro}
\begingroup\sloppy\raggedright
The accompanying single-file script uses a positional command, for example \texttt{python SL-GreenNet.py disk}. The \texttt{anchors} mode checks the Wronskian, Hankel-ratio recurrence, Fourier symbols, closed-form disk density, and field normalization. The \texttt{disk}, \texttt{quad-refine}, \texttt{stab-tail}, and \texttt{ellipse} modes regenerate the disk acceptance sweep, the quadrature-consistency anchors of Table~\ref{tab:quadrefine}, the finite-mode stability check, and the ellipse screening data and figures, respectively. The \texttt{annulus}, \texttt{crescent}, \texttt{flower}, \texttt{square}, and \texttt{lshape} modes cover the corresponding direct-solve tests. The $k{=}185$ spot checks of \S\ref{sec:highfreq} are regenerated by the \texttt{highfreq} mode (disk and flower acceptance-estimate--$N$ curves) and the \texttt{highfreq-bean} mode (dense-resonance refusal; $\sim$5 minutes each, single core), and the L-shape diagnosis behind the \S\ref{sec:polygon} paragraph by the \texttt{corner-audit} mode ($\sim$2 minutes); all three are pure NumPy/SciPy. Neural experiments are separated into \texttt{train-resnet}, \texttt{disk-kalpha}, \texttt{flower-alpha}, and \texttt{asym}; \texttt{dual-demo} requires the checkpoint and test metrics produced by \texttt{asym}. Additional research modes not tied to a reported table (\texttt{benchmark}, \texttt{fem-verify}, \texttt{kidney}, \texttt{geo-star}, \texttt{alpha-sweep}, \texttt{spike}, \texttt{spike-e6}, \texttt{spike-e7}, \texttt{spike-e8}, \texttt{axisym-mms}, \texttt{axisym-eig}) are documented in the script header. Thus the file is a common implementation and reproduction driver, not a claim that one command regenerates every table without intermediate artifacts. The pole-model A/B comparison and the manually assembled cascade summary are documented exceptions to complete table-level automation. Random seeds, discretizations, and training budgets are stated in the corresponding sections and script defaults; timings remain implementation- and hardware-specific.\par\endgroup

\begingroup\sloppy\raggedright

\par\endgroup
\end{document}